\documentclass[a4paper,english,10pt,twoside,reqno]{amsart}

\usepackage{amssymb}
\usepackage{amsmath,amsthm}
\usepackage[english]{babel}
\usepackage[utf8]{inputenc}
\usepackage[T1]{fontenc}
\usepackage{lmodern}
\usepackage{stmaryrd}
\usepackage[hyperindex,colorlinks, allcolors=black]{hyperref}
\usepackage{mathtools}
\DeclareMathAlphabet{\pazocal}{OMS}{zplm}{m}{n}
\usepackage{color}
\usepackage{xcolor}
\usepackage{pgfplots}
\usepackage{subfig}
\usepackage{ulem}

\newtheorem{theorem}{Theorem}[section]
\newtheorem{corollary}[theorem]{Corollary}

\newtheorem{lemma}[theorem]{Lemma}
\newtheorem{proposition}[theorem]{Proposition}

\theoremstyle{definition}

\newtheorem{remark}[theorem]{Remark}

\newcommand{\eps}{\varepsilon}
\newcommand{\norm}[1]{\left\Vert #1 \right\Vert}
\newcommand{\sub}{\subseteq}
\newcommand{\NN}{\mathbb{N}}

\newcommand{\RR}{\mathbb{R}}
\newcommand{\CC}{\mathbb{C}}
\newcommand{\cA}{{\mathcal{A}}}
\newcommand{\cB}{{\mathcal{B}}}
\newcommand{\cC}{{\mathcal{C}}}
\newcommand{\cD}{{\mathcal{D}}}
\newcommand{\cE}{{\mathcal{E}}}
\newcommand{\cF}{{\mathcal{F}}}

\newcommand{\cL}{{\pazocal{L}}}

\newcommand{\cM}{\mathcal{M}}
\newcommand{\cN}{{\mathcal{N}}}
\newcommand{\cO}{{\mathcal{O}}}
\newcommand{\cS}{{\mathcal{S}}}
\newcommand{\cT}{\mathcal{T}}

\newcommand{\cZ}{\mathcal{Z}}

\newcommand{\spanv}{\mathrm{span}}
\newcommand{\dd}{\,\mathrm{d}}
\newcommand{\e}{\mathrm{e}}
\newcommand{\ii}{\mathrm{i}}

\newcommand{\curl}{\operatorname{curl} }
\newcommand{\divv}{\operatorname{div} }
\newcommand{\mE}{\textbf{\textrm{E}}}
\newcommand{\mH}{\textbf{\textrm{H}}}
\newcommand{\mJ}{\textbf{\textrm{J}}}
\newcommand{\tr}{\operatorname{tr}}

\newcommand{\Fext}{\mathcal{F}_{\mathrm{ext}}}
\newcommand{\Fint}{\mathcal{F}_{\mathrm{int}}}
\newcommand{\tFint}{\tilde{\mathcal{F}}_{\mathrm{int}}}
\newcommand{\Finteff}{\mathcal{F}_{\mathrm{int}}^{\mathrm{eff}}}

\newcommand{\Fintj}[1]{\mathcal{F}_{\mathrm{int},#1}}

\newcommand{\supp}{\operatorname{supp}}
\newcommand{\skappa}{\overline{\kappa}}

\newcommand{\inn}{\mathrm{in}}

\newcommand{\sskappa}{\mathring{\kappa}}

\allowdisplaybreaks

\numberwithin{equation}{section} 

\begin{document}

\title[Analysis of a Peaceman-Rachford scheme in heterogeneous media]{Analysis of a Peaceman-Rachford ADI scheme for Maxwell equations in heterogeneous media}

\thanks{Funded by the Deutsche Forschungsgemeinschaft (DFG, German Research Foundation) -- Project-ID 258734477 -- SFB 1173.}

\author{Konstantin Zerulla}\address{Department of Mathematics, Karlsruhe Institute of Technology, Englerstr. 2, 76131 Karlsruhe, Germany.\newline E-mail: \texttt{konstantin.zerulla@email.de}}

\author{Tobias Jahnke}\address{Department of Mathematics, Karlsruhe Institute of Technology, Englerstr. 2, 76131 Karlsruhe, Germany.\newline E-Mail:\texttt{tobias.jahnke@kit.edu}}
\date{November 8, 2022}

\begin{abstract}  
The Peaceman-Rachford alternating direction implicit (ADI) \linebreak scheme 
for linear time-dependent Maxwell equations is analyzed on a heterogeneous cuboid.
Due to discontinuities of the material parameters, the solution of the Maxwell equations is less than $H^2$-regular in space.
For the ADI scheme, we prove a rigorous time-discrete error bound with a convergence rate that is half an order lower than the classical one.
Our statement imposes only assumptions on the initial data and the material parameters, but not on the solution.
To establish this result, we analyze the regularity of the Maxwell equations in detail in an appropriate functional analytical framework.
The theoretical findings are complemented by a numerical experiment indicating that the proven convergence rate is indeed observable and optimal.
\end{abstract}
\subjclass[2020]{35Q61, 47D06, 65M15, 65J08}
\keywords{Maxwell equations, splitting method, heterogeneous media, error analysis, loss of convergence order}
\maketitle

\section{Introduction}\label{section1}

The propagation of electromagnetic waves in media can be modelled by time-dependent Maxwell equations, see \cite{Jackson,BornWolf,Griffiths,DautrayLions}.
A thorough analytical understanding, as well as an efficient and reliable numerical solution of Maxwell equations is hence desirable for many applications, such as the design of antennas and wave\-guides, see \cite{Namiki2000} and Section~9.3 in \cite{SalehTeich} for instance.

On domains with tensor structure, alternating direction implicit (ADI) schemes 
as proposed in \cite{Namiki2000,ZhCZ2000} 
are very attractive for the numerical solution of Maxwell equations.
Instead of approximating the solution of the full Maxwell system at once,
the Maxwell differential operator is split up according to the space dimensions along which derivatives are taken.
The sub-systems associated to these parts are propagated in a certain way and in a certain order in every time-step.
Hereby one alternates between explicit and implicit time integration schemes.
The implicit steps for both sub-systems only amount to the solution of essentially 
one-dimensional elliptic problems which makes the schemes very efficient.
While the original works \cite{Namiki2000,ZhCZ2000} apply a Peaceman-Rachford time integrator to the split problem, an energy conserving scheme is constructed in \cite{ChenLiLiang}.
Another attractive feature of both approaches is the numerical unconditional stability (without CFL restriction on the time step size).
An even more efficient formulation of ADI schemes is derived in \cite{Leong,Tan2}.
A modified ADI scheme that preserves the uniform exponential decay properties of damped Maxwell equations is constructed and analyzed in \cite{Zerulla20}.

In presence of material parameters that are at least Lipschitz continuous, the time discretization errors of the ADI schemes from \cite{Namiki2000,ZhCZ2000,ChenLiLiang} are rigorously analyzed in \cite{HoJaSc,Eilinghoff Diss,EiSc18,EiSc32,EiJaSc}.
By rigorous we mean that the analysis imposes only verifiable assumptions on the data but not on the unknown solution.
The main achievement of this paper is a similar error analysis in a technically much more involved situation.

The above mentioned papers analyze only the semi-discretization in time, and finite differences in space are used in numerical examples to obtain fully discrete systems.
In \cite{HoKoe,KoeDiss}, however, the Peaceman-Rachford ADI scheme is combined with a discontinuous Galerkin (dG) space discretization.
The implicit steps are here shown to decouple into block diagonal systems where the size of the blocks depends only on the polynomial degree of the dG ansatz space.
In \cite{KoeDiss,HoKoe20} the error of the dG-ADI full discretization is additionally analyzed, establishing the classical order of the schemes under assumptions on the data and the solution.

In heterogeneous media like waveguides, the material parameters in Maxwell equations are discontinuous.
This leads to new difficulties for the analysis of ADI schemes.
In \cite{ZerullaDiss}, the abstract time-discrete Peaceman-Rachford ADI scheme from \cite{Namiki2000,ZhCZ2000} is shown to converge with reduced order $3/2$ in $L^2$ on a cuboid that consists of two homogeneous subcuboids with a common interface.
In \cite{Zeru22b}, a more complicated heterogeneous partition of a cuboid is considered.
The initial data for the Maxwell equations are less regular than in \cite{HoJaSc,Eilinghoff Diss,EiSc18,EiSc32,EiJaSc,ZerullaDiss}, such that the previous error analysis does not apply.
For this reason, a different dimension splitting scheme is constructed.
In the current paper, we study a similar material configuration as in \cite{Zeru22b}, but more regular data giving rise to more regular solutions of the Maxwell equations and higher convergence rates for the classical Peaceman-Rachford ADI scheme.
The presented rigorous error results are new, to the best of our knowledge.

We consider the time-dependent linear isotropic Maxwell equations
\begin{align}
	\partial_t\mE&=\frac{1}{\varepsilon}\curl\mH-\frac{1}{\varepsilon}\mJ, &&\partial_t\mH=-\frac{1}{\mu}\curl\mE,\label{Maxwell system}\\
	\mE(0)&=\mE_0, &&\mH(0)=\mH_0,
\end{align}
for $t\geq0$ on a cuboid 
\begin{align*}
	Q=(a_1^-,a_1^+)\times(a_2^-,a_2^+)\times(a_3^-,a_3^+),
\end{align*}
with perfectly conducting boundary conditions
\begin{align*}
	\mE\times\nu=0,\qquad \mu\mH\cdot\nu=0
\end{align*}
on the boundary $\partial Q$.
Here, $\mE=\mE(t,x)\in\RR^3$ is the electric field, $\mH=\mH(t,x)\in\RR^3$ the magnetic field, $\mJ=\mJ(t,x)\in\RR^3$ the external volume current, $\mu=\mu(x)>0$ the magnetic permeability, and $\eps=\eps(x)>0$ the electric permittivity.
We accompany the Maxwell equations (\ref{Maxwell system}) with additional divergence conditions for $\mE$ and the Gauss law $\divv(\mu\mH)=0$, by considering the Maxwell equations on an appropriate state space $X_2$, see (\ref{def X2}).

\begin{figure}%
	\centering
	\subfloat[\centering Three dimensional picture of model domain]{\def\svgwidth{5.2cm}
		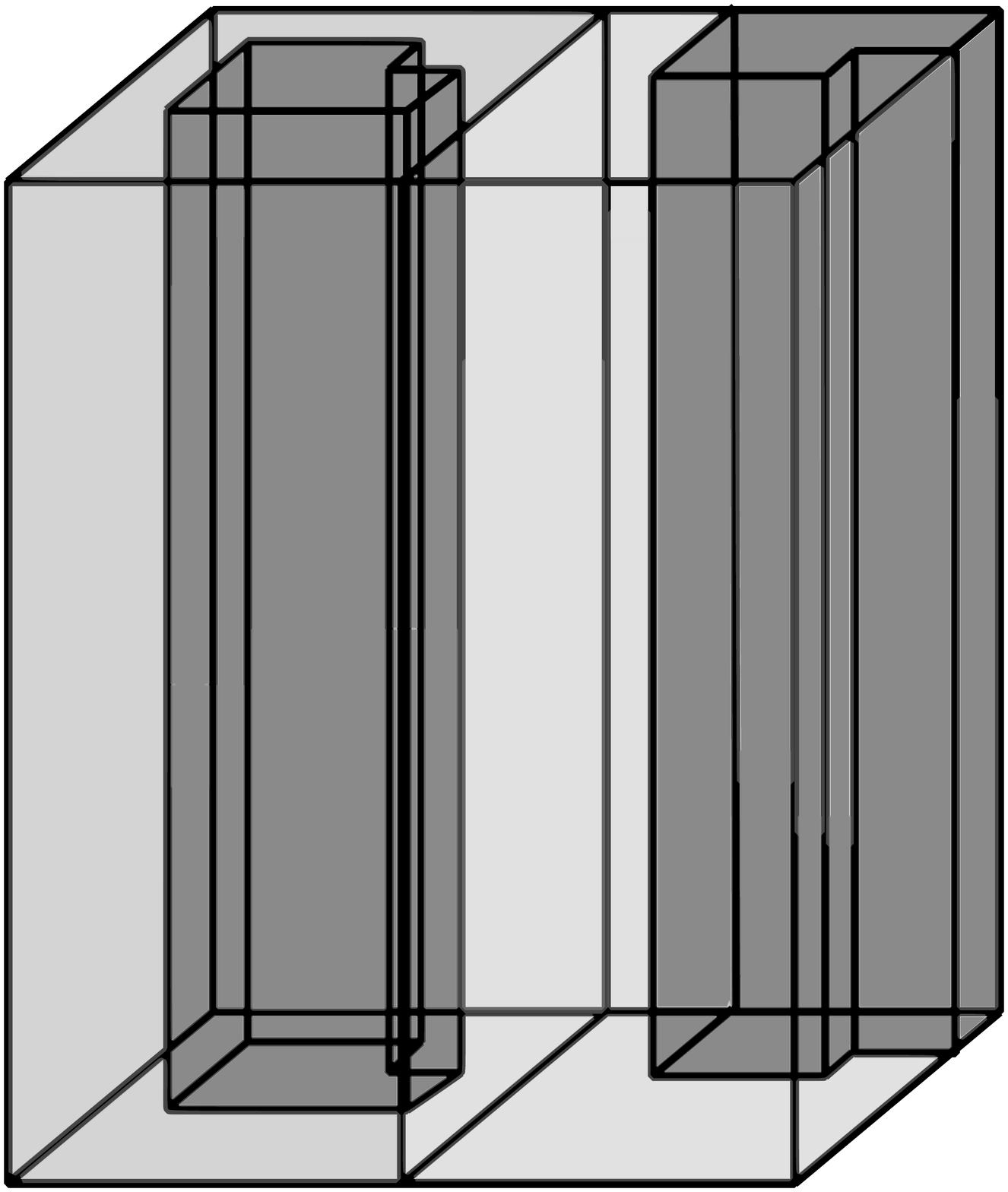 }%
	\qquad\qquad
	\subfloat[\centering Illustration of subdomain notation]{\def\svgwidth{5.2cm}
\begingroup%
  \makeatletter%
  \providecommand\color[2][]{%
    \errmessage{(Inkscape) Color is used for the text in Inkscape, but the package 'color.sty' is not loaded}%
    \renewcommand\color[2][]{}%
  }%
  \providecommand\transparent[1]{%
    \errmessage{(Inkscape) Transparency is used (non-zero) for the text in Inkscape, but the package 'transparent.sty' is not loaded}%
    \renewcommand\transparent[1]{}%
  }%
  \providecommand\rotatebox[2]{#2}%
  \newcommand*\fsize{\dimexpr\f@size pt\relax}%
  \newcommand*\lineheight[1]{\fontsize{\fsize}{#1\fsize}\selectfont}%
  \ifx\svgwidth\undefined%
    \setlength{\unitlength}{595.27559055bp}%
    \ifx\svgscale\undefined%
      \relax%
    \else%
      \setlength{\unitlength}{\unitlength * \real{\svgscale}}%
    \fi%
  \else%
    \setlength{\unitlength}{\svgwidth}%
  \fi%
  \global\let\svgwidth\undefined%
  \global\let\svgscale\undefined%
  \makeatother%
  \begin{picture}(1,1.41428571)%
    \lineheight{1}%
    \setlength\tabcolsep{0pt}%
    \put(0,0){\includegraphics[width=\unitlength,page=1]{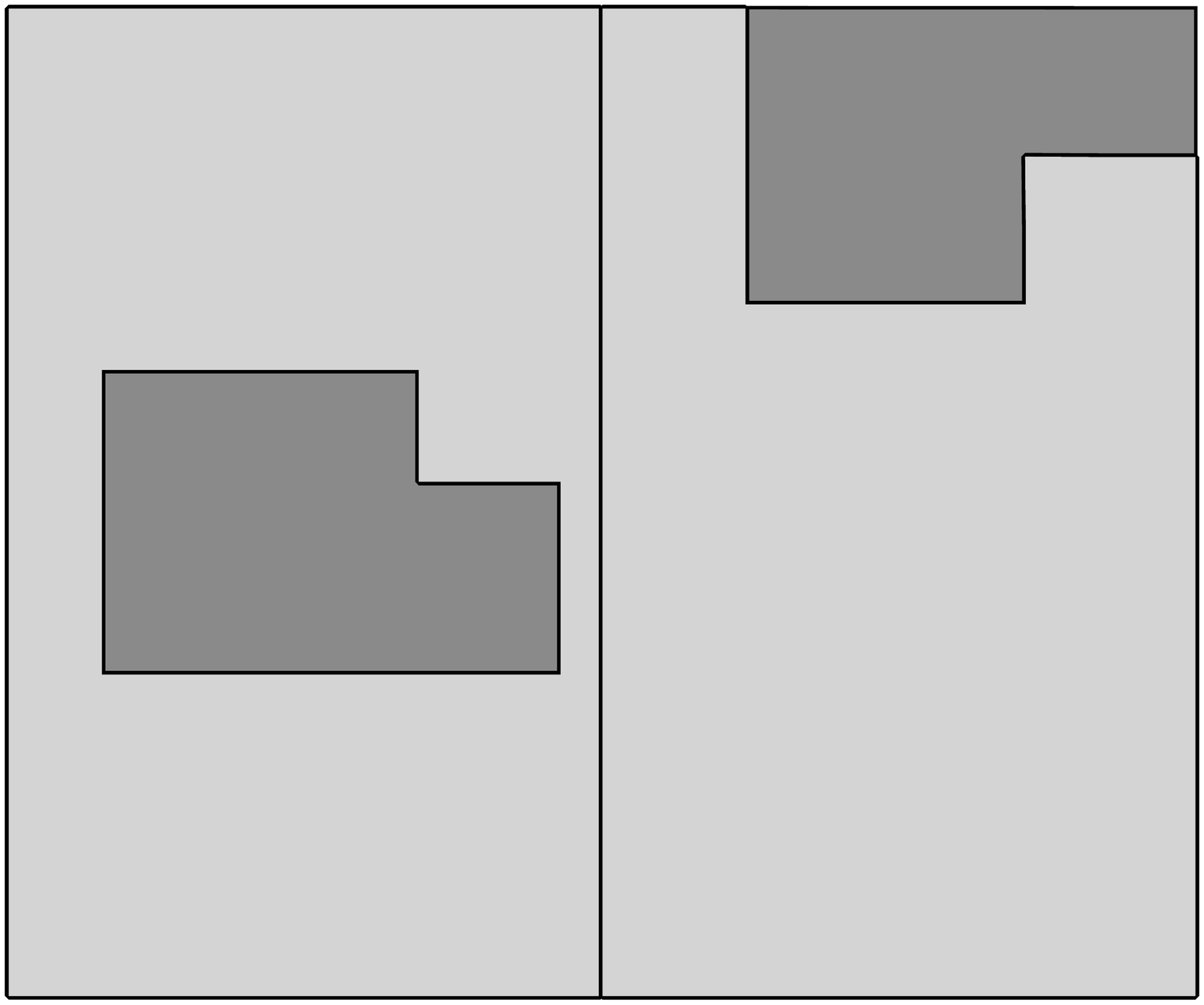}}%
    \put(0.20490665,0.73083939){\color[rgb]{0,0,0}\makebox(0,0)[lt]{\lineheight{1.25}\smash{\begin{tabular}[t]{l}$\tilde{Q}_{1,1}$\end{tabular}}}}%
    \put(0.68980989,1.0833806){\color[rgb]{0,0,0}\makebox(0,0)[lt]{\lineheight{1.25}\smash{\begin{tabular}[t]{l}$\tilde{Q}_{2,1}$\end{tabular}}}}%
    \put(0.71271871,0.82374733){\color[rgb]{0,0,0}\makebox(0,0)[lt]{\lineheight{1.25}\smash{\begin{tabular}[t]{l}$\tilde{Q}_{2,0}$\end{tabular}}}}%
    \put(0.28381478,1.05538094){\color[rgb]{0,0,0}\makebox(0,0)[lt]{\lineheight{1.25}\smash{\begin{tabular}[t]{l}$\tilde{Q}_{1,0}$\end{tabular}}}}%
  \end{picture}%
\endgroup%
}%
	\caption{Example of a heterogeneous model domain $Q$}%
	\label{fig:example}%
\end{figure}

The parameters $\eps$ and $\mu$ describe the material $Q$ consists of.
The geometric conditions for the composition of $Q$ are similar as in the companion paper \cite{Zeru22b}, and they originate from the model of an embedded waveguide, see Section~9.3 in \cite{SalehTeich} for instance:
We divide $Q$ into a chain of smaller cuboids $\tilde{Q}_1,\dots,\tilde{Q}_L$ with interfaces parallel to the $x_2-x_3$-plane.
We collect these interfaces in a set $\tilde{\cF}_{\text{int}}$.
In each cuboid $\tilde{Q}_i$, we additionally have smaller separated subdomains  that touch the faces $\{x_3=a_3^-\}$ and $\{x_3=a_3^+\}$, are distinct from the interfaces in $\tFint$, and satisfy the following property. 
Each subdomain can be represented by a cuboidal grid with all grid elements touching the top and bottom face of $Q$, and adjacent grid elements having a common interface.
The remainder of $\tilde{Q}_i$ is denoted by $\tilde{Q}_{i,0}$.
Note that the subdomains $\tilde{Q}_{i,1},\dots,\tilde{Q}_{i,K}$ represent embedded waveguide structures.
Figure~\ref{fig:example} displays an example for the considered setting with $L=2$ and $K=1$.
The statements and arguments transfer in a straightforward manner to more complicated structures, where a cuboid $\tilde{Q}_{i,j}$, $j\in\{1,\dots,K\}$, contains additional subcuboids that touch the faces $\{x_3=a_3^-\}$ and $\{x_3=a_3^+\}$ but are separated from the other faces of $\tilde{Q}_{i,j}$.
This extension is, however, omitted in order to keep the notation as simple as possible.
We note that the considered composition of $Q$ is adapted to ADI splitting schemes because these methods rely on subdomains of tensor structure for efficiency reasons, see \cite{KoeDiss,HoKoe}.
Each subdomain $\tilde{Q}_{i,j}$ is assumed to consist of a homogeneous material, and $\mu$ is not allowed to change in $\tilde{Q}_i$, meaning throughout
\begin{align}\label{assumptions parameters}
		\eps|_{\tilde{Q}_{i,j}},\mu|_{\tilde{Q}_{i,j}}\in\RR_{>0},\qquad \mu|_{\tilde{Q}_{i,0}}=\mu|_{\tilde{Q}_{i,l}}, 
\end{align}
for $i\in\{1,\dots,L\}$, $j\in\{0,\dots,K\}$, and $l\in\{1,\dots,K\}$.
(In contrast to \cite{Zeru22b}, we do not impose a monotonicity condition on $\eps$, but see also Remark~\ref{remark previous work}.)
At the interfaces between the subdomains, we impose several physical transmission conditions.
In particular, we allow for a nontrivial surface charge, see (\ref{def X2}) and Remark~\ref{remark X2}.

The main result of this paper is Theorem~\ref{theorem 7.5}.
It shows that the Peaceman-Rachford ADI scheme (\ref{eq 7.1}) from \cite{Namiki2000,ZhCZ2000,EiSc18} converges in $L^2$ with order $3/2-$ to the solution of the Maxwell system (\ref{Maxwell system}), provided that the initial data and the external current are sufficiently regular and satisfy appropriate boundary, transmission and divergence conditions.
To be more precise, for every number $\theta\in(0,3/2)$, the error of the ADI scheme converges with order $3/2-\theta$.
Note that we only impose assumptions on the data of the problem, but not on the solution.
Because of the irregularity of the material parameters, our error result provides a smaller convergence rate than the classical order $2$ from \cite{HoJaSc,Eilinghoff Diss,EiSc18,EiSc32,KoeDiss}.
The numerical examples in Section~\ref{section 7} show, however, that the order reduction predicted by our analysis is indeed observable in practice.

To prove Theorem~\ref{theorem 7.5}, we need a detailed knowledge of the regularity of the solutions to the Maxwell equations (\ref{Maxwell system}).
In other papers, the time-harmonic Maxwell equations are analyzed on more general and complicated heterogeneous polyhedral domains, see \cite{BonnetBhenDiaHazardLohrengel,CoDaNi,Ciarlet16,Ciarlet20,CiarletLefevreLohrengelNicaise,BonitoGuermondLuddens} for instance.
In particular, we apply Theorem~7.1 from \cite{CoDaNi} in the proof of Lemma~\ref{lemma 6.3}.
Similar to \cite{Zeru22b}, we present a detailed regularity analysis to establish a sharp explicit link from the maximal relative jump of the electric permittivity $\eps$ to the regularity of the solution, see Corollary~\ref{corollary 6.17} and Remark~\ref{remark 6.18}.
To the best of our knowledge, this regularity statement is new for our problem.
Note also that a time-dependent Maxwell system with additional surface current on a cuboid consisting of two homogeneous subcuboids is analyzed in \cite{DoeZer22}.

During the regularity analysis, we localize near the interior edges in the cuboid, and study an elliptic transmission problem for the first two components of the electric field, see Section~\ref{section3}, Section~\ref{section 5.2} and also \cite{CoDa,CoDaNi,Ciarlet20}.
To express the first two components of the electric field as the sum of a regular and a singular function, we additionally employ an interpolation result for subspaces from \cite{Zeru22a}, see Lemma~\ref{lemma 6.9}.
Arguing similar as in \cite{HoJaSc,ZerullaDiss}, we then show that a space $X_2$, see (\ref{def X2}), embeds into a space of functions with piecewise fractional Sobolev regularity.
The actual regularity and wellposedness result in Corollary~\ref{corollary 6.17} is then obtained by means of the regular state space $X_2$ and semigroup theory on that space, see Lemma~\ref{lemma 6.16}.

Having the regularity and wellposedness results from Corollary~\ref{corollary 6.17} and Remark~\ref{remark 6.18} at hand, the crucial ingredient in the proof for the main result in Theorem~\ref{theorem 7.5} is an estimate for a critical error term.
The needed analysis is established in Section~\ref{section 6.2} by means of a sophisticated $H^{\infty}$-functional calculus approach.

\subsubsection*{Structure of the paper}\

In Section~\ref{section2}, we recall useful function spaces related to Maxwell equations and introduce an appropriate analytical framework in which we can interpret (\ref{Maxwell system}) as a Cauchy problem.
Afterwards, we study an elliptic transmission problem on the unit disc which is useful for the analysis of the first two electric field components, see Section~\ref{section3}.
In the succeeding Section~\ref{section 4}, we present two lemmas on elliptic transmission problems with nontrivial contribution on the interfaces in $\tFint$.
By means of the findings in Sections~\ref{section3}-\ref{section 4}, we show a regularity result for the space $X_2$ in Section~\ref{section 5}.
The wellposedness and regularity of (\ref{Maxwell system}) is concluded in Section~\ref{section 5.4}.
The considered Peaceman-Rachford ADI scheme (\ref{eq 7.1}) is analyzed in Section~\ref{section 6}, and we establish a rigorous error result in Theorem~\ref{theorem 7.5}.
Finally, we present numerical error plots showing that a loss of convergence order for scheme (\ref{eq 7.1}) is observable in our heterogeneous model problem, see Section~\ref{section 7}.

\subsubsection*{Notation}\

We employ the same notation as in \cite{Zeru22b}:
For technical reasons, we use a partition of the cube $Q$ that is subordinate 
to $\overline{Q}=\cup_{i=1}^L\cup_{j=0}^K\overline{\tilde{Q}}_{i,j}$, and that is obtained by appropriate refinement.
The new partition consists of $N$ disjoint cuboids $Q_1,\dots,Q_N$ that touch the top and bottom faces of $Q$.
For subcuboids with a common interface, we additionally assume that the edges of the corresponding faces coincide.
Note that the parameters $\eps$ and $\mu$ are piecewise constant on the new partition. 
For a function $f$ on $Q$, we denote its restriction to a subcuboid $Q_i$ by $f^{(i)}$.

The interfaces of the fine partition $Q_1,\dots,Q_N$ are collected in a set $\Fint$, and the exterior faces are collected in $\Fext$.
The set of effective interfaces is defined as
\begin{align*}
	\Finteff:=\{\cF\sub Q\text{ is a face of }\tilde{Q}_{i,j},\ i\in\{1,\dots,L\},\ j\in\{0,\dots,K\}\}.
\end{align*}
Note that $\Finteff$ consists of all interfaces between submedia with different physical properties.
It is also important to associate a unit normal vector to each interface and exterior face.
For an interface $\cF\in\Fint\cup\Finteff$ being parallel to the $x_j-x_3$-plane, we choose its normal vector $\nu_{\cF}$ as the canonical unit vector $e_l$, $l\not=j\in\{1,2\}$.
For an external face $\tilde{\cF}\in\Fext$, the unit normal vector $\nu_{\tilde{\cF}}$ is chosen as the unit normal vector $\nu$ of $\partial Q$.
At several instances, we also use jumps of functions at an interface $\cF\in\Fint$:
Let $Q_i$ and $Q_j$ be the two cuboids sharing the interface $\cF$, and let $f$ be a function on $Q$ with $f^{(i)}$ and $f^{(j)}$ having well defined traces on $\cF$.
We then define the jump of $f$ at $\cF$ by $\llbracket f\rrbracket_{\cF}:=\tr_{\cF}f^{(i)}-\tr_{\cF}f^{(j)}$.

The open faces of $Q$ are denoted by 
\begin{align*}
	\Gamma_j^{\pm}:=\{x\in\partial Q\ |\ x_j=a_j^{\pm},\ x_l\in(a_l^-,a_l^+)\ \text{for }l\not=j\},\qquad \Gamma_j:=\Gamma_j^+\cup\Gamma_j^-,
\end{align*}
for $j\in\{1,2,3\}$.
The domain of a linear operator $A$ on a normed vector space $(X,\lVert\cdot\rVert)$ is called $\cD(A)$, and its graph norm is $\lVert\cdot\rVert_{\cD(A)}^2:=\lVert\cdot\rVert^2+\lVert A\cdot\rVert^2$.


\section{Analytical preliminaries}\label{section2}

This section is divided into two parts.
First we collect useful function spaces, analytical concepts, as well as an extension result.
In the second part, we introduce a convenient analytical framework in which we can interpret the Maxwell equations as an evolution equation on a state space.

\subsection{Useful function spaces and an extension result}\label{section 2.1}\

Here we follow Section~2.1 of \cite{Zeru22b}.
First, we recap important definitions and facts related to the $\curl$ and $\divv$ operators on the entire cuboid $Q$.
The concepts transfer directly to the considered subdomains of $Q$.
We use the Banach spaces
\begin{align*}
	H(\curl,Q)&:=\{\phi\in L^2(Q)^3\ |\ \curl\phi\in L^2(Q)^3\}, &&\norm{\phi}_{\curl}^2:=\norm{\phi}_{L^2}^2+\norm{\curl\phi}_{L^2}^2,\\
	H(\divv,Q)&:=\{\phi\in L^2(Q)^3\ |\ \divv \phi\in L^2(Q)\}, && \norm{\phi}_{\divv}^2:=\norm{\phi}_{L^2}^2+\norm{\divv \phi}_{L^2}^2.
\end{align*} 
The spaces $H_0(\curl,Q)$ and $H_0(\divv,Q)$ denote the subspaces of $H(\curl,Q)$ and $H(\divv,Q)$, respectively, that are obtained by completing the space $C_c^{\infty}(Q)^3$.
Theorems~I.2.5--I.2.6 in \cite{GiRa} yield that the normal trace operator $v\mapsto v\cdot\nu|_{\partial Q}$ extends continuously to $H(\divv,Q)$ with kernel $H_0(\divv,Q)$.
It maps into $H^{-1/2}(\partial Q)$ and gives rise to Green's formula
\begin{align*}
	\int_Qv\cdot\nabla\varphi\dd x+\int_Q(\divv v)\varphi\dd x=\langle v\cdot\nu,\varphi\rangle_{H^{-1/2}(\partial Q)\times H^{1/2}(\partial Q)},
\end{align*}
for $v\in H(\divv,Q)$, $\varphi\in H^1(Q)$.

Theorems~I.2.11--I.2.12 in \cite{GiRa} provide a continuous extension of the tangential trace operator $v\mapsto v\times\nu|_{\partial Q}$ to $H(\curl,Q)$ with kernel $H_0(\curl,Q)$ and range in $H^{-1/2}(\partial Q)^3$.
The associated Green's formula is given by 
\begin{align*}
	\int_Q(\curl v)\cdot\varphi\dd x-\int_Qv\cdot\curl\varphi\dd x=\langle v\times\nu,\varphi\rangle_{H^{-1/2}(\partial Q)\times H^{1/2}(\partial Q)},
\end{align*}
for $v\in H(\curl,Q)$ and $\varphi\in H^1(Q)^3$.

For the error analysis in Section~\ref{section 6}, we employ extrapolation theory for operators, see Section~V.1.3 in \cite{Amann} and Section~2.10 in \cite{TuWe}.
We sketch the basic concept:
Let $A$ be a closed linear operator with dense domain and nonempty resolvent set $\rho(A)$ on a Banach space $(X,\lVert\cdot\rVert_X)$.
Let $\lambda\in\rho(A)$.
The extrapolation space $X_{-1}^A$ of $A$ is obtained by completing $X$ with respect to the norm $\lVert\cdot\rVert_{X_{-1}^A}=\lVert(\lambda I-A)^{-1}\cdot\rVert_X$.
Note that $A$ has a unique continuous extension $A_{-1}$, mapping from $X$ into $X_{-1}^A$.
The latter mapping is called extrapolation of $A$ to $X$.
Similarly, the resolvent $(\lambda I-A)^{-1}$ extends continuously to $(\lambda I-A_{-1})^{-1}$, mapping from $X_{-1}^A$ to $X$.

At several instances, we moreover use real interpolation spaces, see Section~1.1 in \cite{Lunardi} for instance.
They are in particular useful to define Sobolev spaces of fractional order, see \cite{LionsMagenes,Triebel72}.
Let $d\in\{1,2,3\}$,  and $O\sub\RR^d$ open with a Lipschitz boundary.
The spaces $H^1_0(O)$ and $H_0^2(O)$ denote the closure of the space $C_c^{\infty}(O)$ in $H^1(O)$ and $H^2(O)$, respectively.
We set
\begin{align}
	H^s(O)&:=(L^2(O),H^2(O))_{s/2,2}, &&H^{\theta}_0(O):=(L^2(O),H^1_0(O))_{\theta,2},\label{eq 2.0}\\
	H^{1/2}_{00}(O)&:=(L^2(O),H^1_0(O))_{1/2,2}, &&H^{3/2}_{00}(O):=(L^2(O),H^3(O)\cap H^2_0(O))_{1/2,2},\nonumber
\end{align}
for $s\in(0,2)$ and $\theta\in(0,1)\setminus\{1/2\}$.
Note additionally that $H^{\theta}(O)=H^{\theta}_0(O)$ for $\theta\in(0,1/2)$.
(This fact can be proven with Corollary~1.4.4.5 in \cite{Grisvard}.)

To deal with functions that are regular on each subcuboid but irregular across the interfaces, we employ piecewise Sobolev spaces.
We define
\begin{align*}
	PH^q(Q)&:=\{f\in L^2(Q)\ |\ f^{(i)}\in H^q(Q_i),\; i\in\{1,\dots,N\}\}, &&q\in[0,2],\\
	PH^s_{\Gamma^*}(Q)&:=\{f\in PH^s(Q)\ |\ f^{(i)}=0\text{ on }\partial Q_i\cap\Gamma^*,\ i\in\{1,\dots,N\}\}, &&s\in(1/2,2],\\
	\norm{f}_{PH^q}^2&:=\sum_{i=1}^N\lVert f^{(i)}\rVert_{H^q(Q_i)}^2,\qquad \lVert g\rVert_{PH^s_{\Gamma^*}}:=\lVert g\rVert_{PH^s}, 
\end{align*}
for $f\in PH^q(Q)$, $g\in PH^s_{\Gamma^*}(Q)$, where $\Gamma^*$ is a nonempty union of some of the faces of $Q$.
We close this subsection with a useful extension result that is employed in the proof of Lemma~\ref{lemma 6.3}.
For the proof, we modify extension techniques and arguments from the proof of Lemma~3.1 in \cite{EiSc32}.

\begin{lemma}\label{lemma 2.1}
	Let $\check{Q}=(0,1)^3$, $\nu$ be the unit exterior normal vector of $\partial \check{Q}$, $\cF$ be a face of $\check{Q}$, and $g\in H^{3/2}_{00}(\cF)$.
	There is a function $u\in H^3(\check{Q})\cap H^1_0(\check{Q})$ with $\partial_{\nu}u=0$ on all faces $\cF'\not=\cF$ of $\check{Q}$ and $\partial_{\nu}u=g$ on $\cF$.
	Furthermore, $\Delta u$ has a trace in $H^{1/2}_{00}(\tilde{\cF})$ on all faces $\tilde{\cF}$ of $\check{Q}$.
	The function $u$ can be estimated in norm via $\lVert u\rVert_{H^3(\check{Q})}\leq C\lVert g\rVert_{H^{3/2}_{00}(\cF)}$, involving a uniform constant $C>0$.
	\begin{proof}
		1) Let $\cF=\{0\}\times[0,1]^2\cong[0,1]^2$, and let $\chi:[0,\infty)\to[0,1]$ be a smooth cut-off function with $\chi=1$ on $[0,1/2]$ and $\supp\chi\sub[0,3/4]$.
		We further employ a positive definite selfadjoint operator $L$ on $L^2(\cF)$ with $\cD(L)=H^3(\cF)\cap H^2_0(\cF)$.
		(Such an operator exists, see Section~124 in \cite{RiNa73} and Section~1.2.1 in \cite{LionsMagenes}.) 
		By assumption, $g\in\cD(L^{1/2})$.
		We further note that $L^{1/3}$ generates an analytic semigroup $(\e^{-tL^{1/3}})_{t\geq0}$ on $L^2(\cF)$.
		We claim that
		\begin{align*}
			u(x_1,x_2,x_3):=-\chi(x_1)x_1\big(\e^{-x_1L^{1/3}}g\big)(x_2,x_3),\qquad (x_1,x_2,x_3)\in\check{Q},
		\end{align*}
		is the desired extension of $g$.
		We prove this assertion in the next two steps.
		
		2) In the following, $C>0$ is a uniform constant that is allowed to change from line to line.
		To derive the asserted regularity statement, we note that $u$ is smooth inside $\check{Q}$ as $(\e^{-tL^{1/3}})_{t\geq0}$ is analytic.
		We moreover use the embeddings
		\begin{align*}
			\cD(L)\hookrightarrow H^3(\cF),\quad\cD(L^{2/3})\hookrightarrow H^2(\cF),\quad \cD(L^{1/3})\hookrightarrow H^1(\cF),
		\end{align*}
		and identify $\cF$ with $[0,1]^2$.
		We calculate
		\begin{align}\label{eq 2.1}
			\partial_1u&=-\big(\chi+x_1\chi'-\chi x_1L^{1/3}\big)\e^{-x_1L^{1/3}}g,\\
			\partial_1^2u&=-\big(2\chi'+x_1\chi''-2\chi'x_1L^{1/3}-2\chi L^{1/3}+\chi x_1L^{2/3}\big)\e^{-x_1L^{1/3}}g,\nonumber\\
			\partial_1^3u&=-\big(3\chi''+x_1\chi'''-(3\chi''x_1+6\chi') L^{1/3}+(3\chi'x_1+3\chi)L^{2/3}-\chi x_1L\big)\e^{-x_1L^{1/3}}g.\nonumber
		\end{align}
	By means of Proposition~6.4 in \cite{Lunardi}, we then obtain the desired relations
	\begin{align*}
		\sum_{j=0}^3\int_0^1&\lVert\partial_1^ju(x_1,\cdot)\rVert_{H^{3-j}(\cF)}^2\dd x_1\leq C\sum_{j=0}^3\int_0^1\lVert L^{1-j/3}\partial_1^ju(x_1,\cdot)\rVert_{L^2(\cF)}^2\dd x_1\\
		&\leq C\int_0^1\big(\lVert x_1 L^{2/3}\e^{-x_1L^{1/3}}L^{1/3}g\rVert_{L^2(\cF)}^2+\lVert L^{1/6}\e^{-x_1L^{1/3}}L^{1/2}g\rVert_{L^2(\cF)}^2\big)\dd x_1\\
		&\leq C\lVert L^{1/2}g\rVert_{L^2(\cF)}^2\leq C\lVert g\rVert_{H^{3/2}_{00}(\cF)}^2.
	\end{align*}
	
	3) It remains to study the behavior of $u$ on the boundary of $\check{Q}$.
	By construction, $u=0$ on $\cF$.
	The semigroup $(\e^{-tL^{1/3}})_{t\geq0}$ being analytic, we further have $u(x_1,\cdot)\in\cD(L)\sub H^2_0(\cF)$ for $x_1>0$.
	Taking also the choice of $\chi$ into account, we conclude the asserted homogeneous Dirichlet and Neumann boundary conditions on all faces of $\check{Q}$, except $\cF$.
	Formula~(\ref{eq 2.1}) furthermore shows the desired extension property $\partial_1u(0,\cdot)=-g$.
	
	We finally deal with $\Delta u$, and show that $(\Delta u)|_{\cF}\in H^{1/2}_{00}(\cF)$.
	(The other faces of $\check{Q}$ can be treated in the same way.)
	By part~2), $\partial_j^2u\in H^1(\check{Q})$ for $j\in\{1,2,3\}$.
	Combining the fact $u=0$ on $\cF$ with Lemma~2.1 in \cite{EiSc18}, we consequently infer that $\partial_j^2u=0$ on $\cF$ for $j\in\{2,3\}$.
	In view of the analyticity of $(\e^{-tL^{1/3}})_{t\geq0}$ and (\ref{eq 2.1}), the function $\partial_1^2u$ moreover vanishes on all faces of $\check{Q}$, except $\cF$.
	By the trace method in Section~I.3 of \cite{LionsMagenes}, we thus conclude that $\partial_1^2u\in H^{1/2}_{00}(\cF)$.
	\end{proof}
\end{lemma}

\subsection{Analytical framework for the Maxwell equations}\label{section 2.2}\

In this section, we use concepts from Section~2.2 in \cite{Zeru22b} and Sections~2--3 in \cite{EiSc32}.
The Maxwell system (\ref{Maxwell system}) is interpreted as an evolution equation on the space $X=L^2(Q)^6$, which is equipped with the weighted inner product
\begin{align*}
	\Big(\begin{pmatrix}\mE\\\mH\end{pmatrix},\begin{pmatrix}\tilde{\mE}\\\tilde{\mH}\end{pmatrix}\Big):=\int_Q\eps\mE\cdot\tilde{\mE}+\mu\mH\cdot\tilde{\mH}\dd x,\qquad \begin{pmatrix}\mE\\\mH\end{pmatrix},\begin{pmatrix}\tilde{\mE}\\\tilde{\mH}\end{pmatrix}\in X.
\end{align*}
The induced norm is denoted by $\lVert\cdot\rVert$.
In presence of (\ref{assumptions parameters}), this norm is equivalent to the standard $L^2$-norm.
The linear Maxwell operator is defined as
\begin{align}
	M:=\begin{pmatrix}0 & \tfrac{1}{\eps}\curl\\ -\tfrac{1}{\mu}\curl & 0\end{pmatrix},\qquad \cD(M):=H_0(\curl,Q)\times H(\curl,Q).\label{def Maxwell}
\end{align}
Note that $\cD(M)$ prescribes tangential interface conditions, as well as the electrical boundary condition.
To incorporate also relevant normal transmission, boundary and divergence conditions, we use the spaces
 \begin{align}\label{def X2}
 	X_0=\{(\mE,\mH)\in L^2(Q)^6\ |\ &\divv(\mu\mH)=0,\ \mu\mH\cdot\nu=0\text{ on }\partial Q\},\\
 	X_2:=\{(\mE,\mH)\in \cD(M^2)\cap X_0\ |\ &\divv(\eps\mE|_{\tilde{Q}_{i,l}})\in H^1_{00}(\tilde{Q}_{i,l}),\ \divv(\eps\mE|_{\tilde{Q}_{i,l}})=0\nonumber\\
 	&\text{on }\Gamma_3\cap\partial\tilde{Q}_{i,l}\text{ for }i\in\{1,\dots,L\},l\in\{0,\dots,K\},\nonumber\\
 	& \llbracket\eps\mE\cdot\nu_{\cF}\rrbracket_{\cF}\in H_{00}^{3/2}(\cF)\text{ for }\cF\in\Finteff\}.\nonumber
 \end{align}
(The set of effective interfaces $\Finteff$ and the domains $\tilde{Q}_{i,l}$ are defined in Section~\ref{section1}, and $H^{3/2}_{00}$ is introduced in (\ref{eq 2.0}).)
Here, $H^1_{00}(\tilde{Q}_{i,l})$ denotes the space of $H^1$-functions on $\tilde{Q}_{i,l}$, whose traces on each face $\cF$ of $\tilde{Q}_{i,l}$ belong to $H^{1/2}_{00}(\cF)$.
While $X_0$ is complete with respect to the norm in $X$, the space $X_2$ is complete with the norm
\begin{align*}
	&\norm{(\mE,\mH)}_{X_2}^2:=\norm{(\mE,\mH)}_{\cD(M^2)}^2+\sum_{i=1}^N\norm{\divv(\eps^{(i)}\mE^{(i)})}_{H^1(Q_i)}^2\\
	&\quad+\sum_{i=1}^{L}\sum_{l=0}^K\sum_{\cF\text{ face of }\tilde{Q}_{i,l}}\norm{\divv\mE|_{\tilde{Q}_{i,l}}}_{H_{00}^{1/2}(\cF)}^2+\sum_{\cF\in\Finteff}\norm{\llbracket\eps\mE\cdot\nu_{\cF}\rrbracket_{\cF}}_{H_{00}^{3/2}(\cF)}^2.
\end{align*}

\begin{remark}\label{remark X2}
	By considering the Maxwell system (\ref{Maxwell system}) on $X_2$, we prescribe the transmission conditions
	\begin{align*}
		\llbracket\mE\times\nu_{\cF}\rrbracket_{\cF}=0=\llbracket\mH\times\nu_{\cF}\rrbracket_{\cF},\qquad \llbracket\mu\mH\cdot\nu_{\cF}\rrbracket_{\cF}=0,
	\end{align*}
	on each interface $\cF\in\Fint$.
	The possibly nonzero jump $\llbracket\eps\mE\times\nu_{\cF}\rrbracket_{\cF}$ has the physical meaning of a surface charge density, see Section~4.12 in \cite{Stratton}, Section~1.1.3 in \cite{BornWolf}, and Section~I.5 in \cite{Jackson} for instance.\hfill$\lozenge$
\end{remark}

In Section~\ref{section 5}, we prove that $X_2$ embeds into the space of piecewise $H^{2-\kappa}$-regular functions with appropriate $\kappa>0$, see Proposition~\ref{proposition 6.15}.
To study the Maxwell system (\ref{Maxwell system}) in $X_2$, we define $M_2$ to be the part of $M$ in $X_2$.
The domain of $M_2$ is determined in the next lemma.
To that end, we transfer a part of Proposition~3.2 in \cite{EiSc32} to our setting.

\begin{lemma}\label{lemma 2.2}
	Let $\eps$ and $\mu$ satisfy (\ref{assumptions parameters}).
	The operator $M_2$ has the domain $\cD(M_2)=\cD(M^3)\cap X_2$.
	\begin{proof}
		In view of the definition of $X_2$, it suffices to show $\cD(M^3)\cap X_2\sub\cD(M_2)$.
		Let $(\mE,\mH)\in\cD(M^3)\cap X_2$, and put $(u,v)=M(\mE,\mH)$.
		Then $(u,v)\in\cD(M^2)$ with $\divv(\eps u)=0$ and $\divv(\mu v)=0$.
		In particular, $\llbracket\eps u\cdot\nu_{\cF}\rrbracket_{\cF}=0$ on $\cF\in\Finteff$.
		The boundary condition $\mu v\cdot\nu=0$ on $\partial Q$ is a consequence of $\mE\times\nu=0$ on $\partial Q$, see Remark~I.2.5 in \cite{GiRa}.
	\end{proof}
\end{lemma}

$M_2$ generates a strongly continuous semigroup on $X_2$, see Lemma~\ref{lemma 6.16}.
From this we can then conclude that the Maxwell system (\ref{Maxwell system}) possesses classical solutions with piecewise $H^{2-\kappa}$-regularity, see Remark~\ref{remark 6.18}.

We close this section with an important remark that is employed several times throughout this paper.

\begin{remark}\label{remark previous work}
    Inspecting the arguments in \cite{Zeru22b}, we see that the results from there also hold in the current material configuration from Section~\ref{section1}.
    This is due to the fact that \cite{Zerucalc1} already allows for assumption~(\ref{eq 3.2}).
\end{remark}


\section{Analysis of an elliptic transmission problem}\label{section3}

In this section, we analyze a Laplace operator on the unit disc with transmission conditions that are motivated by the first two components of the electric field.
Note in particular that the first two components of the electric field are continuous in two space dimensions, and discontinuous in the remaining one.
The analysis follows the one in Section~3 of \cite{Zeru22b}, whence we only sketch identical parts.
Our reasoning is inspired by \cite{Kellogg71}, which treats different transmission conditions, see also \cite{CoDaNi,Ciarlet16,Ciarlet20,Jochmann1999}.
Our goal is a representation of the domain of the Laplace operator in (\ref{eq 3.4}) as the direct sum of a space of regular functions and a one-dimensional span of an explicitly given irregular function.
As in \cite{Zeru22b}, we aim for an explicit result in terms of the size of the jumps of the material parameter $\eps$, see Proposition~\ref{proposition 3.5}.

\subsection{Introduction of a two-dimensional Laplacian with transmission conditions}\label{section 3.1}\

We recall some of the assumptions and constructions in Section~3.2 of \cite{Zeru22b}.
We fix an interior edge $e_{\inn}$ in $Q$ at which $\eps$ has a strong discontinuity, and let $Q_{\inn,1},\dots,Q_{\inn,4}$ be the adjacent cuboids to $e_{\inn}$.
Strong discontinuity of $\eps$ means that $\eps$ attains a different value on one cuboid compared to the remaining three, cf.~Definition~3.4 in \cite{Zeru22b}.
After translation and scaling, we can assume the identity
\begin{align}\label{eq 3.1}
	e_{\inn}=\{(0,0)\}\times[0,1].
\end{align}
The notation $\eps_{\inn}^{(i)}$ means the restriction $\eps|_{Q_{\inn,i}}$.
Owing to symmetry, we can additionally assume the configuration
\begin{align}\label{eq 3.2}
	\eps_{\inn}^{(1)}=\eps_{\inn}^{(2)}=\eps_{\inn}^{(3)}\not=\eps_{\inn}^{(4)},
\end{align}
see (\ref{assumptions parameters}).
We denote by $D$ the unit disc, and assume that the cylinder $D\times[0,1]$ touches no second interior edge.
After rotating, the representation
\begin{align}\label{eq 3.3}
	(D\times(0,1))\cap Q_{\inn,i}=&D_{\inn,i}\times(0,1),\\
	D_{\inn,i}:=&\{(r\cos\varphi,r\sin\varphi)\ |\ r\in(0,1),\ \varphi\in I_{\inn,i}\},\\
	I_{\inn,1}:=&(0,\tfrac{\pi}{2}),\quad I_{\inn,2}:=(\tfrac{\pi}{2},\pi),\quad I_{\inn,3}:=(\pi,\tfrac{3}{2}\pi),\quad I_{\inn,4}:=(\tfrac{3}{2}\pi,2\pi),\nonumber
\end{align}
is valid for $i\in\{1,\dots,4\}$.
In the following, we interpret $\eps_{\inn}$ as a piecewise constant function on the partition $D_{\inn,1}\cup\dots\cup D_{\inn,4}$.
The interfaces in $D$ are denoted by
\begin{align*}
	\cF^D_{k}:=\partial D_{\inn,k}\cap\partial D_{\inn,k+1},\quad \cF^D_{\inn,4}:=\partial D_{\inn,1}\cap\partial D_{\inn,4},\quad k\in\{1,2,3\}.
\end{align*}
The notion of restrictions of functions, piecewise regularity, and the jump $\llbracket\cdot\rrbracket$ at an interface are transferred accordingly to this partition of $D$.

In this section, we aim to represent the domain of the two-dimensional Laplacian
\begin{align}
	(L_{\inn}\psi)^{(i)}&:=\Delta\psi^{(i)},\label{eq 3.4}\\
\psi\in\cD(L_{\inn})&:=\{\psi\in PH^1(D)\ |\ \Delta\psi^{(i)}\in L^2(D_{\inn,i}),\nonumber\\
&\qquad\llbracket\eps_{\inn}\psi\rrbracket_{\cF^D_{l}}=0=\llbracket\nabla\psi\cdot\nu\rrbracket_{\cF^D_{l}},\;\llbracket\psi\rrbracket_{\cF^D_{k}}=0=\llbracket\eps_{\inn}\nabla\psi\cdot\nu\rrbracket_{\cF^D_{k}}\nonumber\\
&\qquad\text{for }i\in\{1,\dots,4\},l\in\{2,4\},k\in\{1,3\}\},\nonumber
\end{align}
as the direct sum of piecewise $H^2$-regular functions on $D$ and the span of an explicitly given singular function, see Proposition~\ref{proposition 3.5}.

We note that the operator $L_{\inn}$ is invertible with compact resolvent. It is furthermore selfadjoint with respect to the inner product
\begin{align}
	(f,g)_{\eps_{\inn},D}:=\int_D\eps_{\inn}fg\dd x,\qquad f,g\in L^2(D).\label{eq 3.4.1}
\end{align}
(In fact, bijectivity is obtained via a Lax-Milgram-Lemma argument, and symmetry is concluded by an integration by parts.)

\subsection{A one-dimensional eigenvalue problem}\label{section 3.2}\

To determine the strongest arising radial singularity of functions in the domain $\cD(L_{\inn})$, we study the eigenvalue problem
\begin{align}
	(\psi^{(i)})''(\varphi)&=-\kappa^2\psi^{(i)}(\varphi) \qquad\text{for }\varphi\in I_{\inn,i},\ i\in\{1,\dots,4\},\label{eq 3.5}\\
	\eps_{\inn}^{(1)}\psi^{(1)}(0)&=\eps_{\inn}^{(4)}\psi^{(4)}(2\pi),\quad (\psi^{(1)})'(0)=(\psi^{(4)})'(2\pi),\nonumber\\
	\psi^{(1)}(\tfrac{\pi}{2})&=\psi^{(2)}(\tfrac{\pi}{2}),\quad (\psi^{(1)})'(\tfrac{\pi}{2})=(\psi^{(2)})'(\tfrac{\pi}{2}),\nonumber\\
	\psi^{(2)}(\pi)&=\psi^{(3)}(\pi),\quad (\psi^{(2)})'(\pi)=(\psi^{(3)})'(\pi),\nonumber\\
	\psi^{(3)}(\tfrac{3}{2}\pi)&=\psi^{(4)}(\tfrac{3}{2}\pi)\quad \eps_{\inn}^{(1)}(\psi^{(3)})'(\tfrac{3}{2}\pi)=\eps_{\inn}^{(4)}(\psi^{(4)})'(\tfrac{3}{2}\pi).\nonumber
\end{align}
in the next two lemmas, cf.~\cite{Kellogg74}.

\begin{lemma}\label{lemma 3.2}
	Let $\eps_{\inn}$ satisfy (\ref{eq 3.2}). 
	Then (\ref{eq 3.5}) has a countable set of eigenvalues $0<\kappa_{1}^2\leq\kappa_{2}^2\leq\dots\to\infty$, and associated piecewise smooth eigenfunctions $\psi_{1},\psi_{2},\dots$, forming an orthonormal basis of $L^2(0,2\pi)$ with respect to the inner product $L^2$-inner product with weight $\eps_{\inn}$.
	\begin{proof}
		Consider the closed, symmetric, and positive definite bilinear form
		\begin{align*}
			a(\psi,\tilde{\psi})&:=\sum_{i=1}^4\int_{I_{\inn,i}}\eps_{\inn}\big[\nabla\psi^{(i)}\cdot\nabla\tilde{\psi}^{(i)}+\psi\tilde{\psi}\big]\dd \varphi,\\
			\cD(a)&:=\{\psi\in PH^1(0,2\pi)\ |\ \psi\text{ satisfies the zero-order transmission}\\
			&\qquad\qquad\qquad\qquad \text{conditions in }(\ref{eq 3.5})\}
		\end{align*} 
		on $L^2(0,2\pi)$. The latter space is equipped with the $L^2$-inner product with weight $\eps_{\inn}$.
		The operator 
		\begin{align*}
		(Tf)|_{I_{\inn,i}}&:=f^{(i)}-(f^{(i)})'',\qquad i\in\{1,\dots,4\},\\
		f\in\cD(T)&:=\{f\in PH^2(0,2\pi)\ |\ f\text{ satisfies the transmission conditions in }
		(\ref{eq 3.5}) \}
		\end{align*}
		is then associated with $a$.
		By Theorem~VI.2.6 in \cite{Kato}, $T$ is hence positive definite, selfadjoint, and invertible on $L^2(0,2\pi)$.
		Taking also into account that $T$ has a compact resolvent, the spectral theorem for selfadjoint operators with compact resolvent provides a countable set of positive eigenvalues and an associated  orthonormal basis of eigenfunctions for $T$.
		The definition of $T$ now implies all asserted statements, except the bound for the smallest eigenvalue.
		
		Suppose $\kappa_1^2$ was nonpositive, and let $\psi_1\not=0$ be an associated eigenfunction.
		Integrating the eigenvalue-eigenfunction relation by parts, we then infer the identity
		\begin{align*}
			\sum_{i=1}^4\int_{I_{\inn,i}}\eps_{\inn}\lvert\psi_{1}'\rvert^2\dd\varphi=0.
		\end{align*}
		This means that $\psi_1$ is piecewise constant, contradicting the zero order transmission conditions in (\ref{eq 3.5}).
	\end{proof}
\end{lemma} 

The next statement is a counterpart of Lemma~3.5 in \cite{Zeru22b}. It provides a crucial sharp upper bound for the first eigenvalue of (\ref{eq 3.5}).
The bound is given by the number $\sskappa\in[0,1/3)$ with
\begin{align}
\max_{\substack{i\in\{1,\dots,L\},\\ l\in\{1,\dots,K\}}}\frac{(\eps|_{\tilde{Q}_{i,l}}-\eps|_{\tilde{Q}_{i,0}})^2}{\eps|_{\tilde{Q}_{i,l}}\eps|_{\tilde{Q}_{i,0}}}=\frac{4\sin^2(\sskappa\pi)}{\cos(\tfrac{3}{2}\sskappa\pi)\cos(\tfrac{1}{2}\sskappa\pi)},\label{eq 3.6}
\end{align}
where $(\tilde{Q}_{i,l})_{i,l}$ are the submedia from Section~\ref{section1}.
Note that $\sskappa$ is uniquely determined by (\ref{eq 3.6}).
We moreover observe the structural similarity to the defining relation for $\skappa$, see \cite{Zeru22b}.

\begin{lemma}\label{lemma 3.3}
	Let $\eps_{\inn}$ satisfy (\ref{eq 3.2}). Then the inequalities $\kappa_{1}\leq \sskappa<1\leq\kappa_2$ are true.
\end{lemma}

As the proof of Lemma~\ref{lemma 3.3} merely consists in a long calculation that uses well known techniques, we skip it, and refer to \cite{calculations}.

\subsection{Analysis of a two-dimensional Laplacian with transmission conditions}\label{section 3.3}\

Our goal is a direct decomposition of the domain $\cD(L_{\inn})$ from (\ref{eq 3.4}) into a space of $H^2$-regular functions, and the span of a radially singular function.
To that end, let $\chi:[0,\infty)\to[0,1]$ be a smooth cut-off function which is equal to one on $[0,1/4)$, decreases monotonically, and is supported on $[0,1/2]$.
In the spirit of \cite{Kellogg71}, we define the supplementary spaces
\begin{align}
	\check{M}_{\inn}:=\{&\psi\in\bigcap_{i=1}^4C^1(\overline{D_{\inn,i}})\cap PH^2(D)\ |\ r\partial_x^2\psi^{(i)},\ r\partial_x\partial_y\psi^{(i)},\ r\partial_y^2\psi^{(i)}\in C(\overline{D_{\inn,i}})\text{ and}\nonumber\\
	&\text{tend to }0\text{ as } r\to0,\ \psi^{(i)}\in C^3(\overline{D_{\inn,i}}\setminus\{0\}),\ \psi=0\text{ on }\partial D,\ 
	\nonumber
	\\
	&\llbracket\eps_{\inn}\psi\rrbracket_{\cF^D_{l}}=0=\llbracket\partial_{\nu}\psi\rrbracket_{\cF^D_{l}},
	\llbracket\psi\rrbracket_{\cF^D_{k}}=0=\llbracket\eps_{\inn}\partial_{\nu}\psi\rrbracket_{\cF^D_{k}}
	\nonumber\\
	&\text{for }i\in\{1,\dots,4\},l\in\{2,4\},k\in\{1,3\}\},\nonumber\\
	\check{N}_{\inn}&:=\spanv\{\chi(r)r^{\kappa_{1}}\psi_{1}(\varphi)\}.\label{eq 3.6.1}
\end{align}

For the first space, the next lemma provides a useful a-priori energy estimate involving the Laplacian $L_{\inn}$. 
As the proof is obtained by straightforward adaptions of the one for Lemma~2.2 in \cite{Kellogg71}, we omit it. 

\begin{lemma}\label{lemma 3.4}
	Let $\eps_{\inn}$ satisfy (\ref{eq 3.2}).
	There is a constant $C=C(\eps_{\inn})>0$ with
	\begin{align*}
		\norm{\psi}_{PH^2(D)}\leq C\Big(\norm{\psi}_{L^2(D)}+\norm{L_{\inn}\psi}_{L^2(D)}\Big),\qquad \psi\in\check{M}_{\inn}.
	\end{align*}
\end{lemma}

We close this section with the desired decomposition result for the domain $\cD(L_{\inn})$ from (\ref{eq 3.4}). 
As the proof transfers the reasoning from the three-dimensional setting in Theorem~5.1 of \cite{Kellogg71} with different transmission conditions to the present transmission conditions in the two-dimensional case, we only sketch the relevant different parts. 
Compare also Theorem~2.1 in \cite{Kellogg71} and the reasoning in \cite{Grisvard75,Lemrabet78}.

\begin{proposition}\label{proposition 3.5}
	Let $\eps_{\inn}$ satisfy (\ref{eq 3.2}). The domain $\cD(L_{\inn})$ can be decomposed into
	\begin{align*}
		\cD(L_{\inn})=\overline{\check{M}_{\inn}}^{\norm{\cdot}_{PH^2}}\oplus\check{N}_{\inn}.
	\end{align*}
	\begin{proof}
		1) We show that $L_{\inn}$ maps the space $W:=\overline{\check{M}_{\inn}}\oplus\check{N}_{\inn}$ onto $L^2(D)$. As $L_{\inn}$ is injective, this implies the asserted statement. 
		To derive the surjectivity of $L_{\inn}:W\to L^2(D)$, we prove that the orthogonal complement $\cN$ of the space $\eps_{\inn}L_{\inn}(W)$ in $L^2(D)$ is trivial.
		
		Let $v\in\cN$. 
		Consider the function $v_r:[0,2\pi)\to\RR$, $\varphi\mapsto v(r\cos\varphi,r\sin\varphi)$, $r\in(0,1)$, and abbreviate $s(\varphi):=(\cos\varphi,\sin\varphi)$, $\varphi\in[0,2\pi)$.
		In the eigenbasis $\{\psi_{k}\ |\ k\in\NN\}$ of system (\ref{eq 3.5}), $v_r$ has the expansion
		\begin{align*}
			v_r=\sum_{k=1}^{\infty}\alpha_k(r)\psi_{k},\qquad \alpha_k(r):=\int_{0}^{2\pi}\eps_{\inn} v(rs(\varphi))\psi_{k}(\varphi)\dd \varphi,\qquad r\in(0,1).
		\end{align*}
		Analogously to the proof of Theorem~5.1 in \cite{Kellogg71}, the estimate
		\begin{align}
			\sum_{k=1}^{\infty}\int_0^1r\lvert\alpha_k(r)\rvert^2\dd r<\infty,\label{eq 3.7}
		\end{align}
		can be verified.
		
		2) Let now $\tilde{\chi}\in C_c^{\infty}(0,1)$, $k\in\NN$, and set
		\begin{align*}
			u_k(rs(\varphi)):=\tilde{\chi}(r)\psi_{k}(\varphi),\qquad r\in(0,1),\ \varphi\in[0,2\pi).
		\end{align*}
		By construction, $u_k$ is an element of $\check{M}_{\inn}$.
		The representation of the Laplacian in polar coordinates and the choice of $\psi_{k}$ as an eigenfunction of (\ref{eq 3.5}) yield the identity
		\begin{align*}
			\Delta u_k^{(i)}(r,\varphi)
			=\frac{1}{r^2}(r\tilde{\chi}'(r)+r^2\tilde{\chi}''(r)-\kappa_{k}^2\tilde{\chi}(r))\psi_{k}^{(i)}(\varphi).
		\end{align*}
		As $v$ is an element of $\cN$, we then arrive at the relations
		\begin{align}
			0&=(v,\eps_{\inn}L_{\inn} u_k)=\int_0^1\int_0^{2\pi}\tfrac{\eps_{\inn}}{r}(r\tilde{\chi}'(r)+r^2\tilde{\chi}''(r)-\kappa^2_{k}\tilde{\chi}(r))\psi_{k}(\varphi)v(rs(\varphi))\dd \varphi\dd r\nonumber\\
			&=\int_0^{1}\tfrac{1}{r}(r^2\tilde{\chi}''(r)+r\tilde{\chi}'(r)-\kappa_{k}^2\tilde{\chi}(r))\alpha_k(r)\dd r.\label{eq 3.7.1}
		\end{align}
		An integration by parts now shows the identity
		\begin{align}
			r(r\alpha_k')'-\kappa_{k}^2\alpha_k=0,\qquad r\in(0,1),\label{eq 3.8}
		\end{align}
		as $\tilde{\chi}\in C_c^{\infty}(0,1)$ is chosen arbitrary.
		We consequently infer the formula
		$\alpha_k(r)=a_kr^{\kappa_k}+b_kr^{-\kappa_k}$ with real numbers $a_k, b_k$.
		
		3) We next deduce that $b_k=0$ for $k\geq2$.
		Lemma~\ref{lemma 3.3} implies $\kappa_k\geq 1$ for $k\geq2$, while (\ref{eq 3.7}) gives rise to the relations
		\begin{align*}
			\infty>\sum_{k=1}^{\infty}\int_0^1r\alpha_k(r)^2\dd r=\sum_{k=1}^{\infty}\int_0^1a_k^2r^{1+2\kappa_k}+2a_kb_kr+b_k^2r^{1-2\kappa_r}\dd r.
		\end{align*}
		This shows that $b_k=0$ for $k\geq2$.
		
		Consider now the function $\tilde{u}_1(r,\varphi):=\chi(r)r^{\kappa_{1}}\psi_{1}(\varphi)$ (recall that $\chi$ is the cut-off function from the definition of $\check{N}_{\inn}$).
		As $\psi_1$ is an eigenfunction of (\ref{eq 3.5}), a calculation leads to the formula
		\begin{align*}
			\Delta\tilde{u}_1^{(i)}=\Big(\chi''r^{\kappa_{1}}+\chi'(2\kappa_{1}+1)r^{\kappa_{1}-1}\Big)\psi_{1}^{(i)}.
		\end{align*}
		Using that $v$ is an element of $\cN$, we conclude the equations
		\begin{align*}
			0&=(v,\eps_{\inn}L_{\inn}\tilde{u}_1)_{L^2(D)}\\
			&=\int_0^1\int_0^{2\pi}\Big(rv(rs(\varphi))(\chi''r^{\kappa_{1}}+\chi'(2\kappa_{1}+1)r^{\kappa_{1}-1})\eps_{\inn}\psi_{1}(\varphi)\Big)\dd\varphi\dd r\\
			&=\int_0^1\Big(\chi''r^{\kappa_{1}+1}+\chi'(2\kappa_{1}+1)r^{\kappa_{1}}\Big)\alpha_1(r)\dd r\\
			&=\int_0^1\Big((\chi'r^{2\kappa_{1}+1})'a_1+(\chi''r+\chi'(2\kappa_{1}+1))b_1\Big)\dd r.
		\end{align*}
		An integration by parts and the choice of $\chi$ then finally give rise to the result
		\begin{align*}
			0&=[\chi'r^{2\kappa_{1}+1}a_1]_{r=0}^1+[(\chi'r+\chi(2\kappa_{1}+1))b_1]_{r=0}^1-\int_0^1\chi'b_1\dd r=2\kappa_{1}b_1.
		\end{align*}
		This means that also $b_1$ is zero.
		
		4) We finally show that all numbers $a_k$ are zero.
		Analogously to the proof of Theorem~5.1 in \cite{Kellogg71}, we employ the mapping $\check{u}_k(r,\varphi):=\xi(r)\psi_{k}(\varphi)$ for $k\in\NN$.
		Here, $\xi:=[0,1]\to\RR$ is smooth with $\supp\xi\sub(0,1]$, $\xi(1)=0$ and $\xi'(1)=1$.
		Then $\check{u}_k$ is an element of $\check{M}_{\inn}$.
		As in (\ref{eq 3.7.1}), we obtain
		\begin{align*}
			0=\int_0^1(r\xi''(r)+\xi'(r)-r^{-1}\kappa_{k}^2\xi(r))a_kr^{\kappa_{k}}\dd r.
		\end{align*} 
		Using the boundary conditions and the location of the support of $\xi$ in an integration by parts, the desired result
		\begin{align*}
			0&=[\xi'r^{\kappa_{k}+1}]_{r=0}^1a_k-\int_0^1\Big(\xi'\kappa_ka_kr^{\kappa_{k}}+\kappa_k^2a_k\xi r^{\kappa_k-1}\Big)\dd r=a_k-[\xi r^{\kappa_{k}}]_{r=0}^1\kappa_ka_k=a_k
		\end{align*}
		follows.
	\end{proof}
\end{proposition}

Combining Lemma~\ref{lemma 3.4} and Proposition~\ref{proposition 3.5}, there are bounded projections from $\cD(L_{\inn})$ onto the closed subspaces $\overline{\check{M}_{\inn}}$ and $\check{N}_{\inn}$, respectively.
The subspace $\check{N}_{\inn}$ plays here the role of the singular part of functions in the domain $\cD(L_{\inn})$.
This decomposition is one of the key tools in the regularity analysis of the electric field.


\section{Inhomogeneous elliptic transmission problems}\label{section 4}

The next two lemmas deal with elliptic transmission problems for functions whose normal derivatives have prescribed discontinuities across the interfaces.	
The first lemma is in particular useful to investigate the regularity of the magnetic field, while the second one is used for the electric field.
Although both lemmas are well known to experts in the field, see Section~4 in \cite{CoDaNi} for instance, we provide the corresponding proof at least for the first statement to keep the presentation self-contained.
Note that we use the formulation $H^{\kappa}(\Fint)$, $\kappa\in(0,1/2]$, for the space of all functions on the union $\cup_{\cF\in\Fint}\cF$ that belong on each interface $\cF$ to $H^{\kappa}(\cF)$.

\begin{lemma}\label{lemma 4.1}
	Let $\Gamma^*$ be a union of some of the exterior face pairs $\Gamma_s$, $s\in\{1,2,3\}$, and let $\mu$ satisfy (\ref{assumptions parameters}).
	Let additionally $\kappa\in(1/6,1/2)$, $f\in L^2(Q)$, and $g\in H^{\kappa}(\Fint)$.
	There is a unique function $u\in V:=\{\varphi\in PH^1_{\Gamma^*}(Q)\ |\ [[\mu\varphi]]_{\cF}=0\text{ for }\cF\in\Fint\}$ with
	\begin{align}
		\sum_{i=1}^N\int_{Q_i}\mu^{(i)}\nabla u^{(i)}\cdot\nabla\varphi^{(i)}\dd x=&\int_Q\mu f\varphi\dd x+\sum_{\cF\in\Fint}\int_{\cF}\mu g\varphi\dd \sigma,\qquad \varphi\in V.\label{eq 4.1}
	\end{align}
	The mapping $u$ even belongs to $PH^{3/2+\kappa}(Q)$ with
	\begin{align*}
		\sum_{i=1}^N\lVert u^{(i)}\rVert_{H^{3/2+\kappa}(Q_i)}^2\leq C\big(\lVert f\rVert_{L^2(Q)}^2+\sum_{\cF\in\Fint}\lVert g\rVert_{H^{\kappa}(\cF)}^2\big),
	\end{align*}
	involving a uniform constant $C=C(\mu,\kappa,Q)>0$.
	\begin{proof}
		1) We focus on the case $\Gamma^*=\Gamma_2$ as all remaining configurations can be treated with similar arguments. 
		The Lax-Milgram lemma and the trace theorem for $H^1$-functions provide a unique function $u\in V$ satisfying (\ref{eq 4.1}).
		We investigate the regularity of $u$ in the following. 
		
		2) Consider the elliptic transmission problem
		\begin{align}
			\begin{aligned}
				-\mu^{(i)}\Delta\psi^{(i)}&=0 &&\text{on }Q_i\text{ for }i\in\{1,\dots,N\},\\
				\psi&=0 &&\text{on }\Gamma_2,\\
				\nabla\psi\cdot\nu&=0 &&\text{on }\partial Q\setminus\Gamma_2,\\
				\llbracket\mu\psi\rrbracket_{\cF}&=0,\quad \llbracket\nabla\psi\cdot\nu_{\cF}\rrbracket_{\cF}=-\phi &&\text{on }\cF\in\Fint,
			\end{aligned}\label{eq 4.2} 
		\end{align} 
		involving a mapping $\phi$ on $\Fint$.
		
		Set $H^{-\delta}(\cF):=(H^{\delta}(\cF))^*$ for $\cF\in\Fint$, $\delta\in(0,1/2]$, and define $H^{-\delta}(\Fint)$ as the dual space of $H^{\delta}(\Fint)$.
		Note that $H^{-\delta}(\Fint)$ is isomorphic to the product $\prod_{\cF\in\Fint}H^{-\delta}(\cF)$.
		
		We first assume that $\phi\in H^{-1/2}(\Fint)$. 
		System (\ref{eq 4.2}) then corresponds to the formula
		\begin{align}
			\sum_{i=1}^N\int_{Q_i}\mu^{(i)}(\nabla\psi^{(i)})\cdot(\nabla\varphi^{(i)})\dd x&=\sum_{\cF\in\Fint}\langle\phi,\mu\varphi\rangle_{H^{-1/2}(\cF)\times H^{1/2}(\cF)},\quad \varphi\in V.\label{eq 4.3}
		\end{align}
		The Lax-Milgram lemma provides a unique weak solution $\psi=\psi(\phi)\in V$ of (\ref{eq 4.2}) respectively (\ref{eq 4.3}).
		We furthermore conclude the inequality
		\begin{align}
			\lVert\psi\rVert_{PH^1(Q)}\leq C\lVert\phi\rVert_{H^{-1/2}(\Fint)}.\label{eq 4.4}
		\end{align}
		
		3) Let $H^{1/2}_{00}(\Fint)$ be the space of all functions on $\Fint$, whose restrictions to all interfaces $\cF\in\Fint$ belong to $H^{1/2}_{00}(\cF)$.
		Let now $\phi\in H^{1/2}_{00}(\Fint)$.
		Using the extension results of Propositions~2.2 and 2.3 in \cite{AssousCiarlet97} on every interface, there is a function $\hat{\psi}\in PH^2(Q)$ with
		\begin{align}
			\begin{aligned}
				\llbracket\nabla\hat{\psi}\cdot\nu_{\cF}\rrbracket_{\cF}&=-\phi &&\text{on }\cF\in\Fint,\\
				\hat{\psi}&=0 &&\text{on }\cF\in\Fext\cup\Fint,\\
				\nabla\hat{\psi}^{(i)}\cdot\nu&=0 &&\text{on }\partial Q,\\
				\lVert\hat{\psi}\rVert_{PH^2(Q)}&\leq C\norm{\phi}_{H_{00}^{1/2}(\Fint)}. &&
			\end{aligned}\label{eq 4.5}
		\end{align}
		Proposition~3.1 in \cite{Zeru22b} further provides a unique function $\tilde{\psi}\in PH^2(Q)\cap V$ with $\Delta\tilde{\psi}^{(i)}=\Delta\hat{\psi}^{(i)}$ for $i\in\{1,\dots,N\}$ and $\llbracket\nabla\tilde{\psi}\cdot\nu_{\cF}\rrbracket_{\cF}=0$ for $\cF\in\Fint$.
		This means that the mapping $\psi:=\hat{\psi}-\tilde{\psi}\in PH^2(Q)$ solves (\ref{eq 4.2}) in strong form.
		Combining Proposition~3.1 in \cite{Zeru22b} with (\ref{eq 4.5}), we moreover arrive at the inequality
		\begin{align}
			\norm{\psi}_{PH^{2}(Q)}\leq C\norm{\phi}_{H_{00}^{1/2}(\Fint)}.\label{eq 4.6}
		\end{align}
		
		4) Let $\phi=g$. Remark~12.6 in Chapter~1 of \cite{LionsMagenes} leads to the identity $H^{\kappa}(\Fint)=(H^{-1/2}(\Fint),H^{1/2}_{00}(\Fint))_{\kappa+1/2,2}$.
		Interpolating now between (\ref{eq 4.4}) and (\ref{eq 4.6}), we thus infer that (\ref{eq 4.2}) has a unique solution $\psi=\psi(g)$ in the interpolation space
		\begin{align*}
			(V,PH^2(Q)\cap V)_{\kappa+1/2,2}&\hookrightarrow(PH^1(Q),PH^2(Q))_{\kappa+1/2,2}\cap V=PH^{3/2+\kappa}(Q)\cap V,
		\end{align*}
		satisfying the inequality
		\begin{align*}
			\lVert\psi\rVert_{PH^{3/2+\kappa}(Q)}\leq C\lVert g\rVert_{H^{\kappa}(\Fint)}.
		\end{align*}
	
		5) Proposition~3.1 in \cite{Zeru22b} provides a unique function $\tilde{u}\in PH^2(Q)\cap V$ with $\Delta\tilde{u}^{(i)}=-f^{(i)}$ and $\llbracket\nabla\tilde{u}\cdot\nu_{\cF}\rrbracket_{\cF}=0$ for $\cF\in\Fint$.
		It can be estimated by $\lVert \tilde{u}\rVert_{PH^2(Q)}\leq C\lVert f\rVert_{L^2(Q)}$. 
		Altogether, the function $u:=\tilde{u}+\psi\in PH^{3/2+\kappa}(Q)$ is the unique solution of (\ref{eq 4.1}), and it satisfies the asserted inequality.
	\end{proof}
\end{lemma}

Recall for the next statement that $\tFint$ denotes the collection of the interfaces between the coarse partition $\tilde{Q}_1,\dots,\tilde{Q}_L$, see Section~\ref{section1}.
We employ the spaces 
\begin{align*}
	H^1_{\Gamma^*}(\cF)&:=\{u\in H^1(\cF)\ |\ u=0\text{ on }\Gamma^*\cap\cF\},\\
	H^{1/2}_{\Gamma^*}(\cF)&:=(L^2(\cF),H^{1}_{\Gamma^*}(\cF))_{1/2,2},\quad \cF\in\tFint,\\
	H^{1/2}_{\Gamma^*}(\tFint)&:=\{v\in L^2(\tFint)\ |\ v|_{\cF}\in H^{1/2}_{\Gamma^*}(\cF)\text{ for }\cF\in\tFint\},
\end{align*}
for a nonempty union $\Gamma^*$ of the face pairs $\Gamma_2$ and $\Gamma_3$. 

As the lemma is a straightforward generalization of Lemmas~8.13 and 8.14 in \cite{ZerullaDiss}, we omit its proof.

\begin{lemma}\label{lemma 4.2}
	Let $\Gamma^*$ be a nonempty union of the face pairs $\Gamma_2$ and $\Gamma_3$.
	Let further $g\in H^{1/2}_{\Gamma^*}(\tFint)$, and $f\in L^2(Q)$.
	There is a unique function $u\in H^1_{\Gamma^*}(Q)$ with
	\begin{align}
		\int_Q\nabla u\cdot\nabla\varphi\dd x=\int_Qf\varphi\dd x+\sum_{\cF\in\tFint}\int_{\cF}g\varphi\dd \sigma,\qquad \varphi\in H^1_{\Gamma^*}(Q).\label{id1 second lemma derivative interface single jump}
	\end{align}
	The mapping $u$ is even an element of $PH^2(Q)$ with
	\begin{align*}
		\lVert u\rVert_{PH^2(Q)}\leq C(\lVert f\rVert_{L^2(Q)}+\sum_{\cF\in\tFint}\lVert g\rVert_{H^{1/2}_{\Gamma^*}(\cF)}),
	\end{align*}
	for a uniform constant $C=C(Q)>0$.
\end{lemma}


\section{Regularity analysis for the Maxwell system}\label{section 5}

This section is devoted to a detailed regularity analysis for the Maxwell system (\ref{Maxwell system}).
To establish the desired result in Corollary~\ref{corollary 6.17} and Remark~\ref{remark 6.18}, we derive an embedding statement for the space $X_2$ from (\ref{def X2}) in Proposition~\ref{proposition 6.15}, and employ semigroup theory on $X_2$.

The proof of Proposition~\ref{proposition 6.15} is structured into Sections~\ref{section 5.1}--\ref{section 5.4}.
In Section~\ref{section 5.1}, we use Theorem~7.1 from \cite{CoDaNi} to derive a first regularity statement for the electric field component of functions in $X_2$.
In Sections~\ref{section 5.2} and \ref{section 5.3}, we then separately analyze the electric and magnetic field components of a vector in $X_2$.
Here our findings from Sections~\ref{section3}-\ref{section 4} and \cite{Zeru22b} come into play.

\subsection{A first regularity statement}\label{section 5.1}\

This subsection provides a useful regulaity result for the electric field component of functions in the space $X_2$ from (\ref{def X2}).

\begin{lemma}\label{lemma 6.1}
	Let $\eps$ and $\mu$ satisfy (\ref{assumptions parameters}) and $(\mE,\mH)\in X_2$.
	The functions $\Delta\mE^{(i)}$ and $\Delta\mH^{(i)}$ belong to $L^2(Q_i)^3$ for $i\in\{1,\dots,N\}$.
	\begin{proof}
		Since the coefficients $\eps$ and $\mu$ are piecewise constant, the definition of $X_2$ implies that the function $\divv\mE^{(i)}$ belongs to $H^1(Q_i)$, and that the vector $\curl\curl\mE^{(i)}$ is contained in $L^2(Q_i)^3$. 
		We then calculate 
		\begin{align*}
			\curl\curl\mE^{(i)}=-\Delta\mE^{(i)}+\nabla\divv\mE^{(i)}
		\end{align*}
		in $H^{-1}(Q_i)$. As a result, $\Delta\mE^{(i)}$ belongs to $L^2(Q_i)^3$. The magnetic field component $\mH$ is treated similarly. 
	\end{proof}
\end{lemma}

\begin{remark}\label{remark 6.2}
	Lemma~\ref{lemma 6.1} implies an interior regularity result for fields $(\mE,\mH)\in X_2$ within the subcuboids $Q_1,\dots,Q_N$.
	Indeed, $\mE^{(i)}$ and $\mH^{(i)}$ belong to $H^2_{\text{loc}}(Q_i)^3$ for $i\in\{1,\dots,N\}$ by standard elliptic regularity theory, see for instance Section~6.3.1 in \cite{Evans}. 
	We need, however, regularity statements up to the boundary of every subcuboid for our error analysis.\hfill$\lozenge$ 
\end{remark}  

To establish the desired embedding statement for the space $X_2$ in Proposition~\ref{proposition 6.15}, we employ Theorem~7.1 of \cite{CoDaNi}.
As a preliminary step, we derive in the following a lower bound for the first eigenvalue of a Laplace-Beltrami operator on the lower hemisphere with transmission conditions and homogeneous Dirichlet boundary conditions.

We reuse the framework from Section~\ref{section3}.
In particular, $e_{\inn}=\{(0,0)\}\times[0,1]$ is an interior edge of four subcuboids $Q_{\inn,1},\dots,Q_{\inn,4}$ with $\eps$ satisfying (\ref{eq 3.2}).
Let $\cM:=(0,0,0)$ be one of the intersection points of $e_{\inn}$ with $\partial Q$.
After scaling, we can assume that the ball $B(\cM,1)$ with radius $1$ touches no other interior edge of $Q$.
Using the adjacent subcuboids $Q_{\inn,1},\dots,Q_{\inn,4}$ to $e_{\inn}$, we introduce the spherical regions
\begin{align*}
	G_{\inn}:=\partial B(\cM,1)\cap Q,\qquad G_{\inn,i}:=\partial B(\cM,1)\cap Q_{\inn,i},\qquad i\in\{1,\dots,4\}.
\end{align*} 

We also employ the piecewise constant representative $\eps_{\inn}$ of $\eps$ on $G_{\inn}$, and the notion of restrictions of functions and piecewise regularity is transferred to the partition of $G_{\inn}$ into $G_{\inn,1},\dots,G_{\inn,4}$.
We then study the Laplace-Beltrami operator
\begin{align*}
	L_{G}\psi&:=\tfrac{1}{\eps_{\inn}}\divv(\eps_{\inn}\nabla\psi),\\
	\psi\in\cD(L_{G})&:=\{\psi\in H^1_0(G_{\inn})\ |\ \divv(\eps_{\inn}\nabla \psi)\in L^2(G_{\inn})\},
\end{align*}
on $G_{\inn}$.
We note that $-L_{G}$ is positive definite, has a compact resolvent, and is selfadjoint on $L^2(G_{\inn})$ with respect to the $L^2$-inner product on $G_{\inn}$ with weight $\eps_{\inn}$.

We next determine a lower bound for the first eigenvalue of $-L_{G}$.

\begin{lemma}\label{lemma 5.3}
	Let $\eps_{\inn}$ satisfy (\ref{eq 3.2}).
	The first eigenvalue of $-L_{G}$ is greater than or equal to $2$.
	\begin{proof}
		Let $\psi\in\cD(-L_{G})^{1/2}=H_0^1(G_{\inn})$.
		We then extend $\psi^{(1)}$ in the following way to a function $\tilde{\psi}$ on the unit sphere $S^2$:
		Reflect $\psi^{(1)}$ at the interface $\partial G_{\inn,1}\cap\partial G_{\inn,4}$ to $G_{\inn,4}$. 
		The resulting function is afterwards reflected at the plane $\{x_2=0\}$ to become a function $\tilde{\psi}$ in $H^1_0(G_{\inn})$.
		We finally set
		\begin{align*}
			\tilde{\psi}(x_1,x_2,x_3):=-\tilde{\psi}(x_1,x_2,-x_3),\qquad (x_1,x_2,x_3)\in S^2\cap\{x_3<0\}.
		\end{align*}
		
		By construction, $\tilde{\psi}$ then belongs to $H^1(S^2)$, and it has zero mean.
		As the first positive eigenvalue of the standard Laplace-Beltrami operator on $S^2$ equals $2$, see Section~II.4 in \cite{IsaacChavel} for instance, we infer with the Poincaré inequality
		\begin{align*}
			\lVert\eps_{\inn}^{(1)}\nabla\tilde{\psi}\rVert_{L^2(S^2)}^2\geq 2\lVert\eps_{\inn}^{(1)}\tilde{\psi}\rVert_{L^2(S^2)}^2,
		\end{align*}
	see Section~D.II in Chapter III of \cite{BergerGauduchonMazet71}.
		In view of the definition of $\tilde{\psi}$, we then infer
		\begin{align*}
			\lVert\eps_{\inn}^{(1)}\nabla\psi^{(1)}\rVert_{L^2(G_{\inn,1})}^2\geq 2\lVert\eps_{\inn}^{(1)}\psi^{(1)}\rVert_{L^2(G_{\inn,1})}^2.
		\end{align*}
		Repeating the same argument for the restrictions of $\psi$ to $G_{\inn,2},G_{\inn,3},G_{\inn,4}$, we conclude
		\begin{align*}
			\lVert\eps_{\inn}\nabla\psi\rVert_{L^2(G_{\inn})}^2\geq 2\lVert\eps_{\inn}\psi\rVert_{L^2(G_{\inn})}^2.
		\end{align*}
		The Rayleigh quotient of $-L_{G}$ can now be estimated from below by
		\begin{align*}
			\frac{\int_{G_{\inn}}\eps_{\inn}(-L_{G}\phi)\phi\dd x}{\int_{G_{\inn}}\eps_{\inn}\phi^2\dd x}\geq 2,\qquad \phi\in \cD(-L_{G})\setminus\{0\}.
		\end{align*}
		The Courant-Fischer theorem now implies the asserted estimate.
	\end{proof}
\end{lemma}

We next transform our problem into the setting of \cite{CoDaNi} to apply Theorem~7.1 from there.
This provides a crucial regularity statement for the electric field component.
For the proof, we recall that the arising classes of exterior faces and interior faces are defined in Section~\ref{section1}.

\begin{lemma}\label{lemma 6.3}
	Let $\eps$ and $\mu$ satisfy (\ref{assumptions parameters}), and let $(\mE,\mH)\in X_2$.
	Then $\mE$ belongs to $PH^{3/2}(Q)^3$.
	\begin{proof}
		1) We first manipulate the function $\mE$ so that it satisfies the conditions required in \cite{CoDaNi}.
		Let $\cF\in\Finteff\setminus\tFint$, and $\tilde{Q}_{i,l}$ with $i\in\{1,\dots,L\}$, $l\in\{1,\dots,K\}$ be a cuboid with face $\cF$ (the case $\cF\in\tFint$ is treated similarly). 
		By definition of $X_2$ in (\ref{def X2}), the jump $\llbracket\eps\mE\cdot\nu_{\cF}\rrbracket_{\cF}$ belongs to $H_{00}^{3/2}(\cF)$.
		Lemma~\ref{lemma 2.1} then provides an extension $\tilde{u}_{\cF}\in H^3(\tilde{Q}_{i,l})\cap H_0^1(\tilde{Q}_{i,l})$ with $\nabla \tilde{u}_{\cF}\cdot\nu_{\cF'}=0$ on every face $\cF'$ of $\tilde{Q}_{i,l}$ except $\cF$, and $\partial_{\nu} \tilde{u}_{\cF}=-\tfrac{1}{\eps|_{\tilde{Q}_{i,l}}}\llbracket\eps\mE\cdot\nu_{\cF}\rrbracket_{\cF}$.
		Furthermore, $\Delta\tilde{u}_{\cF}$ is an element of $H^1_{00}(\tilde{Q}_{i,l})$. 
		(Recall that this means that $\Delta\tilde{u}_{\cF}$ has a trace in $H^{1/2}_{00}$ on every face of $\tilde{Q}_{i,l}$.)
		We then extend $\tilde{u}_{\cF}$ trivially to the remainder of $Q$, and arrive at the relation $\llbracket\eps\nabla \tilde{u}_{\cF}\cdot\nu_{\cF}\rrbracket_{\cF}=\llbracket\eps\mE\cdot\nu_{\cF}\rrbracket_{\cF}$.
		In a similar way, we continue for the remaining effective interfaces, obtaining functions $\tilde{u}_{\cF'}$ for $\cF'\in\Finteff$. 
		Summation gives rise to a new function $u:=\sum_{\cF'\in\Finteff}\tilde{u}_{\cF'}$.
		
		By construction, $u$ belongs to $H^1(Q)\cap PH^3(Q)$, the function $\divv(\eps\mE-\eps\nabla u)$ is an element of $PH^1(Q)\cap L^2(Q)$, and $u$ has vanishing derivatives at the exterior faces of $Q$.
		Combining additionally the fact $\Delta u|_{\tilde{Q}_{i,l}}\in H^1_{00}(\tilde{Q}_{i,l})$ with the requirement $\divv(\eps\mE|_{\tilde{Q}_{i,l}})\in H^1_{00}(\tilde{Q}_{i,l})$, we  infer that $\divv(\eps\mE-\eps\nabla u)|_{\tilde{Q}_{i,l}}$ is an element of $H^1_{00}(\tilde{Q}_{i,l})$ for $i\in\{1,\dots,L\}$, $l\in\{0,\dots,K\}$.
		
		2) Denote by $\cF^{\text{eff}}_{\text{int},j}$ the effective interfaces with normal vector $e_j$, $j\in\{1,2\}$.
		The set of exterior faces $\cF_{\text{ext},k}$, $k\in\{1,2,3\}$, is defined analogously. 
		Let $\cF\in\cF^{\text{eff}}_{\text{int},1}\setminus\tFint$ and let $\tilde{Q}_{i,l}$ be a cuboid with face $\cF$.
		Applying Proposition~2.3 in \cite{AssousCiarlet97} and extending trivially by zero, there is a function $v_{\cF,1}\in H^2(\tilde{Q}_{i,l})\cap H_0^1(\tilde{Q}_{i,l})$ with $\llbracket\eps\nabla v_{\cF,1}\cdot\nu_{\cF}\rrbracket_{\cF}=\llbracket\divv(\eps\mE-\eps\nabla u)\rrbracket_{\cF}$ and $\nabla v_{\cF,1}\cdot\nu_{\tilde{\cF}}=0$ for every other face $\tilde{\cF}$ of $\tilde{Q}_{i,l}$.
		We repeat the construction for every other interface $\tilde{\cF}\in\cF^{\text{eff}}_{\text{int},j}$, receiving functions $v_{\tilde{\cF},j}$, $j\in\{1,2\}$.
		At each exterior face $\tilde{\cF}\in\cF_{\text{ext},k}$, we extend the trace $\divv(\eps\mE-\eps\nabla u)|_{\tilde{\cF}}$, obtaining a function $v_{\tilde{\cF},k}$, $k\in\{1,2,3\}$.
		(For the exterior faces and faces in $\tFint$, we additionally use a cut-off argument to localize near the respective faces.) 
		Altogether, we obtain functions 
		\begin{align}
		v_j:=\sum_{\cF\in\cF^{\text{eff}}_{\text{int},j}\cup\cF_{\text{ext},j}}v_{\cF,j},\quad j\in\{1,2\},\quad v_3:=	\sum_{\cF\in\cF_{\text{ext},3}}v_{\cF,3},
		\end{align}
		that belong to $PH^2(Q)\cap H^1_0(Q)$.
		We set $v:=(v_j)_{j=1}^3$.
		
		By construction, the function $\tilde{\mE}:=\mE-\nabla u-v$ then belongs to $H_0(\curl,Q)$, and $\divv(\eps\tilde{\mE})$ is an element of $H^1_0(Q)$.
		Combining the identity $\curl\tilde{\mE}=\curl(\mE-v)$, the interface conditions of $v$ and Lemma~2.1 in \cite{EiSc18}, we further conclude $\curl\tilde{\mE}\times\nu_{\cF}=\curl\mE\times\nu_{\cF}$ for every interface $\cF\in\Fint$.
		This means that the vector $\tfrac{1}{\mu}\curl\tilde{\mE}$ belongs to $H(\curl,Q)$.
		
		3) We next show that the function $\tilde{\mE}$ fits into the setting of \cite{CoDaNi}.
		Let $\omega$ be a positive real number, $\imath$ be the imaginary unit, and set $\tilde{\mH}:=-\tfrac{1}{\imath\omega\mu}\curl\tilde{\mE}$.
		The vector $\tilde{\mH}$ then belongs to $H(\curl,Q)$.
		We also introduce the function $\tilde{\mJ}:=-\imath\omega\eps\tilde{\mE}+\curl\tilde{\mH}$.
		By definition, $(\tilde{\mE},\tilde{\mH})$ then satisfies the time-harmonic Maxwell equations with parameter $\omega$ and inhomogeneity $\tilde{J}$.
		Using the facts $\tfrac{1}{\mu}\curl\tilde{\mE}\in H(\curl,Q)$ and $\divv(\eps\tilde{\mE})\in H^1_0(Q)$ in an integration by parts, we further obtain that $\tilde{\mE}$ satisfies the variational problem (1.5) in \cite{CoDaNi} with $f\in L^2(Q)$.
		
		4) As a consequence of part 3), Theorem~7.1 in \cite{CoDaNi} applies to $\tilde{\mE}$.
		Hence, it remains to estimate the numbers $\sigma_{\eps}^{\text{Dir}}$ and $\sigma_{\mu}^{\text{Neu}}$ from p.~646 in \cite{CoDaNi}.
		In the following, we use the notation from pp.~645 in \cite{CoDaNi}.
		For an edge $e$ where $\eps$ is constant, a reflection argument, Satz~32.2 in \cite{Triebel72} and estimate~(1) in \cite{Lorch93} imply the lower estimate $\lambda_{\eps,e}^{\text{Dir}}>\sqrt{5}$.
		To bound $\lambda_{\eps,e}^{\text{Dir}}$ for an exterior edge $e$ where $\eps$ has a discontinuity, it suffices to consider the operator
		\begin{align*}
			\tilde{L}\psi&:=\tfrac{1}{\eta}\divv(\eta\nabla\psi),\quad\cD(\tilde{L}):=\{\psi\in H^1_0(D)\ |\ \divv(\eta\nabla\psi)\in L^2(Q)\},
		\end{align*} 
		where $\eta$ is a piecewise constant positive function on $\RR^2$ that changes value at $\{x_1=0\}$.
		(Recall that $D$ is the unit disc.)
		Adapting the reasoning for Lemma~3.6 in \cite{Zeru22b}, we find that the spectrum of $-\tilde{L}$ consists of the eigenvalues $(\mu_k^{(\kappa_l)})^2$, $k\in\NN$, $l\in\NN_0$, with $0=\kappa_0^2<\kappa_1^2\dots\to\infty$ denoting the eigenvalues of system (3.12) in \cite{Zeru22b}, and $\mu_k^{(\nu)}$ denoting the $k$-th zero of the Bessel function of order $\nu$.
		Estimate~(1) in \cite{Lorch93} then implies $\lambda_{\eps,e}^{\text{Dir}}>\sqrt{5}$.
		If $e$ is an interior edge where $\eps$ satisfies (\ref{eq 3.2}), Lemma~3.6 in \cite{Zeru22b}, Remark~\ref{remark previous work}, and relation~(1) in \cite{Lorch93} yield the same inequality $\lambda_{\eps,e}^{\text{Dir}}>\sqrt{5}$.
		
		5) We next estimate the number $\lambda_{\eps,c}^{\text{Dir}}$ for every corner $c$.
		Let $c$ first be a corner near which $\eps$ is constant.
		As the smallest positive eigenvalue of the Laplace-Beltrami operator on the unit sphere is $2$, see Section~II.4 in \cite{IsaacChavel}, a reflection argument shows $\lambda_{\eps,c}^{\text{Dir}}\geq1$.
		Let next $c$ be an external corner of an interface in $\tFint$, see Section~\ref{section1}.
		Then a modification of the reasoning in Lemma~\ref{lemma 5.3} yields again $\lambda_{\eps,c}^{\text{Dir}}\geq1.$
		Let finally $c$ be a corner near which $\eps$ satisfies (\ref{eq 3.2}).
		Then Lemma~\ref{lemma 5.3} yields also here $\lambda_{\eps,c}^{\text{Dir}}\geq1.$
		Altogether, we conclude $\sigma_{\eps}^{\text{Dir}}\geq3/2$.
		
		6) We finally bound the number $\sigma_{\mu}^{\text{Neu}}$ from below.
		Let $e$ first be an edge where $\mu$ is constant, and denote by $J_{\nu}$ the Bessel function of order $\nu$.
		A reflection argument, Satz~32.3 in \cite{Triebel72}, the identity $J_0'=-J_1$, estimate~(1) in \cite{Lorch93}, and the results in Section~15.3 of \cite{Watson} imply that $\lambda_{\mu,e}^{\text{Neu}}\geq\sqrt{3}$.
		Let next $e$ be an edge at which $\mu$ is discontinuous.
		Now Lemma~8.6 in \cite{ZerullaDiss} and the just mentioned references yield again $\lambda_{\mu,e}^{\text{Neu}}\geq\sqrt{3}$.
		
		Let $c$ be a corner at which $\mu$ is constant.
		The reasoning in part~5) then shows that $\lambda_{\mu,c}^{\text{Neu}}\geq1$.
		Let finally $c$ be an external corner of an interface in $\tFint$.
		A reflection argument and Lemma~8.9 in \cite{ZerullaDiss} then show that $\lambda_{\mu,c}^{\text{Neu}}>1/2$.
		We then conclude that $\sigma_{\mu}^{\text{Neu}}\geq1$.
		The asserted statement is now a consequence of Theorem~7.1 in \cite{CoDaNi}. 
	\end{proof}
\end{lemma}

\subsection{Analysis of the electric field component}\label{section 5.2}\

In this subsection, we study the regularity of the electric field component $\mE$ of a vector $(\mE,\mH)\in X_2$. 
The subsequent constructions mostly focus on the first component of $\mE$.
The second component can be treated similarly, due to the symmetry of the model problem. 
For the sake of a clear presentation, we elaborate the differences between the two components at the relevant steps.
The regularity of the third component is finally concluded.

Recall that $\tFint$ denotes the collection of all interfaces between the larger subcuboids $\tilde{Q}_1,\dots,\tilde{Q}_L$, see Section~\ref{section1}.
The next statement is directly obtained from Lemmas~9.12 and 9.13 in \cite{ZerullaDiss}, by means of a cut-off argument near an interface in $\tFint$.
We give a detailed account of a similar cut-off argument in the proof of Lemma~\ref{lemma 6.5}, whence we skip the next proof.

\begin{lemma}\label{lemma 6.4}
	Let $\eps,\mu$ satisfy (\ref{assumptions parameters}), let $(\mE,\mH)\in X_2$, and let $\cF\in\tFint$.
	There is an open set $\cO$ with $\cF\sub\cO$ and $\mE_1,\mE_2\in PH^2(\cO\cap Q)$. 
	Furthermore, the estimate
	\begin{align*}
		\lVert(\mE_1,\mE_2)\rVert_{PH^2(\cO\cap Q)}\leq C\Big(&\sum_{i=1}^N\big(\sum_{j=1}^3\big(\lVert\mE_j^{(i)}\rVert_{L^2(Q_i)}+\lVert\Delta\mE_j^{(i)}\rVert_{L^2(Q_i)}\big)+\lVert\curl\mE\rVert_{L^2(Q_i)}\big)\\
		&+\sum_{i=1}^L\sum_{l=0}^K\lVert\divv\mE|_{\tilde{Q}_{i,l}}\rVert_{H^1_{00}(\tilde{Q}_{i,l})}+\lVert\llbracket\eps\mE_1\rrbracket_{\cF}\rVert_{H_{00}^{3/2}(\cF)}\Big),
	\end{align*}
	is valid with a uniform constant $C=C(\eps,\mu,Q)>0$.
\end{lemma}

We next deduce that the first two electric field components are piecewise $H^2$-regular outside tubes around the interior edges of $Q$ where $\eps$ is discontinuous.
Although this statement might be well known to experts, we provide a proof for the sake of a self-contained presentation.

\begin{lemma}\label{lemma 6.5}
	Let $\eps,\mu$ satisfy (\ref{assumptions parameters}), and let $(\mE,\mH)\in X_2$.
	Let $\delta>0$ be smaller than half of the shortest edge of one of the subcuboids of $Q$. 
	Denote by $T(\delta)$ the union of cylinders with radius $\delta$ around the interior edges at which $\eps$ is discontinuous. 
	Then $\mE_1^{(i)},\mE_2^{(i)}$ belong to $H^2(Q_i\setminus T(\delta))$ for $i\in\{1,\dots,N\}$ with
	\begin{align*}
		\lVert\mE_m\rVert_{PH^2(Q\setminus T(\delta))}\leq C&\Big[\sum_{i=1}^N\big(\sum_{j=1}^3(\lVert\mE_j^{(i)}\rVert_{L^2(Q_i)}+\lVert\Delta\mE_j^{(i)}\rVert_{L^2(Q_i)})+\lVert\curl\mE\rVert_{L^2(Q_i)}\big)\\
		&+\sum_{i=1}^L\sum_{l=0}^K\lVert\divv\mE|_{\tilde{Q}_{i,l}}\rVert_{H^1_{00}(\tilde{Q}_{i,l})}+\sum_{\cF\in\Fint}\lVert\llbracket\eps\mE\cdot\nu_{\cF}\rrbracket_{\cF}\rVert_{H^{3/2}_{00}(\cF)}\Big],
	\end{align*}  
	$m\in\{1,2\}$, involving a uniform constant $C=C(\eps,\mu,\delta,Q)>0$.
	\begin{proof}
		1) Note that Lemma~\ref{lemma 6.4} already provides the asserted regularity and estimate in a neighborhood of all interfaces in $\tFint$.
		We first focus on the interfaces in $\Finteff\setminus\tFint$, and employ a cut-off argument.
		Let $e_1,\dots,e_S$ be the interior edges in $Q$ 
		near which $\eps$ is discontinuous. 
		We set $\cE:=\{e_1,\dots,e_S\}$.
		For $l\in\{1,\dots,S\}$, denote the distance function to $e_l$ by $d_l$.
		Additionally, we use a smooth function $\tilde{\chi}:[0,\infty)\to[0,1]$ with $\tilde{\chi}=1$ on $[0,\delta^2/36]$ and $\supp\tilde{\chi}\sub[0,\delta^2/16]$.
		Set
		\begin{align*}
			\hat{\chi}(x):=\prod_{l=1}^S(1-\tilde{\chi}(d^2_l(x))), \qquad x\in\overline{Q}.
		\end{align*}
		This mapping is smooth, vanishes in $T(\delta/6)$, and is equal to one outside of $T(\delta/4)$. In the following, we analyze the product $\hat{\chi}\mE$.
		
		2) Let $\cF\in\Finteff$ with $\nu_{\cF}=(1,0,0)$, and let $\mathring{Q}_1,\mathring{Q}_2\sub Q$ be two adjacent cuboids with interface $\cF$ and side length $\delta/6$ in $x_1$-direction.
		Without loss of generality, we can assume that $\cF=\{0\}\times[0,1]^2$.
		Set also 
		\begin{align*}
			\mathring{Q}:=\mathring{Q}_1\cup\mathring{Q}_2\cup\cF=(-\tfrac{\delta}{6},\tfrac{\delta}{6})\times[0,1]^2.	
		\end{align*}
		We additionally employ a smooth cut-off function $\mathring{\chi}:[-1,1]\to[0,1]$ with $\mathring{\chi}=1$ on $[-\delta/8,\delta/8]$, and $\supp\mathring{\chi}\sub[-\delta/7,\delta/7]$.
		In the following, we show that $g_j:=\mathring{\chi}(x_1)\hat{\chi}\mE_j$ is piecewise $H^2$-regular on $\mathring{Q}$ for $j\in\{1,2\}$.
		
		Let $j\in\{1,2,3\}$, $k\in\{0,1\}$.
		As a consequence of the product rule, we first infer that $\divv(\partial_j^k(\mathring{\chi}\hat{\chi})\mE)\in L^2(\mathring{Q}_i)$, $\curl(\partial_j^k(\mathring{\chi}\hat{\chi})\mE)\in L^2(\mathring{Q})$, $\llbracket\eps\partial_j^k(\mathring{\chi}\hat{\chi})\mE\cdot\nu_{\cF}\rrbracket_{\cF}\in H^{3/2}_{00}(\cF)$, and that $\partial_j^k(\mathring{\chi}\hat{\chi})\mE\times\nu_{\partial\mathring{Q}}=0$ on $\partial\mathring{Q}$.
		By Proposition~9.8 in \cite{ZerullaDiss}, $\partial_j^k(\mathring{\chi}\hat{\chi})\mE\in PH^1(\mathring{Q})$ with
		\begin{align}
			\lVert\partial_j^k(\mathring{\chi}\hat{\chi})\mE\rVert_{PH^1(\mathring{Q})}\leq C\big(\lVert\mE\rVert_{L^2(\mathring{Q})}+\lVert\curl\mE\rVert_{L^2(\mathring{Q})}+\sum_{i=1}^2\lVert\divv\mE\rVert_{L^2(\mathring{Q}_i)}\big).\label{eq 6.1}
		\end{align}
		
		Combining the product rule with Lemma~\ref{lemma 6.1}, we infer that $\Delta g_j|_{\mathring{Q}_i}$ belongs to $L^2(\mathring{Q}_i)$ for $i,j\in\{1,2\}$, and that $\eps g_1-\eps\psi$ is contained in $H^1_0(\mathring{Q})$ with an appropriate extension $\psi$, see Lemma~\ref{lemma 2.1}.
		Using the identity 
		\begin{align*}
			\partial_1g_1=\partial_1(\mathring{\chi}\hat{\chi})\mE_1+\mathring{\chi}\hat{\chi}\divv\mE-\mathring{\chi}\hat{\chi}(\partial_2\mE_2+\partial_3\mE_3),
		\end{align*}
		the definition of $X_2$ in (\ref{def X2}) and the fact $\nabla\mE_2,\nabla\mE_3\in H(\curl,\mathring{Q})$, we conclude the relation
		\begin{align*}
			\llbracket\eps(\partial_1 g_1-\partial_1\psi)\varphi\rrbracket_{\cF}=\llbracket\mathring{\chi}\hat{\chi}\divv\mE-\partial_1\psi\rrbracket_{\cF}(\eps\varphi)
		\end{align*}
		for all $\varphi\in\{\varphi\in PH^1(\mathring{Q})\ |\ \eps\varphi\in H^1_0(\mathring{Q})\}$.
		Furthermore, $\llbracket\mathring{\chi}\hat{\chi}\divv(\eps\mE)-\partial_1\psi\rrbracket_{\cF}\in H^{1/2}_{00}(\cF)$ by the trace theorem.
		Consequently, Lemma~8.14 in \cite{ZerullaDiss} shows that $g_1$ is piecewise $H^2$-regular on $\mathring{Q}$ and that it satisfies
		\begin{align*}
			\lVert g_1\rVert_{PH^2(\mathring{Q})}\leq C\Big[&\sum_{j=1}^2\big(\lVert\mE_1\rVert_{L^2(\mathring{Q}_j)}+\lVert\Delta\mE_1\rVert_{L^2(\mathring{Q}_j)}+\lVert\curl\mE\rVert_{L^2(\mathring{Q}_j)}\\
			&+\lVert\divv\mE\rVert_{H^1_{00}(\mathring{Q}_{j})}\big)+\lVert\llbracket\eps\mE_1\rrbracket\rVert_{H^{3/2}_{00}(\cF)}\Big],
		\end{align*}
		where we also use (\ref{eq 6.1}) and Lemma~\ref{lemma 2.1}.
		
		3) We next study $g_2=\mathring{\chi}(x_1)\hat{\chi}\mE_2$ on $\mathring{Q}$, and restrict ourselves to the case that $\cF$ does not touch the boundary faces of $Q$ in $\Gamma_2$. 
		(The other case is obtained by a straightforward modification.)
		The mapping $g_2$ belongs to $H^1(\mathring{Q})$, and vanishes on the faces with normal vector $e_j$, $j\in\{1,3\}$.
		We further infer the formula
		\begin{align*}
			\llbracket(\partial_1g_2)\tilde{\varphi}\rrbracket_{\cF}=\llbracket\mathring{\chi}\hat{\chi}(\curl\mE-\partial_2\mE_1)\rrbracket_{\cF}\tilde{\varphi},\qquad\tilde{\varphi}\in H^1(\mathring{Q}),
		\end{align*}
		and that $\llbracket\mathring{\chi}\hat{\chi}(\curl\mE-\partial_2\mE_1)\rrbracket_{\cF}\in H^{1/2}_{00}(\cF)$.
		(Note that $\curl\mE\in H^1(Q)^3$ by Proposition~4.6 in \cite{Zeru22b}.)
		Let $\mathring{\Gamma}_2$ be the union of faces of $\mathring{Q}$ with normal vector $e_2$.
		We note that $\partial_2g_2=0$ on $\mathring{\Gamma}_2$.
		By Lemma~8.13 in \cite{ZerullaDiss}, $g_2\in PH^2(\mathring{Q})$ with
		\begin{align*}
			\lVert g_2\rVert_{PH^2(\mathring{Q})}\leq C&\sum_{j=1}^2\big(\lVert\mE_2\rVert_{L^2(\mathring{Q}_j)}+\lVert\Delta\mE_2\rVert_{L^2(\mathring{Q}_j)}+\lVert\curl\mE\rVert_{L^2(\mathring{Q}_j)}\\
			&+\lVert\divv\mE\rVert_{H^1_{00}(\mathring{Q}_{j})}+\lVert\llbracket\eps\mE_1\rrbracket_{\cF}\rVert_{H^{3/2}_{00}(\cF)}\big),
		\end{align*}
	where we also use (\ref{eq 6.1}), Proposition~4.6 in \cite{Zeru22b}, and Remark~\ref{remark previous work}.
	
	4) Due to symmetry, the results of 2) and 3) transfer directly to all interfaces in $\Finteff$ with normal vector $\nu_{\cF}=(0,1,0)$.
	Hence, it remains to consider the first two components of $\mE$ on $Q\setminus\cS_{\delta/10}$, with $\cS_{\delta/10}$ being an open neighborhood around all interfaces in $Q$ with distance at most $\delta/10$ to an interface.
	
	Let $m\in\{1,2\}$, $i\in\{1,\dots,N\}$, and $\cS_{\delta/16}$ be an open neighborhood around the interfaces with maximal distance $\delta/16$. 
	Let additionally $\overline{\chi}:\overline{Q}_i\to[0,1]$ be a smooth cut-off function with the following properties.
	The function is equal to $1$ on $\overline{Q}_i\setminus\cS_{\delta/10}$, and supported in $\overline{Q}_i\setminus\cS_{\delta/16}$.
	Similar reasoning as in parts~2) and 3) yields that the product $\overline{\chi}\mE_m^{(i)}$ fits into the framework of Lemma~3.1 in \cite{EiSc32}.
	Consequently, $\overline{\chi}\mE_m^{(i)}$ belongs to $H^2(Q_i)$ with
	\begin{align*}
		\lVert\overline{\chi}\mE_m^{(i)}\rVert_{H^2(Q_i)}\leq C\big(&\lVert\mE_m^{(i)}\rVert_{L^2(Q_i)}+\lVert\Delta\mE_m^{(i)}\rVert_{L^2(Q_i)}+\lVert\curl\mE^{(i)}\rVert_{L^2(Q_i)}\\
		&+\lVert\divv\mE^{(i)}\rVert_{H^1_{00}(Q_i)}\big).
	\end{align*}
	\end{proof}
\end{lemma}

We next analyze $\mE$ near interior edges where $\eps$ has a discontinuity. 
Only the second component of $\mE$ is analyzed in detail as the first one can be treated with similar arguments, due to symmetry.
In the following, we use the notation from Section~\ref{section3}.
In particular, we fix an interior edge $e$ at which $\eps$ has a discontinuity, and assume
\begin{align*}
	e=\{(0,0)\}\times[0,1].
\end{align*}
The adjacent subcuboids are denoted by $Q_{\inn,1},\dots,Q_{\inn,4}$, and the representative $\eps_{\inn}$ of $\eps$ is assumed to satisfy (\ref{eq 3.2}).

As in \cite{CoDa,CoDaNi}, we use a cylindrical coordinate system centered at $e$.
In particular, we employ the cylinder
\begin{align*}
	\cZ:=D\times[-1,2],
\end{align*}
around $e$ with radius $1$.
After scaling, we can assume that $\cZ$ touches no other interior edge.
We set
\begin{align*}
	\tilde{Q}_{\inn,i}:=\{(x_1,x_2,x_3)\ |\ (x_1,x_2,\tfrac{1}{2})\in Q_{\inn,i},\ x_3\in(-1,2)\},\qquad \cZ_i:=\cZ\cap \tilde{Q}_{\inn,i},
\end{align*} 
for $i\in\{1,\dots,4\}$,
and transfer the notion of restrictions of functions and piecewise regularity to this partition of $\cZ$.
The parameter $\eps_{\inn}$ is extended to $\cZ$ by even reflection at the top and bottom face of $Q$.
We assume that
\begin{align*}
	\cZ_i=D_{\inn,i}\times[-1,2],\qquad i\in\{1,\dots,4\},
\end{align*}
with the disk parts $D_{\inn,i}$ from (\ref{eq 3.3}).
The interfaces
\begin{align*}
	\cF_k^{\cZ}:=\overline{\cZ_k}\cap\overline{\cZ_{k+1}},\quad \cF_4^{\cZ}:=\overline{\cZ_1}\cap\overline{\cZ_4},\qquad k\in\{1,2,3\},
\end{align*}
are also employed.

Let $(\mE,\mH)\in X_2$, and denote the even reflection at the top and bottom face of $Q$ by $w^+$ for $w\in L^2(Q)$.
We extend the electric field component $\mE$ of $(\mE,\mH)\in X_2$ to the function
\begin{align}
	\tilde{\mE}&:=\mE &&\text{on }\cZ\cap Q,\label{eq 6.1.1}\\
	\tilde{\mE}&:=\begin{pmatrix}-\mE_1^+\\-\mE_2^+\\\mE_3^+\end{pmatrix} &&\text{on }\cZ\setminus Q.
\end{align}
Additionally, we employ smooth cut-off functions $\chi_{1,2},\chi_3:[0,\infty)\to[0,1]$ with $\chi_{1,2}=1$ on $[0,1/2]$, $\chi_3=1$ on $[0,2/3]$, and $\supp\chi_{1,2}\sub[0,\tfrac{5}{8}]$, $\supp\chi_3\sub[0,3/4]$.
We then study the function
\begin{align}
	v(x_1,x_2,x_3):=\chi_{1,2}(\lvert(x_1,x_2)\rvert)\chi_3(\lvert x_3-\tfrac{1}{2}\rvert)\tilde{\mE}_2(x_1,x_2,x_3)\label{eq 6.2}
\end{align}
for $(x_1,x_2,x_3)\in\cZ$. 

With Lemma~\ref{lemma 6.3} and the boundary condition $\mE\times\nu=0$ on $\partial Q$, $v$ is piecewise $H^1$-regular, satisfies homogeneous Dirichlet boundary conditions on the boundary $\partial\cZ$, and $\Delta v^{(i)}$ is an element of $L^2(\cZ_i)$.

Consider the Laplacian
\begin{align}
	(\Delta_{\cZ}u)^{(i)}&:=\Delta u^{(i)},\qquad i\in\{1,2,3,4\},\label{eq 6.3}\\
	\cD(\Delta_{\cZ})&:=\{\psi\in PH^1(\cZ)\ |\ \Delta\psi^{(i)}\in L^2(\cZ_i),\ \psi=0 \text{ on }\partial\cZ,\ \llbracket\psi\rrbracket_{\cF_k^{\cZ}}=0\nonumber\\
	&\hphantom{:=\{\psi\in PH^1_0(\cZ)\ |\ }=\llbracket\eps_{\inn}\nabla\psi\cdot\nu_{\cF_k^{\cZ}}\rrbracket_{\cF_k^{\cZ}},\ \llbracket\eps_{\inn}\psi\rrbracket_{\cF_l^{\cZ}}=0=\llbracket\nabla\psi\cdot\nu_{\cF_l^{\cZ}}\rrbracket_{\cF_l^{\cZ}}\nonumber\\
	&\hphantom{:=\{\psi\in PH^1_0(\cZ)\ |\ }\text{for }i\in\{1,\dots,4\},\ k\in\{1,3\},l\in\{2,4\}\},\nonumber
\end{align}
on $\cZ$. 
After reflecting and rotating, we can assume that $v$ satisfies the zero order transmission conditions on $\cF_1^{\cZ}$ and $\cF_3^{\cZ}$ in the domain $\cD(\Delta_{\cZ})$.
Note that $-\Delta_{\cZ}$ is positive definite and selfadjoint on $L^2(\cZ)$ with respect to the $L^2$-inner product with weight $\eps_{\inn}$.
Useful is also the formula
\begin{align*}
	\cD(-\Delta_{\cZ})^{1/2}=\{\phi\in PH^1(\cZ)\ |\ &\psi=0\text{ on }\partial \cZ,\ \llbracket\phi\rrbracket_{\cF_k^{\cZ}}=0,\ \llbracket\eps_{\inn}\phi\rrbracket_{\cF_l^{\cZ}}=0,\\ &k\in\{1,3\},l\in\{2,4\}\}.
\end{align*}

We first provide two technical preliminary lemmas.
The first one uses to approximate functions in $\cD(-\Delta_{\cZ})^{1/2}$ by regular functions.

\begin{lemma}\label{lemma 6.6}
	Let $\eps_{\inn}$ satisfy (\ref{eq 3.2}).
	The space
	\begin{align*}
		W:=\{&\psi\in PH^2(\cZ)\cap\cD(-\Delta_{\cZ})^{1/2}\ |\ \psi^{(i)}\text{ is smooth},\ \supp\psi\cap\partial\cZ=\emptyset,\ \partial_{\nu_{\cF}}\psi^{(i)}=0\\
		&\forall\ \text{faces }\cF\text{ of }\cZ_i,\ \varphi\text{ vanishes in a neighborhood of all edges of }\cZ_1,\dots,\cZ_4\}
	\end{align*}
	is dense in $\cD(-\Delta_{\cZ})^{1/2}$.
	\begin{proof}
		Let $\phi\in\cD(-\Delta_{\cZ})^{1/2}$, and let $\tilde{\chi}_m:[0,\infty)\to[0,1]$ be a smooth cut-off function with $\tilde{\chi}_m=1$ on $[0,1-1/m]$, $\supp\tilde{\chi}_m\sub[0,1-\tfrac{1}{2m}]$, and $\lVert\tilde{\chi}_m'\rVert_{\infty}\leq Cm$ with a uniform constant $C>0$ for $m\in\NN$.
		Consider then the mapping
		\begin{align*}
			\psi_m(x_1,x_2,x_3):=\tilde{\chi}_m(\lvert(x_1,x_2)\rvert)\phi(x_1,x_2,x_3),\qquad (x_1,x_2,x_3)\in\cZ.
		\end{align*}
		Taking into account that $\phi$ vanishes on $\partial \cZ$, a similar reasoning as in the proof of Lemma~2.1 in \cite{EiSc18} shows that $\psi_m^{(i)}\to\phi^{(i)}$ in $H^1(\cZ_i)$ as $m\to\infty$.
		
		Fix $m\geq3$ in the following. We extend $\psi_m$ trivially in $x_1-x_2$-direction to the cuboid $\tilde{Q}:=[-1,1]^2\times[-1,2]$, keeping its piecewise $H^1$-regularity and transmission behavior.
		Adapting Lemma~2.1 in \cite{Zeru22b} to $\tilde{Q}$ as well as to the present boundary and transmission conditions, there is a sequence $(\tilde{\psi}_l)_l$ of piecewise smooth functions converging on $\tilde{Q}_{\inn,i}$ to $\psi_m$ in $H^1$ and fulfilling the following properties.
		Each mapping $\tilde{\psi}_l$ vanishes on $\partial\tilde{Q}$ and in a neighborhood of all interior and exterior edges of $\tilde{Q}$ (meaning the intersections of $\partial\cZ_{\inn,i}$ with $\Fint$ for $i\in\{1,\dots,4\}$).
		Additionally, it obeys the transmission conditions in $\cD(-\Delta_{\cZ})^{1/2}$, and the normal derivative of $\tilde{\psi}_l$ vanishes in a neighborhood of each interface in $\tilde{Q}$.
		
		We then define
		\begin{align*}
			\phi_l(x_1,x_2,x_3):=\tilde{\chi}_{3m}(\lvert(x_1,x_2)\rvert)\tilde{\psi}_l(x_1,x_2,x_3),\qquad (x_1,x_2,x_3)\in\tilde{Q},\ l\in\NN.
		\end{align*}
		Each function $\phi_l$ belongs to $W$, and we deduce the convergence statement
		\begin{align*}
			\lVert\phi_l^{(i)}-\psi_m^{(i)}\rVert_{H^1(\cZ_i)}&=\lVert\tilde{\chi}_{3m}^{(i)}(\lvert(x_1,x_2)\rvert)(\tilde{\psi}_l^{(i)}-\psi_m^{(i)})\rVert_{H^1(\cZ_i)}\\
			&\leq\lVert\tilde{\chi}_{3m}^{(i)}\rVert_{W^{1,\infty}(\cZ_i)}\lVert\tilde{\psi}_l^{(i)}-\psi_m^{(i)}\rVert_{H^1(\cZ_i)}\to 0,\qquad l\to\infty,
		\end{align*}
		by using the construction of $\psi_m$ (in particular the fact $\tilde{\chi}_{3m}\tilde{\chi}_m=\tilde{\chi}_m$).
		As $(\psi_m)_m$ converges to $\phi$, we conclude the desired statement.
	\end{proof}
\end{lemma}

We next construct a regular function that extends the jumps of the second component of the electric field across the interfaces.
The statement is valid up to possible reduction of the radius of $\cZ$, which is implicitly assumed in the following.

\begin{lemma}\label{lemma 6.7}
	Let $\eps,\mu$ satisfy (\ref{assumptions parameters}), and let $\eps_{\inn}$ satisfy (\ref{eq 3.2}).
	Let $(\mE,\mH)\in X_2$, and define $v$ by (\ref{eq 6.2}).
	There is a function $\psi\in PH^2(\cZ)$ with $v+\psi\in\cD(-\Delta_{\cZ})^{1/2}$ and
	\begin{align*}
		\lVert\psi\rVert_{PH^2(\cZ)}\leq C\sum_{\cF\in\Finteff}\lVert\llbracket\eps\mE\cdot\nu_{\cF}\rrbracket_{\cF}\rVert_{H^{3/2}_{00}(\cF)},	
	\end{align*}
	for a uniform constant $C=C(\eps,\mu,Q)>0$.
	\begin{proof}
		It suffices to extend the jump $\llbracket\eps_{\inn}v\rrbracket_{\cF_4^{\cZ}}$ to the cylinder $\cZ$.
		As $\eps_{\inn}|_{\cZ_{4}}\not=\eps_{\inn}|_{\cZ_j}$ for $j\in\{1,2,3\}$, there is an effective interface $\cF$ with $\cF_4^{\cZ}\cap Q\sub\cF$.
		After adaption of the radius $\cZ$, there is a cuboid $\hat{Q}\sub Q$ with face $\cF$, $\cZ_4\cap Q\sub\hat{Q}$ and $\eps$ being constant on $\hat{Q}$.
		
		As $\llbracket\eps\mE\cdot\nu_{\cF}\rrbracket_{\cF}\in H_{00}^{3/2}(\cF)$, Lemma~\ref{lemma 2.1}  provides $\hat{\psi}\in H^3(\hat{Q})\cap H^1_0(\hat{Q})$ with 
		\begin{align*}
			\partial_{2}\hat{\psi}|_{\cF}=\frac{1}{\eps_{\inn}^{(4)}}\llbracket\eps\mE_2\rrbracket_{\cF},\quad
			\lVert\hat{\psi}\rVert_{H^3(\hat{Q})}\leq C\lVert\llbracket\eps\mE\cdot\nu_{\cF}\rrbracket_{\cF}\rVert_{H^{3/2}_{00}(\cF)}.
		\end{align*}
		On all other faces of $\hat{Q}$, the function $\hat{\psi}$ satisfies homogeneous Neumann boundary conditions.
		We then define $\tilde{\psi}:=\partial_2\hat{\psi}$ on $\hat{Q}$ and $\tilde{\psi}:=0$ on $Q\setminus\hat{Q}$.
		Lemma~2.1 from \cite{EiSc18} and the construction of $\hat{\psi}$ then imply that $\tilde{\psi}$ and its normal derivatives vanish on all interior and exterior faces of $Q$, except $\cF$.
		We then extend $\tilde{\psi}$ to $\cZ$ by
		\begin{align*}
			\check{\psi}:=\tilde{\psi}\quad\text{on }\cZ\cap Q,\qquad \check{\psi}:=-\tilde{\psi}^+\quad\text{on }\cZ\setminus Q,
		\end{align*}
		where $\tilde{\psi}^+$ denotes the even reflection of $\tilde{\psi}$ at the top and bottom face of $Q$.
		By definition, $\check{\psi}^{(i)}$ is $H^2$-regular on $\cZ_i$, and $\llbracket\eps_{\inn}\check{\psi}\rrbracket_{\cF_4^{\cZ}}=-\llbracket\eps_{\inn}\tilde{\mE}_2\rrbracket_{\cF_4^{\cZ}}$.
		 
		Recall next the smooth cut-off functions $\chi_{1,2}$ and $\chi_3$ from definition~(\ref{eq 6.2}). 
		We define the function
		\begin{align}
			\psi(x_1,x_2,x_3):=\chi_{1,2}(\lvert(x_1,x_2)\rvert)\chi_3(\lvert x_3-\tfrac{1}{2}\rvert)\check{\psi}(x_1,x_2,x_3),\quad (x_1,x_2,x_3)\in\cZ.\label{*const psi}
		\end{align}
		We note that the desired extension property $\llbracket\eps_{\inn}\psi\rVert_{\cF_4^{\cZ}}=-\llbracket\eps_{\inn} v\rrbracket_{\cF_4^{\cZ}}$ is valid, that the sum $v+\psi$ belongs to $\cD(-\Delta_{\cZ})^{1/2}$, and that $\psi$ satisfies the asserted estimate.
	\end{proof}
\end{lemma}

Using ideas from the proof of Lemma~3.7 in \cite{HoJaSc}, we next derive an appropriate elliptic transmission problem for the function $v$ on $\cZ$.

\begin{lemma}\label{lemma 6.8}
	Let $\eps,\mu$ satisfy (\ref{assumptions parameters}), and let $\eps_{\inn}$ satisfy (\ref{eq 3.2}). 
	Let also $(\mE,\mH)\in X_2$, $v$ be defined by (\ref{eq 6.2}), and take $\psi$ from Lemma~\ref{lemma 6.7}.
	There is a function $g\in \bigcap_{j\in\{3,4\}}H^{1/2}(\cF_j^{\cZ})$ with
	\begin{align*}
		\sum_{i=1}^4\int_{\cZ_i}\eps_{\inn}\nabla(v^{(i)}+\psi^{(i)})\cdot\nabla\phi^{(i)}\dd x&=-\sum_{i=1}^4\int_{\cZ_i}\eps_{\inn}\Delta(v^{(i)}+\psi^{(i)})\phi^{(i)}\dd x\\
		&\quad-\int_{\cF_3^{\cZ}}g\phi\dd\sigma-\int_{\cF_4^{\cZ}}\eps_{\inn}g\phi\dd\sigma
	\end{align*}
	for $\phi\in\cD(-\Delta_{\cZ})^{1/2}$.
	The mapping $g$ satisfies $g=\chi_{1,2}(\tfrac{2}{3}\lvert(x_1,x_2)\rvert)\chi_3(\tfrac{2}{3}\lvert x_3-\tfrac{1}{2}\rvert)g$, and can be estimated according to
	\begin{align*}
		\norm{g}_{H^{1/2}(\cF_j^{\cZ})}&\leq C\Big[\sum_{i=1}^N\big(\sum_{j=1}^2(\lVert\mE_j^{(i)}\rVert_{L^2(Q_i)}+\lVert\Delta\mE_j^{(i)}\rVert_{L^2(Q_i)})+\lVert\curl\mE\rVert_{H^1(Q_i)}\big)\\
		&\qquad\qquad+\sum_{i=1}^L\sum_{l=0}^K\lVert\divv\mE\rVert_{H^1_{00}(\tilde{Q}_{i,l})}+\sum_{\cF\in\Finteff}\lVert\llbracket\eps\mE\cdot\nu_{\cF}\rrbracket_{\cF}\rVert_{H^{3/2}_{00}(\cF)}\Big],
	\end{align*}
	with a uniform constant $C=C(\eps,\mu,Q)>0$.
	\begin{proof}
		1) We approximate $\cZ_i$ from the interior by means of the sets
		\begin{align*}
			\cZ_{i,n}&:=\{(r\cos\varphi,r\sin\varphi,z)\ |\ r\in[1/n,1-1/n],\ \varphi\in I_{\inn,i}^n,\ z\in[-1+1/n,2-1/n]\},\\
			I_{\inn,1}^n&:=[1/n,\pi/2-1/n],\quad I_{\inn,2}^n:=[\pi/2+1/n,\pi-1/n],\\
			I_{\inn,3}^n&:=[\pi+1/n,3/2\pi-1/n],\quad
			I_{\inn,4}^n:=[3/2\pi+1/n,2\pi-1/n],
		\end{align*}
		for $i\in\{1,\dots,4\}$ and $n\geq n_0\in\NN$, compare (\ref{eq 3.3}).
		The plane face of $\cZ_{i,n}$ with normal vector $e_j=(\delta_{jk})_k$ is denoted by $\cF_{i,n}^{j}$ for $j\in\{1,2\}$.
		The reasoning in Lemma~\ref{lemma 6.1} and Remark~\ref{remark 6.2} then implies that $u:=v+\psi$ is $H^2$-regular on $\cZ_{i,n}$.
		
		2) Recall the space $W$ from Lemma~\ref{lemma 6.6}.
		Let $\phi\in W$.
		We write in the following $\divv\tilde{\mE}$ and $\divv(\eps_{\inn}\tilde{\mE})$ for convenience, meaning the functions that are defined piecewise on $\cZ$.
		Using that $\nabla u^{(i)}$ vanishes near the boundary of $\cZ$, an integration by parts leads to the relations
		\begin{align}
			\sum_{i=1}^4&\int_{\cZ_i}\eps_{\inn}^{(i)}\nabla u^{(i)}\cdot\nabla\phi^{(i)}\dd x=\lim_{n\to\infty}\sum_{i=1}^4\int_{\cZ_{i,n}}\eps_{\inn}^{(i)}\nabla u^{(i)}\cdot\nabla\phi^{(i)}\dd x\nonumber\\
			=&\lim_{n\to\infty}\Big[\sum_{i=1}^4-\int_{\cZ_{i,n}}\eps_{\inn}^{(i)}(\Delta u^{(i)})\phi^{(i)}\dd x-\int_{\cF_{1,n}^{2}}\eps_{\inn}^{(1)}(\partial_2u^{(1)})\phi^{(1)}\dd\sigma\nonumber\\
			&-\int_{\cF_{1,n}^{1}}\eps_{\inn}^{(1)}(\partial_1u^{(1)})\phi^{(1)}\dd\sigma+\int_{\cF_{2,n}^{1}}\eps_{\inn}^{(2)}(\partial_1u^{(2)})\phi^{(2)}\dd\sigma-\int_{\cF_{2,n}^{2}}\eps_{\inn}^{(2)}(\partial_2u^{(2)})\phi^{(2)}\dd\sigma\nonumber\\
			&+\int_{\cF_{3,n}^{2}}\eps_{\inn}^{(3)}(\partial_2u^{(3)})\phi^{(3)}\dd\sigma+\int_{\cF_{3,n}^{1}}\eps_{\inn}^{(3)}(\partial_1u^{(3)})\phi^{(3)}\dd\sigma-\int_{\cF_{4,n}^{1}}\eps_{\inn}^{(4)}(\partial_1u^{(4)})\phi^{(4)}\dd\sigma\nonumber\\
			&+\int_{\cF_{4,n}^{2}}\eps_{\inn}^{(4)}(\partial_2u^{(4)})\phi^{(4)}\dd\sigma\Big].\label{eq 6.4}
		\end{align} 
		We treat the terms on the right hand side of (\ref{eq 6.4}) separately in the next two steps.
		
		3.i) We first handle the second and last summand.
		Note that we omit the arguments of the cut-off functions $\chi_{1,2}$ and $\chi_3$ from (\ref{eq 6.2}) in the following. 
		For $n$ sufficiently large, an integration by parts then leads to the identities 
		\begin{align*}
			&-\int_{\cF_{1,n}^{2}}\eps_{\inn}^{(1)}(\partial_2u^{(1)})\phi^{(1)}\dd\sigma+\int_{\cF_{4,n}^{2}}\eps_{\inn}^{(4)}(\partial_2u^{(4)})\phi^{(4)}\dd\sigma\\
			&=-\int_{\cF_{1,n}^{2}}\eps_{\inn}^{(1)}(\partial_2\psi^{(1)})\phi^{(1)}\dd\sigma+\int_{\cF_{4,n}^{2}}\eps_{\inn}^{(4)}(\partial_2\psi^{(4)})\phi^{(4)}\dd\sigma\nonumber\\
			&-\int_{\cF_{1,n}^{2}}\eps_{\inn}^{(1)}\chi_3(\tilde{\mE}_2^{(1)}\partial_2\chi_{1,2}+\chi_{1,2}\divv\tilde{\mE}^{(1)}-\chi_{1,2}(\partial_1\tilde{\mE}_1^{(1)}+\partial_3\tilde{\mE}_3^{(1)}))\phi^{(1)}\dd\sigma\nonumber\\
			&+\int_{\cF_{4,n}^{2}}\eps_{\inn}^{(4)}\chi_3(\tilde{\mE}_2^{(4)}\partial_2\chi_{1,2}+\chi_{1,2}\divv\tilde{\mE}^{(4)}-\chi_{1,2}(\partial_1\tilde{\mE}_1^{(4)}+\partial_3\tilde{\mE}_3^{(4)}))\phi^{(4)}\dd\sigma\nonumber\\
			&=-\int_{\cF_{1,n}^{2}}\eps_{\inn}^{(1)}(\partial_2\psi^{(1)})\phi^{(1)}\dd\sigma+\int_{\cF_{4,n}^{2}}\eps_{\inn}^{(4)}(\partial_2\psi^{(4)})\phi^{(4)}\dd\sigma\nonumber\\
			&-\int_{\cF_{1,n}^{2}}\eps_{\inn}^{(1)}\chi_3\big[(\tilde{\mE}_2^{(1)}\partial_2\chi_{1,2}+\chi_{1,2}\divv\tilde{\mE}^{(1)})\phi^{(1)}+\sum_{j\in\{1,3\}}\partial_j(\chi_{1,2}\phi^{(1)})\tilde{\mE}_j^{(1)}\big]\dd\sigma\nonumber\\
			&+\int_{\cF_{4,n}^{2}}\eps_{\inn}^{(4)}\chi_3\big[(\tilde{\mE}_2^{(4)}\partial_2\chi_{1,2}+\chi_{1,2}\divv\tilde{\mE}^{(4)})\phi^{(1)}+\sum_{j\in\{1,3\}}\partial_j(\chi_{1,2}\phi^{(1)})\tilde{\mE}_j^{(4)}\big]\dd\sigma.
		\end{align*}
		Using the definition of $X_2$ in (\ref{def X2}), $\divv\tilde{\mE}$ is piecewise $H^1$-regular on $\cZ$.
		By Lemma~\ref{lemma 6.3}, the same is true for $\tilde{\mE}_1$ and $\tilde{\mE}_2$.
		We then infer the formula
		\begin{align}
			\lim_{n\to\infty}\Big[&-\int_{\cF_{1,n}^{2}}\eps_{\inn}^{(1)}(\partial_2u^{(1)})\phi^{(1)}\dd\sigma+\int_{\cF_{4,n}^{2}}\eps_{\inn}^{(4)}(\partial_2u^{(4)})\phi^{(4)}\dd\sigma\Big]\label{eq 6.5}\\
			&=-\int_{\cF_4^{\cZ}}\big(\eps_{\inn}\phi\llbracket\partial_2\psi+\tilde{\mE}_2\chi_3\partial_2\chi_{1,2}+\chi_{1,2}\chi_3\divv\tilde{\mE}\rrbracket_{\cF_4^{\cZ}}\big)\dd\sigma.\nonumber
		\end{align}
		
		For the fifth and sixth summands on the right hand side of (\ref{eq 6.4}), we employ the above reasoning again. 
		Taking additionally the relation $\eps_{\inn}^{(2)}=\eps_{\inn}^{(3)}$, the construction of $\psi$ in (\ref{*const psi}) and the condition $\divv(\eps\mE)\in H^1(\tilde{Q}_{i,l})$ for $i\in\{1,\dots,L\}$, $l\in\{0,\dots,K\}$, into account, we derive the result
		\begin{align}
			\lim_{n\to\infty}\Big[-\int_{\cF_{2,n}^{2}}\eps_{\inn}^{(2)}(\partial_2u^{(2)})\phi^{(2)}\dd\sigma+\int_{\cF_{3,n}^{2}}\eps_{\inn}^{(3)}(\partial_2u^{(3)})\phi^{(3)}\dd\sigma\Big]=0.\label{eq 6.6}
		\end{align}
		
		3.ii) The focus now lies on the seventh and eighth terms on the right hand side of (\ref{eq 6.4}).
		In view of the definition of $W$ in Lemma~\ref{lemma 6.6}, an integration by parts leads to the relations
		\begin{align*}
			&-\int_{\cF_{4,n}^{1}}\eps_{\inn}^{(4)}(\partial_1u^{(4)})\phi^{(4)}\dd\sigma+\int_{\cF_{3,n}^{1}}\eps_{\inn}^{(3)}(\partial_1u^{(3)})\phi^{(3)}\dd\sigma\\
			&=-\int_{\cF_{4,n}^{1}}\eps_{\inn}^{(4)}\big(\chi_3\tilde{\mE}_2^{(4)}\partial_1\chi_{1,2}+\chi_3\chi_{1,2}(\curl\tilde{\mE}^{(4)})_3+\partial_1\psi^{(4)}-\chi_3\chi_{1,2}\partial_2\tilde{\mE}_1^{(4)}\big)\phi^{(4)}\dd\sigma\nonumber\\
			&\quad+\int_{\cF_{3,n}^{1}}\eps_{\inn}^{(3)}\big(\chi_3\tilde{\mE}_2^{(3)}\partial_1\chi_{1,2}+\chi_3\chi_{1,2}(\curl\tilde{\mE}^{(3)})_3+\partial_1\psi^{(3)}-\chi_{1,2}\chi_3\partial_2\tilde{\mE}_1^{(3)}\big)\phi^{(3)}\dd\sigma.
		\end{align*}
		Employing Lemma~\ref{lemma 6.5}, the mapping $\tilde{\mE}_2\partial_1\chi_{1,2}$ is piecewise $H^2$-regular on $\cZ$.
		By definition of $W$ and Lemma~7.1 in \cite{ZerullaDiss}, $\llbracket\partial_j(\chi_{1,2}\phi)\rrbracket_{\cF_3^{\cZ}}=0$ for $j\in\{2,3\}$, and $\tr_{\cF_3^{\cZ}}\partial_1\psi=0$, by definition of $\psi$ in (\ref{*const psi}).
		Using Proposition~4.6 in \cite{Zeru22b}, $(\curl\tilde{\mE})_3$ is furthermore piecewise $H^1$-regular on $\cZ$.
		Altogether, we conclude
		\begin{align}\label{eq 6.7}
			&\lim_{n\to\infty}\Big[-\int_{\cF_{4,n}^{1}}\eps_{\inn}^{(4)}(\partial_1u^{(4)})\phi^{(4)}\dd\sigma+\int_{\cF_{3,n}^{1}}\eps_{\inn}^{(3)}(\partial_1u^{(3)})\phi^{(3)}\dd\sigma\Big]\nonumber\\
			&=-\int_{\cF_3^{\cZ}}\!\!\big(\llbracket\eps_{\inn}\rrbracket_{\cF_3^{\cZ}}\chi_3(\tilde{\mE}_2\partial_1\chi_{1,2}+\chi_{1,2}(\curl\tilde{\mE})_3)-\chi_{1,2}\chi_3\partial_2\llbracket\eps_{\inn}\tilde{\mE}_1\rrbracket_{\cF_3^{\cZ}}\big)\phi\dd\sigma.
		\end{align}
		(Note here that the function $\partial_2\llbracket\eps_{\inn}\tilde{\mE}_1\rrbracket_{\cF_3^{\cZ}}$ belongs to $H^{1/2}(\cF_3^{\cZ})$ as a result of the definition of $\tilde{\mE}_1$ and the condition $\llbracket\eps\mE\cdot\nu_{\cF}\rrbracket_{\cF}\in H^{3/2}_{00}(\cF)$ for $\cF\in\Finteff$.)
		The same procedure and the fact $\eps_{\inn}^{(1)}=\eps_{\inn}^{(2)}$ also yield the limit
		\begin{align}
			\lim_{n\to\infty}\Big[-\int_{\cF_{1,n}^{1}}\eps_{\inn}^{(1)}(\partial_1u^{(1)})\phi^{(1)}\dd\sigma+\int_{\cF_{2,n}^{1}}\eps_{\inn}^{(2)}(\partial_1u^{(2)})\phi^{(2)}\dd\sigma\Big]=0
			\label{eq 6.8}
		\end{align}
		for the third and fourth summand on the right hand side of (\ref{eq 6.4}).
		
		4) Combining (\ref{eq 6.4})--(\ref{eq 6.8}), we arrive at the desired result
		\begin{align*}
			\sum_{i=1}^4\int_{\cZ_i}&\eps_{\inn}^{(i)}\nabla u^{(i)}\cdot\nabla\phi^{(i)}\dd x
			=-\sum_{i=1}^4\int_{\cZ_i}\eps_{\inn}(\Delta u^{(i)})\phi^{(i)}\dd x\\
			&\quad-\int_{\cF_3^{\cZ}}\Big[\llbracket\eps_{\inn}\rrbracket_{\cF_3^{\cZ}}\chi_3(\tilde{\mE}_2\partial_1\chi_{1,2}+\chi_{1,2}(\curl\tilde{\mE})_3)-\chi_{1,2}\chi_3\partial_2\llbracket\eps_{\inn}\tilde{\mE}_1\rrbracket_{\cF_3^{\cZ}}\Big]\phi\dd\sigma\\
			&\quad-\int_{\cF_4^{\cZ}}\llbracket\partial_2\psi+\tilde{\mE}_2\chi_3\partial_2\chi_{1,2}+\chi_{1,2}\chi_3\divv\tilde{\mE}\rrbracket_{\cF_4^{\cZ}}(\eps_{\inn}\phi)\dd\sigma.
		\end{align*}
	In view of Lemma~\ref{lemma 6.6}, this identity also holds for all $\phi\in\cD(-\Delta_{\cZ})^{1/2}$.
	The asserted energy estimate is due to Lemma~\ref{lemma 6.5} and due to the choice of $\chi_{1,2}$, $\chi_3$.
	The formula $g=\chi_{1,2}(\tfrac{2}{3}\lvert(x_1,x_2)\rvert)\chi_3(\tfrac{2}{3}\lvert x_3-\tfrac{1}{2}\rvert)g$ follows from (\ref{*const psi}).
	\end{proof}
\end{lemma}

We next analyze the transmission problem in Lemma~\ref{lemma 6.8}.
The proof follows the same strategy as the one for Lemma~\ref{lemma 4.1}, whence we only sketch the relevant arguments.
From Section~\ref{section 3.3}, we recall the one-dimensional space
\begin{align*}
	\check{N}_{\inn}=\spanv\{\chi(r)r^{\kappa_{1}}\psi_{1}(\varphi)\},
\end{align*} 
with $\kappa_{1}$ and $\psi_1$ denoting the first eigenvalue and an associated eigenfunction of (\ref{eq 3.5}). 
The next statement also uses the smooth cut-off functions $\chi_{1,2}$ and $\chi_3$ that arise in (\ref{eq 6.2}), as well as the number $\sskappa$ from (\ref{eq 3.6}).

\begin{lemma}\label{lemma 6.9}
	Let $\eps,\mu$ satisfy (\ref{assumptions parameters}), and let $\eps_{\inn}$ satisfy (\ref{eq 3.2}).
	Let additionally $g\in\bigcap_{j\in\{3,4\}}H^{1/2}(\cF_j^{\cZ})$, $f\in L^2(\cZ)$, and $\theta\in(\sskappa,1)$.
	There is a unique function $u\in\cD(-\Delta_{\cZ})^{1/2}$ with
	\begin{align}
		\sum_{i=1}^4\int_{\cZ_i}\eps_{\inn}\nabla u^{(i)}\cdot\nabla\phi^{(i)}\dd x=&\int_{\cZ}\eps_{\inn}f\phi\dd x-\int_{\cF_3^{\cZ}}\chi_{1,2}(\tfrac{2}{3}\lvert x_2\rvert)\chi_3(\tfrac{2}{3}\lvert x_3-\tfrac{1}{2}\rvert)g\phi\dd\sigma\nonumber\\
		&-\int_{\cF_4^{\cZ}}\chi_{1,2}(\tfrac{2}{3}\lvert x_1\rvert)\chi_3(\tfrac{2}{3}\lvert x_3-\tfrac{1}{2}\rvert)g(\eps_{\inn}\phi)\dd\sigma,\label{eq 6.9}
	\end{align}
	for $\phi\in\cD(-\Delta_{\cZ})^{1/2}$.
	The product $\chi_3(\lvert x_3-\tfrac{1}{2}\rvert)u$ belongs to $L^2((-1,2),PH^{1+\theta}(D)\oplus\check{N}_{\inn})$ with
	\begin{align*}
		\lVert\chi_3(\lvert x_3-\tfrac{1}{2}\rvert)u\rVert_{L^2((-1,2),PH^{1+\theta}(D)\oplus\check{N}_{\inn})}\leq C\big(\lVert f\rVert_{L^2(\cZ)}+\sum_{j=3}^4\lVert g\rVert_{H^{1/2}(\cF_j^{\cZ})}\big),
	\end{align*}
	where $C=C(\eps,\mu,\theta,Q)>0$ is a positive constant.
	\begin{proof}
		1) Throughout, $C=C(\eps,\mu,\theta,Q)>0$ is a positive constant that is allowed to change from line to line.
		We first consider (\ref{eq 6.9}) for a function $g\in L^2(\cF_{3}^{\cZ})\cap L^2(\cF_4^{\cZ})$.
		The Lax-Milgram lemma then provides a unique solution $u=u(g|_{\cF_3^{\cZ}},g|_{\cF_{4}^{\cZ}},f)\in\cD(-\Delta_{\cZ})^{1/2}$ of (\ref{eq 6.9}), giving rise to a well-defined linear solution operator 
		\begin{align}
			T:\prod_{j=3}^{4}L^2(\cF_j^{\cZ})\times L^2(\cZ)&\to L^2((-1,2),PH^1(D)),\label{eq 6.9.1}\\
			(g|_{\cF_{3}^{\cZ}},g|_{\cF_{4}^{\cZ}},f)&\mapsto \chi_3(\lvert x_3-\tfrac{1}{2}\rvert)u.\nonumber
		\end{align}
		Inserting $u=\phi$ into (\ref{eq 6.9}), the trace theorem for $H^1$-functions, the Poincaré inequality, and the Cauchy-Schwarz inequality give additionally rise to the estimate
		\begin{align}
			\lVert\chi_3(\lvert x_3-\tfrac{1}{2}\rvert) u\Vert_{PH^1(\cZ)}\leq C\Big(\lVert f\rVert_{L^2(\cZ)}+\sum_{j=3}^4\lVert g|_{\cF_j^{\cZ}}\rVert_{L^2(\cF_j^{\cZ})}\Big).\label{eq 6.10}
		\end{align}   
		In particular, $T$ is bounded.
		
		2) Fix $f\in L^2(\cZ)$, and suppose additionally that $g\in H^{1/2}_{00}(\cF_{3}^{\cZ})\cap H^{1/2}_{00}(\cF_{4}^{\cZ})$.
		First an appropriate extension of $g$ is constructed.
		We abbreviate $g_j:=g|_{\cF_j^{\cZ}}$, and consider the polyhedron
		\begin{align*}
			P:=\{(x,y)\ |\ x\in[0,1],\ y\in[x-1,0]\}\times[-1,2]\sub\overline{\cZ_4}.
		\end{align*}
		The interfaces $\cF_{3}^{\cZ}$ and $\cF_{4}^{\cZ}$ are in particular faces of $P$.
		Propositions~2.2--2.3 in \cite{AssousCiarlet97} provide an extension $h\in H^2(P)\cap H^1_0(P)$ with $\partial_{\nu}h=-\tfrac{1}{\eps_{\inn}^{(4)}}g_3$ on $\cF_{3}^{\cZ}$, $\partial_{\nu}h=0$ on all other faces of $P$, and $\lVert h\rVert_{H^2(P)}\leq C_1\lVert g\rVert_{H^{1/2}_{00}(\cF_3^{\cZ})}$, involving a uniform constant $C_1>0$.
		By means of Stein's extension operator, $h$ can be extended to a $H^2$-regular function on $\cZ_4$, still denoted by $h$.
		We then set $h:=0$ on $\cZ\setminus\cZ_4$.
		This gives rise to the desired formula $\llbracket\eps_{\inn}\nabla h\cdot\nu_{\cF_3^{\cZ}}\rrbracket_{\cF_3^{\cZ}}=g_3$ on $\cF_{3}^{\cZ}$.
		Similarly, we obtain a function $\tilde{h}\in PH^2(\cZ)$ with $\llbracket\nabla\tilde{h}\cdot\nu_{\cF_4^{\cZ}}\rrbracket_{\cF_4^{\cZ}}=g_4$, $\tilde{h}=0$ on $\cF_4^{\cZ}$, and $\tilde{h}=0=\partial_{\nu_{\cF_j^{\cZ}}}\tilde{h}$ for $j\in\{1,2,3\}$.
		We define
		\begin{align*}
			u_1(x_1,x_2,x_3):=\chi_{1,2}(\tfrac{2}{3}\lvert(x_1,x_2)\rvert)\chi_3(\tfrac{2}{3}\lvert x_3-\tfrac{1}{2}\rvert)\big(h(x_1,x_2,x_3)+\tilde{h}(x_1,x_2,x_3)\big),
		\end{align*}
		for $(x_1,x_2,x_3)\in\cZ$.
		This function is piecewise $H^2$-regular with
		\begin{align}
			\lVert u_1\rVert_{PH^2(\cZ)}&\leq C\sum_{j=3}^4\lVert g_j\rVert_{H^{1/2}_{00}(\cF_j^{\cZ})}.\label{eq 6.11}
		\end{align}
		
		3) The operator $\Delta_{\cZ}$ being invertible, there is a unique function $u_2\in\cD(\Delta_{\cZ})$ with $-\Delta_{\cZ}u_2=f+\Delta u_1$. 
		(Here, $\Delta u_1$ denotes the piecewise defined Laplacian of $u_1$ on $\cZ$.)
		Then $u:=u_1+u_2$ is the unique solution of (\ref{eq 6.9}).
		
		We analyze the product $\chi_3(\lvert x_3-\tfrac{1}{2}\rvert)u$ in the following, using ideas and techniques from \cite{Dauge,CoDa}.
		For convenience, we omit the argument of $\chi_3$ in the following.
		The functions $\chi_3u,\chi_3u_1,\chi_3u_2$ are trivially extended in $x_3$-direction to the infinite cylinder $D\times\RR$.
		Put now
		\begin{align}
			-\Delta(\chi_3 u_2^{(i)})=:\tilde{f}^{(i)}\in L^2(D_{\inn,i}\times \RR),\qquad i\in\{1,\dots,4\}.\label{eq 6.12}
		\end{align}
		The above extension procedure then leads to the estimate
		\begin{align}
			\lVert\tilde{f}\rVert_{L^2(D\times \RR)}\leq C\big(\lVert f\rVert_{L^2(\cZ)}+\sum_{j=3}^4\lVert g_j\rVert_{H^{1/2}_{00}(\cF_j^{\cZ})}\big).\label{eq 6.13}
		\end{align}
		
		In the following, we employ the partial Fourier transform with respect to the $x_3$-variable. 
		Its application to a function $w\in L^2(D\times \RR)$ is denoted by $\hat{w}$, the inverse transform is called $\check{w}$, and the new variable in Fourier space is $\xi$.
		Transforming (\ref{eq 6.12}) then leads to the formula
		\begin{align}
			(\xi^2-\partial_{1}^2-\partial_{2}^2)\,\hat{}\, (\chi_3 u_2)(x_1,x_2,\xi)=\hat{}\,(\tilde{f}),\qquad (x_1,x_2,\xi)\in D\times\RR.\label{eq 6.14}
		\end{align}
		The variable $\xi$ is considered to be fixed in the next steps (and the statements are then tacitly valid for almost all $\xi$).
		
		4) In the following, we employ several times that  $\hat{}\,(\chi_3u_2)(\cdot,\xi)$ belongs to $\cD(L_{\inn})$, see (\ref{eq 3.4}).
		Taking the relation
		\begin{align*}
			0\leq -\Re\big(\eps_{\inn}\xi^2\,\hat{}\,(\chi_3u_2)(\cdot,\xi),L_{\inn}\hat{}\,(\chi_3u_2)(\cdot,\xi)\big)_{L^2(D)}
		\end{align*}
		into account that follows from the positivity of $-L_{\inn}$, we obtain the inequality
		\begin{align}
			\lVert\sqrt{\eps_{\inn}}&L_{\inn}\,\hat{}\,(\chi_3u_2)(\cdot,\xi)\rVert_{L^2(D)}^2\leq\lVert\sqrt{\eps_{\inn}}\xi^2\,\hat{}\,(\chi_3u_2)(\cdot,\xi)\rVert_{L^2(D)}^2\nonumber\\
			&-2\Re\big(\eps_{\inn}\xi^2\,\hat{}\,(\chi_3u_2)(\cdot,\xi),L_{\inn}\hat{}\,(\chi_3u_2)(\cdot,\xi)\big)_{L^2(D)}+\lVert\sqrt{\eps_{\inn}}L_{\inn}\hat{}\,(\chi_3u_2)(\cdot,\xi)\rVert_{L^2(D)}^2\nonumber\\
			&=\lVert\sqrt{\eps_{\inn}}\,\hat{}\,(\tilde{f})(\cdot,\xi)\rVert_{L^2(D)}^2.\label{eq 6.14.1}
		\end{align}
		We next integrate this estimate with respect to $\xi$, and use Plancherel's Theorem.
		Taking also Proposition~\ref{proposition 3.5} into account, we derive the estimate
		\begin{align*}
			\lVert\chi_3 u_2\rVert_{L^2(\RR,PH^2(D)\oplus\check{N}_{\inn})}\leq C\lVert\tilde{f}\rVert_{L^2(D\times\RR)}.
		\end{align*}
		The definition $u=u_1+u_2$, (\ref{eq 6.11}) and (\ref{eq 6.13}) then result in the desired inequality
		\begin{align*}
			\lVert \chi_3u\rVert_{L^2((-1,2),PH^{2}(D)\oplus\check{N}_{\inn})}\leq C\big(\lVert f\rVert_{L^2(\cZ)}+\sum_{j=3}^4\lVert g_j\rVert_{H^{1/2}_{00}(\cF_j^{\cZ})}\big).
		\end{align*}
		This means that the operator $T$ from (\ref{eq 6.9.1}) is bounded from $\prod_{j\in\{3,4\}}H^{1/2}_{00}(\cF_j^{\cZ})\times L^2(\cZ)$ into $L^2((-1,2),PH^2(D)\times\check{N}_{\inn})$.
		
		5) We now interpolate between the results of parts~1) and 4).
		The relation $H^{s}_0(\cF_j^{\cZ})=H^{s}(\cF_j^{\cZ})$, $s\in(0,1/2)$, first yields the identity
		\begin{align*}
			\Big(\prod_{j=3}^4&L^2(\cF_j^{\cZ})\times L^2(\cZ),\prod_{j=3}^4H^{1/2}_{00}(\cF_j^{\cZ})\times L^2(\cZ)\Big)_{\theta,2}=\prod_{j=3}^4H^{\theta/2}(\cF_j^{\cZ})\times L^2(\cZ).
		\end{align*}
		Remark~14.4 in \cite{LionsMagenes} and Theorem~1 in \cite{Zeru22a} further lead to the formula
		\begin{align*}
			\big(L^2((-1,2),PH^1(D)),L^2((-1,2)&,PH^2(D)\oplus\check{N}_{\inn})\big)_{\theta,2}\nonumber\\
			&=L^2((-1,2),(PH^1(D),PH^2(D)\oplus\check{N}_{\inn})_{\theta,2})\\
			&=L^2((-1,2),PH^{1+\theta}(D)\oplus\check{N}_{\inn}).
		\end{align*}
		We then conclude that $T$ is a bounded operator from $\prod_{j=3}^{4}H^{\theta/2}(\cF_j^{\cZ})\times L^2(\cZ)$ into the space $L^2((-1,2),PH^{1+\theta}(D)\oplus\check{N}_{\inn})$.
	\end{proof}
\end{lemma}

Combining our findings in Lemmas~\ref{lemma 6.7}--\ref{lemma 6.9}, we can now provide the desired regularity statement for the second component of the electric field near interior edges.
For the statement, recall the number $\sskappa$ from (\ref{eq 3.6}).

\begin{lemma}\label{lemma 6.10}
	Let $\eps,\mu$ satisfy (\ref{assumptions parameters}), $\eps_{\inn}$ satisfy (\ref{eq 3.2}), $\theta\in (\sskappa,1)$, and $(\mE,\mH)\in X_2$.
	Let $v$ be defined by (\ref{eq 6.2}). 
	Then the function $\chi_3(\lvert x_3-\tfrac{1}{2}\rvert)v$ belongs to $PH^{1+\theta}(\cZ)$ with
	\begin{align*}
		\lVert\chi_3(\lvert x_3-\tfrac{1}{2}\rvert) v\rVert_{PH^{1+\theta}(\cZ)}\leq &C\Big(\sum_{i=1}^N\big(\lVert\mE\rVert_{L^2(Q_i)}+\lVert\Delta\mE\rVert_{L^2(Q_i)}+\lVert\curl\mE\rVert_{H^1(Q_i)}\big)\\
		&+\sum_{i=1}^L\sum_{l=0}^K\lVert\divv\mE\rVert_{H^1_{00}(\tilde{Q}_{i,l})}+\sum_{\cF\in\Finteff}\lVert\llbracket\eps\mE\cdot\nu_{\cF}\rrbracket_{\cF}\rVert_{H^{3/2}_{00}(\cF)}\Big),
	\end{align*}
	involving a uniform constant $C=C(\eps,\mu,\theta,Q)>0$.	
	\begin{proof}
		1) Throughout the proof, $C=C(\eps,\mu,\theta,Q)>0$ is a constant that is allowed to change from line to line.
		Let $\psi\in PH^2(\cZ)$ be the function from Lemma~\ref{lemma 6.7}.
		By Lemma~\ref{lemma 6.8}, $v+\psi$ satisfies (\ref{eq 6.9}) with appropriate right hand sides $f\in L^2(\cZ)$ and $g\in \prod_{j\in\{3,4\}}H^{1/2}(\cF_j^{\cZ})$. 
		Lemma~\ref{lemma 6.9} then shows that $\chi_3(\lvert x_3-\tfrac{1}{2}\rvert)(v+\psi)$ belongs to $L^2((-1,2),PH^{1+\theta}(D)\oplus\check{N}_{\inn})$.
		Since $v$ is piecewise $H^{3/2}$-regular, see Lemma~\ref{lemma 6.3}, we infer $\chi_3(\lvert x_3-\tfrac{1}{2}\rvert)(v+\psi)\in L^2((-1,2),PH^{3/2}(D))$.
		
		Lemma~\ref{lemma 3.3} implies that $\kappa_{1}<1/3$, whence $\check{N}_{\inn}$ is no subset of $PH^{3/2}(D)$, see \cite{BabuskaAnderssonGuoMelenkOh1996} for instance.
		In view of the piecewise $H^2$-regularity of $\psi$, this altogether means that $\chi_3(\lvert x_3-\tfrac{1}{2}\rvert)v$ belongs to $L^2((-1,2),PH^{1+\theta}(D))$.
		For convenience, we omit the argument of $\chi_3$ in the following.
		Combining the triangle inequality with the estimates from Lemmas~\ref{lemma 6.5}, \ref{lemma 6.7}--\ref{lemma 6.9} and Proposition~4.6 in \cite{Zeru22b}, we infer the relations
		\begin{align}
			\lVert\chi_3 &v\rVert_{L^2((-1,2),PH^{1+\theta}(D))}\leq \lVert\chi_3(v+\psi)\rVert_{L^2((-1,2),PH^{1+\theta}(D))}+C\lVert\psi\rVert_{PH^2(\cZ)}\nonumber\\
			&\leq C\Big(\sum_{i=1}^N\big(\lVert\mE\rVert_{L^2(Q_i)}+\lVert\Delta\mE\rVert_{L^2(Q_i)}+\lVert\curl\mE\rVert_{H^1(Q_i)}\big)\nonumber\\
			&\qquad\quad +\sum_{i=1}^L\sum_{l=0}^K\lVert\divv\mE\rVert_{H^1_{00}(\tilde{Q}_{i,l})}+\sum_{\cF\in\Finteff}\lVert\llbracket\eps\mE\cdot\nu_{\cF}\rrbracket_{\cF}\rVert_{H_{00}^{3/2}(\cF)}\Big).\label{eq 6.15}
		\end{align} 
	
		2) To derive other useful estimates for $\chi_3v$, we employ the constructions of the proof for Lemma~\ref{lemma 6.9}, see also \cite{Dauge,CoDa}.
		Let $f\in L^2(\cZ)$, $g\in\prod_{j=3}^4L^2(\cF_j^{\cZ})$, and let $u=u(g,f)\in PH^1(\cZ)$ be the unique solution of (\ref{eq 6.9}).
		Extend the function $\chi_3(\lvert x_3-\tfrac{1}{2}\rvert)u$ trivially by zero in $x_3$-direction, and denote the extended function again by $\chi_3u$.
		As in (\ref{eq 6.10}), the extension-solution mapping 
		\begin{align*}
			T:\prod_{j=3}^4L^2(\cF_j^{\cZ})\times L^2(\cZ)&\to H^1(\RR,L^2(D))\cap L^2(\RR,PH^1(D)),\\
			(g,f)&\mapsto \chi_3u(g,f)
		\end{align*}
		is bounded. 
		
		3) Assume now additionally that $g\in\prod_{j=3}^4H^{1/2}_{00}(\cF_j^{\cZ})$.
		Recall also $u_1$ and $u_2$ from parts~2) and 3) of the proof for Lemma~\ref{lemma 6.9}, as well as $u=u_1+u_2$.
		The functions $\chi_3u_1$, and $\chi_3u_2$ are as above considered on the infinite cylinder $D\times \RR$ after trivial extension by zero.
		Note moreover that $\hat{w}$ denotes the partial Fourier transform in $x_3$ of a function $w\in L^2(\RR)$, and the variable in Fourier space is $\xi$.
		
		Fix $\xi\in\RR$ in the following.
		The reasoning for (\ref{eq 6.14.1}) and (\ref{eq 6.13}) also gives rise to the estimate
		\begin{align*}
			\lVert\xi^2\,\hat{}\,(\chi_3u_2)\rVert_{L^2(D\times\RR)}\leq C\Big(\lVert f\rVert_{L^2(\cZ)}+\sum_{j=3}^4\lVert g\rVert_{H^{1/2}_{00}(\cF_j^{\cZ})}\Big).
		\end{align*}
		With the triangle inequality and (\ref{eq 6.11}), we then deduce the inequality
		\begin{align}
			\lVert\chi_3u\rVert_{H^2(\RR,L^2(D))}\leq C\big(\lVert f\rVert_{L^2(\cZ)}+\sum_{j=3}^4\lVert g_j\rVert_{H^{1/2}_{00}(\cF_j^{\cZ})}\big).\label{eq 6.16}
		\end{align}
		We next employ the equivalence of the piecewise $H^1$-norm and the graph norm in $\cD(-L_{\inn})^{1/2}$ several times.
		With the symmetry of the operator $(-L_{\inn})^{1/2}$, and the Cauchy-Schwarz estimate, we then obtain the relations
		\begin{align*}
			\lVert&\sqrt{\eps_{\inn}}\lvert\xi\rvert(-L_{\inn})^{1/2}\,\hat{}\,(\chi_3u)(\cdot,\xi)\rVert_{L^2(D\times\RR)}\leq \lVert\sqrt{\eps_{\inn}}\lvert\xi\rvert(-L_{\inn})^{1/2}\,\hat{}\,(\chi_3u_2)(\cdot,\xi)\rVert_{L^2(D\times\RR)}\\
			&\quad+C\lVert\chi_3u_1\rVert_{H^1(\RR,PH^1(D))}\\
			&=\lVert\big(\eps_{\inn}\lvert\xi\rvert^2\,\hat{}\,(\chi_3u_2)(\cdot,\xi),L_{\inn}\,\hat{}\,(\chi_3u_2)(\cdot,\xi)\big)_{L^2(D)}^{1/2}\rVert_{L^2(\RR)}+C\lVert\chi_3u_1\rVert_{H^1(\RR,PH^1(D))}\\
			&\leq\Big\lVert\lVert\sqrt{\eps_{\inn}}\xi^2\,\hat{}\,(\chi_3u_2)(\cdot,\xi)\rVert_{L^2(D)}^{1/2}\lVert\sqrt{\eps_{\inn}}L_{\inn}\,\hat{}\,(\chi_3u_2)(\cdot,\xi)\rVert_{L^2(D)}^{1/2}\Big\rVert_{L^2(\RR)}\\
			&\quad+C\lVert\chi_3u_1\rVert_{H^1(\RR,PH^1(D))}\\
			&\leq\lVert\sqrt{\eps_{\inn}}\xi^2\,\hat{}\,(\chi_3u_2)(\cdot,\xi)\rVert_{L^2(D\times\RR)}^{1/2}\lVert\sqrt{\eps_{\inn}}L_{\inn}\,\hat{}\,(\chi_3u_2)(\cdot,\xi)\rVert_{L^2(D\times\RR)}^{1/2}\\
			&\quad+C\lVert\chi_3u_1\rVert_{H^1(\RR,PH^1(D))}.
		\end{align*}
	Combining (\ref{eq 6.16}), (\ref{eq 6.14.1}), (\ref{eq 6.13}) and (\ref{eq 6.11}), we conclude the bound
	\begin{align*}
		\lVert\chi_3u\rVert_{H^1(\RR,PH^1(D))}\leq C\big(\lVert f\rVert_{L^2(\cZ)}+\sum_{j=3}^4\lVert g\rVert_{H^{1/2}_{00}(\cF_j^{\cZ})}\big).
	\end{align*}	
	In view of (\ref{eq 6.16}), this means that $T$ is also bounded as
	\begin{align*}
		T:\prod_{j=3}^4H^{1/2}_{00}(\cF_j^{\cZ})\times L^2(\cF_j^{\cZ})&\to H^2(\RR,L^2(D))\cap H^1(\RR,PH^1(D)).
	\end{align*}

	4) Interpolating between the boundedness results of parts~2) and 3), we infer that the trivial extension $\chi_3(v+\psi)$ is an element of the space
	\begin{align*}
		\Big(H^1(\RR,L^2(D))&\cap L^2(\RR,PH^1(D)),H^2(\RR,L^2(D))\cap H^1(\RR,PH^1(D))\Big)_{\theta,2}\\
		&\sub H^{1+\theta}(\RR,L^2(D))\cap H^1(\RR,PH^{\theta}(D))\cap H^{\theta}(\RR,PH^1(D)).
	\end{align*}
	Taking also Lemma~\ref{lemma 6.7} and (\ref{eq 6.15}) into account, we conclude that $\chi_3v$ is an element of $PH^{1+\theta}(\cZ)$ satisfying the asserted estimate.
	\end{proof}
\end{lemma}

Combining Lemmas~\ref{lemma 6.5} and \ref{lemma 6.10}, we deduce the desired global regularity statement for the electric field component of a vector in $X_2$.

\begin{lemma}\label{lemma 6.11}
	Let $\eps$ and $\mu$ satisfy (\ref{assumptions parameters}), $\theta\in(0,1)$, and let $(\mE,\mH)\in X_2$. 
	The vector $\mE^{(i)}$ belongs to $H^{1+\theta}(Q_i)^3$ with
	\begin{align*}
		\lVert\mE\rVert_{PH^{1+\theta}(Q)}\leq C\Big(&\sum_{i=1}^N\big(\lVert\mE^{(i)}\rVert_{L^2(Q_i)}+\lVert\Delta\mE^{(i)}\rVert_{L^2(Q_i)}+\lVert\curl\mE\rVert_{H^1(Q_i)}\big)\\
		&+\sum_{i=1}^L\sum_{l=0}^K\lVert\divv\mE\rVert_{H^1_{00}(\tilde{Q}_{i,l})}+\sum_{\cF\in\Finteff}\lVert\llbracket\eps\mE\cdot\nu_{\cF}\rrbracket_{\cF}\rVert_{H^{3/2}_{00}(\cF)}\Big),
	\end{align*}	
	involving a uniform number $C=C(\theta,\eps,\mu,Q)>0$.
	\begin{proof}
		An equivalent version of Lemma~\ref{lemma 6.10} is also valid for $\mE_1$
		due to symmetry.
		By Lemmas~\ref{lemma 6.5} and \ref{lemma 6.10}, it suffices to study only the third component $\mE_3$.
		To that end, we use the formulas
		\begin{align*}
			\partial_1\mE_3^{(i)}&=-(\curl\mE^{(i)})_2+\partial_3\mE_1^{(i)},\qquad \partial_2\mE_3^{(i)}=(\curl\mE^{(i)})_1+\partial_3\mE_2^{(i)},\\
			\partial_3\mE_3^{(i)}&=\divv\mE^{(i)}-\partial_1\mE_1^{(i)}-\partial_2\mE_2^{(i)},
		\end{align*}
		on $Q_i$ for $i\in\{1,\dots,N\}$.
		By definition of $X_2$, Proposition~4.6 of \cite{Zeru22b}, Remark~\ref{remark previous work}, Lemmas~\ref{lemma 6.5} and \ref{lemma 6.10} (respectively an equivalent version for $\mE_1$), the expressions on the right hand sides are at least $H^{\theta}$-regular on $Q_i$.
		This implies the asserted regularity statement.
		The estimate follows in the same way.
	\end{proof}
\end{lemma}

\subsection{Study of the magnetic field component}\label{section 5.3}\

In this section, we analyze the regularity of the magnetic field component $H$ of a vector $(\mE,\mH)\in X_2$, see (\ref{def X2}).
We first state a result that is obtained from Lemma~9.15 in \cite{ZerullaDiss}, Proposition~4.6 in \cite{Zeru22b} and Remark~\ref{remark previous work} via localization with respect to the interfaces in $\tFint$, see Section~\ref{section1}.
As we elaborate several cut-off arguments in Section~\ref{section 5.2}, we skip the proof here for the sake of brevity.

\begin{lemma}\label{lemma 6.12}
	Let $\eps,\mu$ satisfy (\ref{assumptions parameters}), let $(\mE,\mH)\in X_2$, and let $\cF\in\tFint$.
	There is an open set $\cO$ with $\cF\sub\cO$ and $\mH_1\in PH^2(\cO\cap Q)$.
	Furthermore, the estimate
	\begin{align*}
		\lVert\mH_1\rVert_{PH^2(\cO\cap Q)}\leq C\sum_{i=1}^N\big(\lVert\mH\rVert_{L^2(Q_i)}+\lVert\Delta\mH\rVert_{L^2(Q_i)}+\lVert\curl\mH\rVert_{L^2(Q_i)}\big)
	\end{align*}
	holds with a uniform constant $C=C(\eps,\mu,Q)>0$.
\end{lemma}

We next establish a regularity result for the first magnetic field component on the entire cuboid $Q$.
To reach our goal, we use ideas and techniques from the proof of Lemma~3.7 in \cite{HoJaSc}, and Lemma~9.15 in \cite{ZerullaDiss}.
See also the proof for Proposition~3.2 in \cite{EiSc32}.

\begin{lemma}\label{lemma 6.13}
	Let $\eps,\mu$ satisfy (\ref{assumptions parameters}), let $(\mE,\mH)\in X_2$, and $\theta\in(0,1)$. Then the function $\mH_1$ belongs to $PH^{1+\theta}(Q)$ with
	\begin{align*}
		\lVert \mH_1\rVert_{PH^{1+\theta}(Q)}\leq C\big(\sum_{i=1}^N\lVert\Delta\mH_1^{(i)}\rVert_{L^2(Q_i)}+\lVert(\curl\mH)_3\rVert_{H^1(Q_i)}\big),
	\end{align*}
	where $C=C(\eps,\mu,\theta,Q)>0$ is a uniform constant.
	\begin{proof}
		1) Define the spaces 
		\begin{align*}
			\tilde{V}&:=\{\varphi\in PH^1_{\Gamma_1}(Q)\ |\ \llbracket\mu\varphi\rrbracket_{\cF}=0,\ \text{for }\cF\in\Fintj{1},\ \llbracket\varphi\rrbracket_{\cF'}=0\text{ for }\cF'\in\Fintj{2}\},\\
			\tilde{W}&:=\{\varphi\in PH^2(Q)\cap\tilde{V}\ |\  \supp(\varphi)\cap\overline{\Gamma_1}=\emptyset,\ \varphi\text{ vanishes in a neighborhood of all}\\
			&\hphantom{:=\{\varphi\in PH^2(Q)\cap\tilde{V}\ |\ } \text{edges in }Q\}.
		\end{align*}
		Similar to Lemma~2.1 in \cite{Zeru22b}, the space $\tilde{W}$ is dense in $\tilde{V}$ with respect to the piecewise $H^1$-norm.
		Using the representation $Q_i:=(a_1^{-,i},a_1^{+,i})\times(a_2^{-,i},a_2^{+,i})\times(a_3^{-},a_3^{+})$, we introduce the smaller cuboids
		\begin{align}
			Q_{i,n}:=(a_1^{-,i}+\tfrac{1}{n},a_1^{+,i}-\tfrac{1}{n})\times(a_2^{-,i}+\tfrac{1}{n},a_2^{+,i}-\tfrac{1}{n})\times(a_3^{-}+\tfrac{1}{n},a_3^{+}-\tfrac{1}{n})\label{eq 6.16.1}
		\end{align}
		for $n\geq n_0$ sufficiently large.
		The corresponding faces of these cuboids are denoted analogously to the faces of $Q_1,\dots,Q_N$, meaning that $Q_{i,n}$ possesses the face pairs $\Gamma_{j,n}^{(i)}$ for $j\in\{1,2,3\}$.
		The exterior unit normal vector of $Q_{i,n}$ is denoted by $\nu_{i,n}$ with components $(\nu_{i,n})_l$ for $l\in\{1,2,3\}$.
		
		2) Proposition~4.6 in \cite{Zeru22b} and Remark~\ref{remark 6.2} show that $\mH|_{\tilde{Q}_j}\in H^1(\tilde{Q}_j)^3$, and that $\mH|_{Q_i}\in H_{\text{loc}}^2(Q_i)^3$ for $j\in\{1,\dots,L\}$, $i\in\{1,\dots,N\}$. 
		(See Section~\ref{section1} for the definition of $\tilde{Q}_j$.)
		Proposition~4.6 in \cite{Zeru22b} and Remark~\ref{remark previous work} furthermore yield $(\curl\mH)_3\in PH^1(Q)$, and $(\curl\mH)_j\in PH^{2/3}(Q)$ for $j\in\{1,2\}$.
		
		We next derive an elliptic transmission problem for $\mH_1$ fitting into the framework of Lemma~\ref{lemma 4.1}. 
		Let $\varphi\in\tilde{W}$.
		Integrating by parts, we infer the identities
		\begin{align}\label{eq 6.17}
			&\sum_{i=1}^N\int_{Q_i}\mu^{(i)}(\nabla\mH_1^{(i)})\cdot(\nabla\varphi^{(i)})\dd x=\lim_{n\to\infty}\sum_{i=1}^N\int_{Q_{i,n}}\mu^{(i)}(\nabla\mH_1^{(i)})\cdot(\nabla\varphi^{(i)})\dd x\nonumber\\
			&=\lim_{n\to\infty}\sum_{i=1}^N\Big(\int_{Q_{i,n}}-\mu^{(i)}(\Delta\mH_1^{(i)})\varphi^{(i)}\dd x+\int_{\partial Q_{i,n}}\!\!\mu^{(i)}(\nabla\mH_1^{(i)}\cdot\nu_{i,n})\varphi^{(i)}\dd\sigma\Big).
		\end{align}
		We next analyze the limit of the boundary integral terms on the right hand side.
		By inserting $\curl\mH$ and $\divv(\mu\mH)=0$, we obtain the formula
		\begin{align}
			&\lim_{n\to\infty}\sum_{i=1}^N\int_{\partial Q_{i,n}}\mu^{(i)}(\nabla\mH_1^{(i)}\cdot\nu_{i,n})\varphi^{(i)}\dd\sigma\nonumber\\
			&=\lim_{n\to\infty}\sum_{i=1}^N\Big(-\int_{\Gamma_{1,n}^{(i)}}\mu^{(i)}(\nu_{i,n})_{1}(\partial_2\mH_2^{(i)}+\partial_3\mH_3^{(i)})\varphi^{(i)}\dd\sigma\nonumber\\
			&\quad-\int_{\Gamma_{2,n}^{(i)}}\mu^{(i)}(\nu_{i,n})_2(\curl\mH^{(i)})_3\varphi^{(i)}\dd\sigma+\int_{\Gamma_{3,n}^{(i)}}\mu^{(i)}(\nu_{i,n})_3(\curl\mH^{(i)})_2\varphi^{(i)}\dd\sigma\nonumber\\
			&\quad+\int_{\Gamma_{2,n}^{(i)}}\mu^{(i)}(\nu_{i,n})_2(\partial_1\mH_2^{(i)})\varphi^{(i)}\dd\sigma+\int_{\Gamma_{3,n}^{(i)}}\mu^{(i)}(\nu_{i,n})_3(\partial_1\mH_3^{(i)})\varphi^{(i)}\dd\sigma\Big).\label{eq 6.18}
		\end{align}
		In view of the boundary condition $\tfrac{1}{\eps}\curl\mH\times\nu=0$ on $\partial Q$, we infer that the third summand on the right hand side tends to zero as $n\to\infty$.
		The boundary condition also implies the relation
		\begin{align}
			\lim_{n\to\infty}\sum_{i=1}^N\int_{\Gamma_{2,n}^{(i)}}\mu(\nu_{i,n})_2(\curl\mH)_3\varphi\dd \sigma\!=-\!\!\sum_{\cF\in\Fintj{2}}\int_{\cF}\llbracket\mu(\curl\mH)_3\rrbracket_{\cF}\varphi\dd\sigma.\label{eq 6.19}
		\end{align}
		
		We further note that the remaining terms on the right hand side of (\ref{eq 6.18}) all converge to zero.
		To show this claim, it suffices to focus on the fourth one. (All other terms can be treated with similar arguments.) 
		Recall that $\mu\mH\cdot\nu=0$ on $\partial Q$.
		Using the definition of $\tilde{W}$ and Lemma~7.1 in \cite{ZerullaDiss}, an integration by parts yields the identities
		\begin{align*}
			\lim_{n\to\infty}\int_{\Gamma_{2,n}^{(i)}}\mu^{(i)}(\nu_{i,n})_2(\partial_1\mH_2^{(i)})\varphi^{(i)}\dd\sigma=-\lim_{n\to\infty}\int_{\Gamma_{2,n}^{(i)}}\mu^{(i)}(\nu_{i,n})_2\mH_2^{(i)}\partial_1\varphi^{(i)}\dd\sigma=0.
		\end{align*}
		
		Combining (\ref{eq 6.17})--(\ref{eq 6.19}), we conclude the formula
		\begin{align*}
			\sum_{i=1}^N\int_{Q_i}\mu^{(i)}(\nabla\mH_1^{(i)})\cdot(\nabla\varphi^{(i)})\dd x=&-\sum_{i=1}^N\int_{Q_i}\mu^{(i)}(\Delta\mH_1^{(i)})\varphi^{(i)}\dd x\\
			&+\sum_{\cF\in\Fintj{2}}\int_{\cF}\llbracket\mu(\curl\mH)_3\rrbracket_{\cF}\varphi\dd\sigma.
		\end{align*}
		By density, the same formula holds also for all functions $\varphi\in \tilde{V}$. 
		Lemma~\ref{lemma 4.1} now yields the asserted statement. 
	\end{proof}
\end{lemma}

We finally treat the remaining two magnetic field components of vectors in $X_2$.
For the statement, we employ the number $\skappa\in(2/3,1)$ that is defined via the relation
\begin{align}
	\max_{\substack{i\in\{1,\dots,L\},\\ l\in\{1,\dots,K\}}}\frac{(\eps|_{\tilde{Q}_{i,l}}-\eps|_{\tilde{Q}_{i,0}})^2}{\eps|_{\tilde{Q}_{i,l}}\eps|_{\tilde{Q}_{i,0}}}=-\frac{4\sin^2(\skappa\pi)}{\sin(\tfrac{\skappa}{2}\pi)\sin(\tfrac{3\skappa}{2}\pi)},\label{eq 6.20}
\end{align}
compare formula~(3.4) in \cite{Zeru22b}.

\begin{lemma}\label{lemma 6.14}
	Let $\eps,\mu$ satisfy (\ref{assumptions parameters}), $(\mE,\mH)\in X_2$, $\theta\in(1/4,1/2)$, and  $\kappa\in(1-\skappa,1/3)$. 
	Then $\mH_2\in PH^{3/2+\theta}(Q)$ and $\mH_3\in PH^{2-\kappa}(Q)$ with
	\begin{align*}
		\lVert\mH_2\rVert_{PH^{3/2+\theta}(Q)}&\leq C\sum_{i=1}^N\big(\lVert\Delta\mH_2\rVert_{L^2(Q_i)}+\lVert(\curl\mH)_3\rVert_{H^1(Q_i)}+\lVert\mH_1\rVert_{H^{3/2+\theta}(Q_i)}\big),\\
		\lVert\mH_3\rVert_{PH^{2-\kappa}(Q)}&\leq C\sum_{i=1}^N\big(\lVert\Delta\mH_3\rVert_{L^2(Q_i)}+\lVert\mH_1\rVert_{H^{2-\kappa}(Q_i)}+\sum_{j=1}^2\lVert(\curl\mH)_j\rVert_{H^{1-\kappa}(Q_i)}\big),
	\end{align*}
	with a uniform constant $C=C(\theta,\kappa,\eps,\mu,Q)>0$.
	\begin{proof}
		1) The reasoning is for both components similar to the proof of Lemma~\ref{lemma 6.13} (see also the proofs for Lemma~3.7 in \cite{HoJaSc}, Proposition~3.2 in \cite{EiSc32}, and Lemma~9.15 in \cite{ZerullaDiss}).
		Hence we only analyze the third component $\mH_3$.
		We use the spaces of test functions
		\begin{align*}
			\tilde{V}&:=H^1_{\Gamma_3}(Q),\\
			\tilde{W}&:=\{\varphi\in PH^2(Q)\cap H^1_{\Gamma_3}(Q)\ |\ \varphi^{(i)}\text{ is smooth on }\overline{Q_i},\ \supp(\varphi)\cap\overline{\Gamma_3}=\emptyset,\\
			&\hphantom{:=\{\varphi\in PH^2(Q)\cap H^1_{\Gamma_3}(Q)\ |\ }\varphi\text{ vanishes in a neighborhood of all edges in }Q\}.
		\end{align*}
		Note that a similar reasoning as in Lemma~2.1 of \cite{ZerullaDiss} shows that $\tilde{W}$ is dense in $\tilde{V}$ with respect to the $H^1$-norm.
		We moreover reuse the small cuboids $Q_{i,n}$ with exterior unit normal $\nu_{i,n}$ and faces $\Gamma_{j,n}^{\pm,(i)}$ from (\ref{eq 6.16.1}).
		The effective interfaces in $\Finteff$ with unit normal vector $e_j$ are collected in a set $\cF^{\textrm{eff}}_{\textrm{int},j}$, $j\in\{1,2\}$.
		
		2) Let $\varphi\in\tilde{W}$.
		By Remark~\ref{remark 6.2}, the function $\mH_3^{(i)}$ belongs to $H^2_{\text{loc}}(Q_i)$. 
		The definition of $X_2$ in (\ref{def X2}), Proposition~4.6 in \cite{Zeru22b}, and Remark~\ref{remark previous work} furthermore imply that the first two components of $\curl\mH$ are elements of $PH^{1-\kappa}(Q)\sub PH^{2/3}(Q)$.
		
		In the following, we derive an equation for $\mH_3$ that fits into the framework of Lemma~\ref{lemma 4.1}.
		As in the proof of Lemma~\ref{lemma 6.13}, we first infer the formula
		\begin{align}\label{eq 6.21}
			\sum_{i=1}^N\int_{Q_i}(\nabla\mH_3^{(i)})\cdot(\nabla\varphi^{(i)})\dd x=\lim_{n\to\infty}&\sum_{i=1}^N\Big(-\int_{Q_{i,n}}(\Delta\mH_3^{(i)})\varphi^{(i)}\dd x\nonumber\\
			&\qquad+\int_{\partial Q_{i,n}}(\nabla\mH_3^{(i)}\cdot\nu_{i,n})\varphi^{(i)}\dd\sigma\Big).	
		\end{align}
		
		In the next steps, we analyze the limit of the second summand on the right hand side of (\ref{eq 6.21}).
		
		3.i) In view of the location of the support of $\varphi$, we first infer
		\begin{align}\label{eq 6.22}
			\lim_{n\to\infty}\sum_{i=1}^N\int_{\Gamma_{3,n}^{\pm,(i)}}(\nu_{i,n})_3(\partial_3\mH_3^{(i)})\varphi^{(i)}\dd \sigma=0.
		\end{align}
		
		3.ii) Inserting the vector $\curl\mH$, we moreover derive the formula	
		\begin{align}\label{eq 6.23}
			\sum_{i=1}^N\int_{\Gamma_{1,n}^{\pm,(i)}}(\nu_{i,n})_1(\partial_1\mH_3^{(i)})\varphi^{(i)}\dd\sigma=\sum_{i=1}^N\Big(&\int_{\Gamma_{1,n}^{\pm,(i)}}(\nu_{i,n})_{1}(\partial_3\mH_1^{(i)})\varphi^{(i)}\dd\sigma\\
			&-\int_{\Gamma_{1,n}^{\pm,(i)}}(\nu_{i,n})_1(\curl\mH^{(i)})_2\varphi^{(i)}\dd\sigma\Big).\nonumber
		\end{align}
		Taking the boundary conditions $\curl\mH\times\nu=0$, $\mu\mH\cdot\nu=0$ on $\partial Q$, and the definition of $\tilde{W}$ into account, an integration by parts on the right hand side of (\ref{eq 6.23}) shows that the integrals over faces that approach the boundary face $\Gamma_1$ tend to zero as $n\to\infty$.
		To cover all other face integrals for the first summand on the right hand side of (\ref{eq 6.23}), we denote by $\cA^+$ and $\cA^-$ the set of all numbers $i$ with $\Gamma_{1,n}^{+,(i)}$ and $\Gamma_{1,n}^{-,(i)}$ tending to a subface of an interface in $\tFint$ as $n\to\infty$, respectively.
		All numbers $i$ with $\Gamma_{1,n}^{+,(i)}$ or $\Gamma_{1,n}^{-,(i)}$ tending to an interface in $\Finteff\setminus\tFint$, are contained in the sets $\cB^+$ or $\cB^-$, respectively.
		Taking the transmission condition $\llbracket\mu\mH\cdot\nu_{\cF}\rrbracket_{\cF}=0=\llbracket\varphi\rrbracket_{\cF}$ for $\cF\in\Fint$ into account, using Lemma~\ref{lemma 6.13} and integrating by parts, we then conclude the relations
		\begin{align}
			\lim_{n\to\infty}\Big(\sum_{i\in\cA^+}\int_{\Gamma_{1,n}^{+,(i)}}(\partial_3\mH_1^{(i)})\varphi^{(i)}\dd\sigma-&\sum_{i\in\cA^-}\int_{\Gamma_{1,n}^{-,(i)}}(\partial_3\mH_1^{(i)})\varphi^{(i)}\dd\sigma\Big)\nonumber\\
			&=-\sum_{\cF\in\tFint}\int_{\cF}\llbracket\partial_3\mH_1\rrbracket_{\cF}\varphi\dd\sigma,\nonumber\\
			\lim_{n\to\infty}\Big(\sum_{i\in\cB^+}\int_{\Gamma_{1,n}^{+,(i)}}(\partial_3\mH_1^{(i)})\varphi^{(i)}\dd\sigma-&\sum_{i\in\cB^-}\int_{\Gamma_{1,n}^{-,(i)}}(\partial_3\mH_1^{(i)})\varphi^{(i)}\dd\sigma\Big)\nonumber\\
			&=\sum_{\cF\in\cF^{\text{eff}}_{\text{int},1}\setminus\tFint}\int_{\cF}\llbracket\mH_1\partial_3\varphi\rrbracket\dd\sigma=0.\label{eq 6.24}
		\end{align}		
		
		For the second summand on the right hand side of (\ref{eq 6.23}), we use the interface condition $\llbracket\tfrac{1}{\eps}\curl\mH\times\nu_{\cF}\rrbracket_{\cF}=0$ for $\cF\in\Fint$ to derive the identity
		\begin{align}
			-\lim_{n\to\infty}\sum_{i=1}^N\int_{\Gamma_{1,n}^{\pm,(i)}}(\nu_{i,n})_1(\curl\mH^{(i)})_2\varphi^{(i)}\dd\sigma=\sum_{\cF\in\cF^{\text{eff}}_{\text{int},1}}\int_{\cF}\llbracket(\curl\mH)_2\rrbracket_{\cF}\varphi\dd\sigma.\label{eq 6.25}
		\end{align}
	
		4) A straightforward modification of the arguments in part~3.ii) leads to the formula
		\begin{align}
			\lim_{n\to\infty}\sum_{i=1}^N\int_{\Gamma_{2,n}^{\pm,(i)}}(\nu_{i,n})_2(\partial_2\mH_3^{(i)})\varphi^{(i)}\dd\sigma=-\sum_{\cF\in\cF^{\text{eff}}_{\text{int},2}}\int_{\cF}\llbracket(\curl\mH)_1\rrbracket_{\cF}\varphi\dd\sigma.\label{eq 6.26}
		\end{align}
		
		Combining (\ref{eq 6.21})--(\ref{eq 6.26}), we altogether arrive at the equation 
		\begin{align}
			\int_{Q}(\nabla\mH_3)\cdot&(\nabla\varphi)\dd x=-\sum_{i=1}^N\int_{Q_i}(\Delta\mH_3^{(i)})\varphi^{(i)}\dd x-\sum_{\cF\in\tFint}\int_{\cF}\llbracket\partial_3\mH_1\rrbracket_{\cF}\varphi\dd\sigma\nonumber\\
			&\quad+\sum_{\cF\in\cF_{\text{int},1}^{\text{eff}}}\int_{\cF}\llbracket(\curl\mH)_2\rrbracket_{\cF}\varphi\dd\sigma-\sum_{\cF\in\cF_{\text{int},2}^{\text{eff}}}\int_{\cF}\llbracket(\curl\mH)_1\rrbracket_{\cF}\varphi\dd\sigma.\label{eq 6.27}
		\end{align}
		
		Since $\tilde{W}$ is dense in $\tilde{V}$, (\ref{eq 6.27}) also holds for all functions $\varphi\in \tilde{V}$. 
		Lemma~\ref{lemma 4.1} (with $\mu=1$) then shows that $\mH_3$ belongs to $PH^{2-\kappa}(Q)$ with the desired estimate.
	\end{proof}
\end{lemma}

\subsection{Conclusion of the regularity statement for the Maxwell system}\label{section 5.4}\

The next proposition summarizes the findings of Sections~\ref{section 5.1}--\ref{section 5.3}. 
It provides the desired embedding result for the space $X_2$ from (\ref{def X2}). 
Recall for the result the number $\skappa$ from (\ref{eq 6.20}).

\begin{proposition}\label{proposition 6.15}
	Let $\eps$ and $\mu$ satisfy (\ref{assumptions parameters}), and let $\kappa>1-\skappa$.
	The space $X_2$ embeds continuously into $PH^{2-\kappa}(Q)^6$.
	\begin{proof}
		Lemmas~\ref{lemma 6.11} and \ref{lemma 6.13}--\ref{lemma 6.14} show that $X_2$ is a subspace of $PH^{2-\kappa}(Q)^6$. 
		It remains to show the embedding property.
		Throughout, $C=C(\eps,\mu,\kappa,Q)>0$ is a uniform constant that is allowed to change from line to line.
		
		1) Let $(\mE,\mH)\in X_2$. 
		We first recall the estimate
		\begin{align}
			\lVert\mE\rVert_{PH^{2-\kappa}(Q)}\leq &C\Big(\lVert\mE\rVert_{L^2(Q)}+\sum_{i=1}^N\big(\lVert\Delta\mE^{(i)}\rVert_{L^2(Q_i)}+\lVert\curl\mE\rVert_{H^1(Q_i)}\big)\nonumber\\
			&+\sum_{i=1}^L\sum_{l=0}^K\lVert\divv\mE\rVert_{H^1_{00}(\tilde{Q}_{i,l})}+\sum_{\cF\in\Finteff}\lVert\llbracket\eps\mE\cdot\nu_{\cF}\rrbracket_{\cF}\rVert_{H_{00}^{3/2}(\cF)}\Big),\label{eq 6.28}
		\end{align}
		from Lemma~\ref{lemma 6.11}.
		In the following, we bound all expressions on the right hand side of (\ref{eq 6.28}) in terms of the norm in $X_2$, see (\ref{def X2}).
		By definition of $X_2$, $\curl\tfrac{1}{\mu}\curl\mE\in L^2(Q)^3$ and $\curl\mE\in H_0(\divv,Q)$.
		Hence Proposition~4.6 in \cite{Zeru22b} yields the bound
		\begin{align}
			\lVert\curl\mE\rVert_{PH^{1}(Q)}\leq C\big(\lVert\tfrac{1}{\mu}\curl\mE\rVert_{L^2(Q)}+\lVert\tfrac{1}{\eps}\curl\tfrac{1}{\mu}\curl\mE\rVert_{L^2(Q)}\big).\label{eq 6.29}
		\end{align}  
		Employing the formula $\curl\curl\mE^{(i)}=-\Delta\mE^{(i)}+\nabla\divv\mE^{(i)}$, we moreover infer the relation
		\begin{align}
			\lVert\Delta\mE^{(i)}\rVert_{L^2(Q_i)}\leq C\big(\lVert\tfrac{1}{\eps^{(i)}}\curl\tfrac{1}{\mu^{(i)}}\curl\mE^{(i)}\rVert_{L^2(Q_i)}+\lVert\divv(\eps^{(i)}\mE^{(i)})\rVert_{H^1(Q_i)}\big).\label{eq 6.30}
		\end{align}
		Taking also the embedding of $\cD(M^2)$ into $\cD(M)$ into account, (\ref{eq 6.28})--(\ref{eq 6.30}) lead us to the desired estimate $\lVert\mE\rVert_{PH^{2-\kappa}(Q)}\leq C\lVert(\mE,\mH)\rVert_{X_2}$.
		
		2) We proceed with the magnetic field component $\mH$.
		Combining Lemmas~\ref{lemma 6.13}--\ref{lemma 6.14}, we infer the relation
		\begin{align*}
			\lVert\mH\rVert_{PH^{2-\kappa}(Q)}\leq C\sum_{i=1}^N\Big(&\lVert\mH\rVert_{L^2(Q_i)}+\lVert\Delta\mH\rVert_{L^2(Q_i)}+\lVert(\curl\mH)_3\rVert_{H^1(Q_i)}\\
			&+\sum_{j=1}^2\lVert(\curl\mH)_j\rVert_{H^{1-\kappa}(Q_i)}\Big). 
		\end{align*}
		By definition of $X_2$ in (\ref{def X2}), the vector $(\tfrac{1}{\eps}\curl\mH,0)$ is an element of the space $X_1$ from (2.9) in \cite{Zeru22b}.
		Hence Proposition~4.6 in \cite{Zeru22b} and Remark~\ref{remark previous work} provide the bound
		\begin{align*}
			\lVert\curl\mH\rVert_{PH^{1-\kappa}(Q)^2\times PH^1(Q)}\leq C\Big(\lVert\tfrac{1}{\eps}\curl\mH\rVert_{L^2(Q)}+\lVert\tfrac{1}{\mu}\curl\tfrac{1}{\eps}\curl\mH\rVert_{L^2(Q)}\Big).
		\end{align*}
		Similar reasoning as in part~1) now leads to the desired inequality.
	\end{proof}
\end{proposition}

A straightforward adaption of the proof for Proposition~5.1 in \cite{Zeru22b} establishes that the part $M_2$ of $M$ on $X_2$, see (\ref{def X2}), generates a strongly continuous semigroup on $X_2$.
Hence, we skip the proof for brevity.

\begin{lemma}\label{lemma 6.16}
	Let $\eps$ and $\mu$ satisfy (\ref{assumptions parameters}).
	The part $M_2$ of $M$ generates a contractive $C_0$-semigroup $(\e^{tM_2})_{t\geq0}$ on $X_2$. 
	The latter is the restriction of $(\e^{tM})_{t\geq0}$ to $X_2$.
\end{lemma}

Using semigroup theory, we can now provide the desired wellposedness statement for the Maxwell system (\ref{Maxwell system}) on $X_2$.
The resulting regularity statement is elaborated in the subsequent Remark~\ref{remark 6.18} for further reference.
Note that the statement transfers parts of Proposition~3.3 from \cite{EiSc32} and Corollary~9.24 in \cite{ZerullaDiss} to the present setting of discontinous material parameters.
The formula for the surface charge $\rho_{\cF}$ on an interface $\cF$ is also contained in Section~1 of \cite{SchnSpi} and Corollary~5.2 of \cite{Zeru22b}.
To formulate the assumptions on the external current $\mJ$, we use the space
\begin{align*}
	W_T:=L^1([0,T],\cD(M_2))+W^{1,1}([0,T],X_2)
\end{align*}
for fixed $T>0$.
It is equipped with the norm $\lVert\cdot\rVert_{W}$ which is the canonical norm for sums of spaces.
The proof follows the same lines as the one for Corollary~5.2 in \cite{Zeru22b}, whence we skip it.
(To be more precise, one replaces $X_1$ by $X_2$, $M_1$ by $M_2$, and $V(\cF)$ by $H_{00}^{3/2}(\cF)$ in the proof.)

\begin{corollary}\label{corollary 6.17}
	Let $\eps$ and $\mu$ satisfy (\ref{assumptions parameters}).
	Let $T>0$, $w_0=(\mE_0,\mH_0)\in\cD(M_2)=\cD(M^3)\cap X_2$, and let $g:=(\tfrac{1}{\eps}\mJ,0):[0,T]\to X_2$ be continuous, and contained in $W_T$.
	The following statements are valid.
	
	a) System (\ref{Maxwell system}) has a unique classical solution $w=(\mE,\mH)$.
	It is contained in $C([0,T],\cD(M_2))\cap C^1([0,T],X_2)$ with
	\begin{align*}
		\norm{w(t)}_{X_2}&\leq\norm{w_0}_{X_2}+\norm{(\tfrac{1}{\eps}\mJ,0)}_{L^1([0,t],X_2)},\\
		\norm{Mw(t)}_{X_2}&\leq\norm{w_0}_{\cD(M_2)}+(\tfrac{2}{T}+3)\norm{g}_{W_T},\qquad t\in[0,T].
	\end{align*}
	
	b) The volume charge density $\rho^{(i)}$ on $Q_i$ and the surface charge $\rho_{\cF}$ on an interface $\cF\in\Finteff$ are given via the formulas
	\begin{align*}
		\rho^{(i)}(t)&=\divv(\eps^{(i)}\mE^{(i)}(t))=\divv(\eps^{(i)}\mE_0^{(i)})-\int_0^t\divv(\mJ^{(i)}(s))\dd s,\\
		\rho_{\cF}(t)&=\llbracket\eps\mE(t)\cdot\nu_{\cF}\rrbracket_{\cF}=\llbracket\eps\mE_0\cdot\nu_{\cF}\rrbracket_{\cF}-\int_0^t\llbracket\mJ(s)\cdot\nu_{\cF}\rrbracket_{\cF}\dd s,
	\end{align*}
	for $t\in[0,T]$, and $i\in\{1,\dots,N\}$.
\end{corollary}

\begin{remark}\label{remark 6.18}
	Combining Proposition~\ref{proposition 6.15} and \ref{corollary 6.17}, the Maxwell system (\ref{Maxwell system}) has a unique classical solution in the space $C^1([0,T],PH^{2-\kappa}(Q)^6)$ for $\kappa>1-\skappa$. (The number $\skappa$ is defined in (\ref{eq 6.20}).)
	\hfill$\lozenge$
\end{remark}


\section{Error analysis for an ADI scheme}\label{section 6}

The goal of this section is an error result for time-discrete approximations of the system (\ref{Maxwell system}) that are obtained via a Peaceman-Rachford ADI splitting scheme. 
The most important ingredient is our regularity and wellposedness analysis from Section~\ref{section 5}. 

In Section~\ref{section 6.1}, we analyze the involved splitting operators. 
As a preliminary step for the local error analysis, we then estimate a complicated term by means of a sophisticated $H^{\infty}$-functional calculus argument, see Section~\ref{section 6.2}. 
Combining all our established results, we can finally deduce the desired error estimate in Section~\ref{section 6.3}.

\subsection{Preliminary analysis of a Peaceman-Rachford ADI scheme}\label{section 6.1}\

The $\curl$ operator is splitted into the difference
\begin{align*}
	\curl=\begin{pmatrix}
		0&-\partial_3&\partial_2\\\partial_3&0&-\partial_1\\-\partial_2&\partial_1&0
	\end{pmatrix}
	=\begin{pmatrix}0&0&\partial_2\\ \partial_3&0&0\\ 0&\partial_1&0\end{pmatrix}-\begin{pmatrix}0&\partial_3&0\\ 0&0&\partial_1\\ \partial_2&0&0\end{pmatrix}
	=:\cC_1-\cC_2.
\end{align*}
The operators on the right hand side are accompanied with maximal domains
\begin{align*}
	\cD(\cC_j):=&\{u\in L^2(Q)^3\ |\ \cC_ju\in L^2(Q)^3\},\qquad j\in\{1,2\},
\end{align*}
and give rise to a splitting of the Maxwell operator $M$ into the sum $M=A+B$ with
\begin{align*}
	A:=\begin{pmatrix}
		0&\tfrac{1}{\eps}\cC_1\\\tfrac{1}{\mu}\cC_2&0
	\end{pmatrix}\quad\text{and}\quad 
	B:=\begin{pmatrix}
		0&-\tfrac{1}{\eps}\cC_2\\-\tfrac{1}{\mu}\cC_1&0
	\end{pmatrix}.
\end{align*}
The latter splitting operators are considered on the domains
\begin{align*}
	\cD(A):=&\{(\mE,\mH)\in L^2(Q)^6\ |\ (\cC_1\mH,\cC_2\mE)\in L^2(Q)^6,\ \mE_1=0\text{ on }\Gamma_2,\ \mE_2=0\text{ on }\Gamma_3,\\
	&\hphantom{\{(\mE,\mH)\in L^2(Q)^6\ |\ } \mE_3=0\text{ on }\Gamma_1\},\\
	\cD(B):=&\{(\mE,\mH)\in L^2(Q)^6\ |\ (\cC_2\mH,\cC_1\mE)\in L^2(Q)^6,\ \mE_1=0\text{ on }\Gamma_3,\ \mE_2=0\text{ on }\Gamma_1,\\
	&\hphantom{\{(\mE,\mH)\in L^2(Q)^6\ |\ } \mE_3=0\text{ on }\Gamma_2\},
\end{align*}
see \cite{ZhCZ2000,Namiki2000} and Section~2.2 in \cite{HoJaSc}.
Note that the boundary condition for the electric field is distributed onto the domains of the splitting operators.
Hereby, all traces are well-defined due to the imposed partial regularity.

We additionally point out that both splitting operators $A$ and $B$ are defined on domains on which we can apply the involved differential operators on the entire cuboid $Q$.
In other words, we implicitly impose interface conditions in $\cD(A)$ and $\cD(B)$.
The current framework then ensures that the inverses $(I-\tau A)^{-1}$ and $(I-\tau B)^{-1}$ exist.
The latter operators are needed to formulate the ADI scheme (\ref{eq 7.1}).

\begin{lemma}\label{lemma 7.1}
	Let $\eps$ and $\mu$ satisfy (\ref{assumptions parameters}). The operators $A$ and $B$ are skewadjoint on $X$.
	In particular, the operators $(I-\tau L)^{-1}$ and $\gamma_{\tau}(L):=(I+\tau L)(I-\tau L)^{-1}$ are contractive on $X$ for $L\in\{A,B\}$ and $\tau>0$.
\end{lemma}

Lemma~\ref{lemma 7.1} is a special case of Lemma~4.3 in \cite{HoJaSc}.
Note that the result is crucial for our analysis, and in particular for the unconditional stability of the considered ADI scheme.

Fix a step size $\tau>0$, and initial data $(\mE^0,\mH^0)\in\cD(B)$. 
We abbreviate the $n$th time step by $t_n:=n\tau$ for $n\in\NN$.
The solution $(\mE(t_n),\mH(t_n))$ of (\ref{Maxwell system}) with initial data $(\mE^0,\mH^0)$ is approximated by the Peaceman-Rachford ADI scheme
\begin{align}\label{eq 7.1}
	&\begin{pmatrix}\mE^n\\\mH^n
	\end{pmatrix}=\cT_{\tau}\left(\begin{pmatrix}
		\mE^{n-1}\\\mH^{n-1}
	\end{pmatrix}\right)\\
	&:=(I-\tfrac{\tau}{2}B)^{-1}(I+\tfrac{\tau}{2}A)\Big[(I-\tfrac{\tau}{2}A)^{-1}(I+\tfrac{\tau}{2}B)\begin{pmatrix}
		\mE^{n-1}\\\mH^{n-1}\end{pmatrix}-\tfrac{\tau}{2\eps}\begin{pmatrix}\mJ(t_n)+\mJ(t_{n-1})\\0\end{pmatrix}\Big].\nonumber
\end{align}
This scheme was first proposed in \cite{ZhCZ2000,Namiki2000} for the homogeneous case $\mJ=0$, and extended in \cite{EiSc18} to the inhomogeneous case $\mJ\not=0$.

The error analysis in \cite{HansenOstermann,HoJaSc,EiSc18,EiSc32,EiJaSc,KoeDiss} strongly relies on the embedding of the domain $\cD(M_2)=\cD(M^3)\cap X_2$ into $\cD(AB)$. 
The latter statement is valid if the material parameters $\eps$ and $\mu$ are sufficiently regular, meaning they belong to $W^{1,\infty}(Q)$.  
In presence of (\ref{assumptions parameters}), however, this embedding is in general not valid anymore, see Remark~10.5 in \cite{ZerullaDiss}.
Our error analysis shows that the failure of the embedding is the main reason for the order reduction.
Nevertheless, we can at least provide the following weaker result, involving the number $\skappa$ from (\ref{eq 6.20}).

\begin{lemma}\label{lemma 7.2}
	Let $\eps$ and $\mu$ satisfy (\ref{assumptions parameters}), and let $\kappa>1-\skappa$.
	The space $X_2$ embeds continuously into $PH^{2-\kappa}(Q)^6\cap\cD(A)\cap\cD(B)$.
	\begin{proof}
		Proposition~\ref{proposition 6.15} shows that $X_2$ embeds into $PH^{2-\kappa}(Q)^6$.
		The definition of $\cD(M)$ in (\ref{def Maxwell}) and the fact $X_2\sub\cD(M)$ additionally guarantee that all transmission and boundary conditions in $\cD(A)$ and $\cD(B)$ are satisfied by each element in $X_2$.
	\end{proof}
\end{lemma}

We finally introduce the operators
\begin{align}\label{eq 7.2}
	\Lambda_j(t)w:=\frac{1}{t^j(j-1)!}\int_0^t(t-s)^{j-1}\e^{sM}w\dd s,\qquad \Lambda_0(t):=\e^{t M},
\end{align}
for $w\in X, t\geq0$ and $j\in\NN$.
They are useful to expand the semigroup $(\e^{tM})_{t\geq0}$, see \cite{HansenOstermann,HoJaSc} for instance.
It is crucial that $\Lambda_j(t)$ leaves the space $X_2$ invariant, see Lemma~\ref{lemma 6.16}. 
Using standard semigroup theory, we obtain the relations
\begin{align}
	\norm{\Lambda_j(t)}_{\cL(X_2)}&\leq\frac{1}{j!},\quad t M\Lambda_{j+1}(t)=\Lambda_j(t)-\frac{1}{j!}I &&\text{on }\cD(M),\label{eq 7.3}\\
	\Lambda_0(t)=I+tM\Lambda_1(t)&=I+tM+\tfrac{1}{2}t^2M^2+t^3M^3\Lambda_3(t) &&\text{on }\cD(M^3),\label{eq 7.3.1}
\end{align}
for $t\geq0$ and $j\in\NN_0$, see Section~4 in \cite{HoJaSc}.

\subsection{Bound for a critical error term}\label{section 6.2}\

In the next two statements, we establish an estimate on a complicated term, that arises within our error analysis.
Besides an $H^{\infty}$-functional calculus approach, our regularity results from Section~\ref{section 5} come into play.

As introduced in Section~\ref{section2}, we denote the extrapolation of $L\in\{A,B\}$ to $L^2(Q)^6$ by $L_{-1}$. Moreover, we often abbreviate the electromagnetic wave $(\mE,\mH)$ by $w$, for convenience.

\begin{lemma}\label{lemma 7.3}
	Let $\eps$ and $\mu$ satisfy (\ref{assumptions parameters}), $\theta\in(0,1/4)$, $w\in PH^{1/2-2\theta}(Q)^6$, $\tau\in(0,1/4)$, and $L\in\{A,B\}$.
	The estimate
	\begin{align*}
		\lVert(I-\tfrac{\tau}{2}L_{-1})^{-1}L_{-1}w\rVert_{L^2(Q)}\leq C\tau^{-1/2-2\theta}\lVert w\rVert_{PH^{1/2-2\theta}(Q)}
	\end{align*}
	is valid with a uniform constant $C=C(\theta,\eps,\mu,Q)>0$.
	\begin{proof}
		1) We restrict ourselves to the case $L=A$, as the other one can be treated with similar arguments.
		Let $v\in\cD(A^2)$. 
		Using Lemma~\ref{lemma 7.2}, we first infer the relation
		\begin{align}
			\lvert ((I-\tfrac{\tau}{2}A_{-1})^{-1}A_{-1}w,v)\rvert=\lvert(w,A(I+\tfrac{\tau}{2}A)^{-1}v)\rvert.\label{eq 7.4}
		\end{align}
		
		To have a convenient functional calculus, we next switch to $A^2$, and collect two useful facts:
		The operator $I-A^2$ being positive definite and selfadjoint, it has well-defined positive definite selfadjoint fractional powers $(I-A^2)^{\gamma}$ for $\gamma>0$. 
		Using the skewadjointness of $A$, the relation
		\begin{align}\label{eq 7.5}
			\lVert\lambda Az\rVert^2\leq ((I-\lambda^2 A^2)z,z)=\lVert(I-\lambda^2 A^2)^{1/2}z\rVert^2,
		\end{align}
		is moreover valid for $z\in\cD(A^2)$ and $\lambda\in\RR\setminus\{0\}$.
		
		2) We next verify that the space $PH^{1/2-2\theta}(Q)^6$ embeds into $\cD((I-A^2)^{1/4-\theta})$.
		To that end, we note that the space
		\begin{align*}
			PH^2_0(Q)^6:=\{u\in L^2(Q)^6\ |\ u^{(i)}\in H^2_0(Q_i)^6=(\overline{C_c^{\infty}(Q_i)}^{\norm{\cdot}_{H^2}})^6\}
		\end{align*} 
		embeds into $\cD(I-A^2)=\cD(A^2)$.
		Furthermore, we employ the identity 
		\begin{align*}
			PH^{1/2-2\theta}(Q)^6=\big(L^2(Q)^6,PH^2_0(Q)^6\big)_{1/4-\theta,2}.
		\end{align*}
		(The latter is obtained by combining Proposition~2.11 in \cite{JerisonKenig} with Corollary~1.4.4.5 in \cite{Grisvard}, for instance.)
		Altogether, we conclude the desired embedding
		\begin{align}
			PH^{1/2-2\theta}(Q)^6\subset (L^2(Q)^6,\cD(I-A^2))_{1/4-\theta,2}=\cD((I-A^2)^{1/4-\theta}).\label{eq 7.6}
		\end{align}
	
		3) Using (\ref{eq 7.4})--(\ref{eq 7.6}) together with the Cauchy-Schwarz inequality, we now infer the estimate
		\begin{align}\label{eq 7.7}
			\lvert((I-\tfrac{\tau}{2}&A_{-1})^{-1}A_{-1}w,v)\rvert=\lvert((I-A^2)^{1/4-\theta}w,(I-A^2)^{-1/4+\theta}A(I+\tfrac{\tau}{2}A)^{-1}v)\rvert\nonumber\\
			&\leq C_{\theta}\lVert w\rVert_{PH^{1/2-2\theta}(Q)}\lVert (I-A^2)^{-1/4+\theta}\tfrac{2}{\tau}(I-(I+\tfrac{\tau}{2}A)^{-1})v\rVert,
		\end{align}
		with a uniform constant $C_{\theta}=C_{\theta}(\eps,\mu,\theta,Q)>0$.
		We next analyze the second factor on the right hand side of (\ref{eq 7.7}).
		Employing the commutativity of the resolvents of $A$ and $A^2$, we first obtain 
		\begin{align}\label{eq 7.8}
			\lVert(I-A^2)^{-1/4+\theta}\tfrac{2}{\tau}&(I-(I+\tfrac{\tau}{2}A)^{-1})v\rVert\nonumber\\
			&=\lVert A(I-\tfrac{\tau}{2}A)(I-\tfrac{\tau^2}{4}A)^{-1}(I-A^2)^{-1/4+\theta}v\rVert.
		\end{align}
		The skewadjointness of $A$ now leads to the formulas
		\begin{align}\label{eq 7.9}
			&\lVert A(I-\tfrac{\tau}{2}A)(I-\tfrac{\tau^2}{4}A^2)^{-1}(I-A^2)^{-1/4+\theta}v\rVert^2\nonumber\\
			&\quad =\big(A(I+\tfrac{\tau}{2}A)^{-1}(I-A^2)^{-1/4+\theta}v,A(I-\tfrac{\tau}{2}A)(I-\tfrac{\tau^2}{4}A^2)^{-1}(I-A^2)^{-1/4+\theta}v\big)\nonumber\\
			&\quad=\big((I-A^2)^{-1/2+2\theta}v,-A^2(I-\tfrac{\tau^2}{4}A^2)^{-1}v\big).
		\end{align}	
		Taking additionally into account that $((I-A^2)^{-1/2+2\theta}v,(I-\tfrac{\tau^2}{4}A^2)^{-1}v)$ is nonnegative, we deduce
		\begin{align}\label{eq 7.10}
			&\big((I-A^2)^{-1/2+2\theta}v,-A^2(I-\tfrac{\tau^2}{4}A^2)^{-1}v\big)\leq \big((I-A^2)^{-1/2+2\theta}v,(I-A^2)(I-\tfrac{\tau^2}{4}A^2)^{-1}v\big)\nonumber\\
			&\quad=\big((I-A^2)^{1/4+\theta}(I-\tfrac{\tau^2}{4}A^2)^{-1/2}v,(I-A^2)^{1/4+\theta}(I-\tfrac{\tau^2}{4}A^2)^{-1/2}v\big).
		\end{align}
		In view of (\ref{eq 7.7})--(\ref{eq 7.10}), we arrive at the estimate	
		\begin{align}
			\big\lvert ((I-\tfrac{\tau}{2}A_{-1})^{-1}&A_{-1}w,v)\big\rvert\nonumber\\
			 &\leq C_{\theta}\lVert w\rVert_{PH^{1/2-2\theta}(Q)}\lVert(I-A^2)^{1/4+\theta}(I-\tfrac{\tau^2}{4}A^2)^{-1/2}v\rVert.\label{eq 7.11}
		\end{align} 
		
		4) To bound the right hand side of (\ref{eq 7.11}), we employ a bounded $H^{\infty}(\Sigma_{3/4\pi})$-functional calculus for $-A^2$, see Theorem~11.5 in \cite{KuWe}. 
		(The set $\Sigma_{\phi}$ denotes the open sector of angle $\phi\in(0,\pi)$.)
		By means of the function
		\begin{align*}
			\varphi_{\tau}(z):=\frac{(1+z)^{1/4+\theta}}{(1+\tfrac{\tau^2}{4}z)^{1/2}},\qquad z\in\Sigma_{3/4\pi},
		\end{align*}
		we can estimate (\ref{eq 7.11}) via the relations
		\begin{align}
			\lVert(I-A^2)^{1/4+\theta}(I-\tfrac{\tau^2}{4}A^2)^{-1/2}v\rVert=\lVert\varphi_{\tau}(-A^2)v\rVert\leq \lVert\varphi_{\tau}\rVert_{H^{\infty}(\Sigma_{\pi/2})}\norm{v},\label{eq 7.12}
		\end{align}
		see Theorem~11.5 in \cite{KuWe}. 
		
		5) It remains to deal with the number $\lVert\varphi_{\tau}\rVert_{H^{\infty}(\Sigma_{\pi/2})}$.
		By means of the maximum modulus principle for holomorphic functions, we infer the formula
		\begin{align*}
			\lVert\varphi_{\tau}\rVert_{H^{\infty}(\Sigma_{\pi/2})}=\sup_{\Re z=0}\lvert\varphi_{\tau}(z)\rvert.	
		\end{align*}
		As a result, it suffices to estimate $\sup_{\Re z=0}\lvert\varphi_{\tau}(z)\rvert$. Let $z=\ii x\in\ii\RR$.
		We first note the identity
		\begin{align}
			\lvert\varphi_{\tau}(z)\rvert=\frac{(1+x^2)^{1/8+\theta/2}}{(1+\tfrac{\tau^4}{16}x^2)^{1/4}}.\label{eq 7.13}
		\end{align}
		The mapping on the right hand side has the critical values
		\begin{align*}
			x_1=0,\qquad x_{2/3}=\pm\sqrt{-\frac{1}{32(\tfrac{1}{64}-\tfrac{\theta}{16})}+\frac{\tfrac{1}{4}+\theta}{(\tfrac{1}{64}-\tfrac{\theta}{16})\tau^4}},
		\end{align*}
		whence the bound
		\begin{align}
			\sup_{z\in\CC_+}\lVert\varphi_{\tau}(z)\rVert=\max\{1,\lvert\varphi_{\tau}(\ii x_{2/3})\rvert\}\label{eq 7.14}
		\end{align}
		follows.		
		Hence, it remains to study the case that $\lvert\varphi_{\tau}\rvert$ attains its maximum at $\ii x_{2/3}$.
		Combining (\ref{eq 7.13}) with the inequalities 
		\begin{align*}
			\frac{\tau^4}{32(\tfrac{1}{4}-\theta)}<\frac{1}{512(\tfrac{1}{4}-\theta)},\qquad
			1<\frac{16}{\tau^4}=\frac{\tfrac{1}{4}}{\tfrac{1}{64}\tau^4}<\frac{\tfrac{3}{4}-\theta}{(\tfrac{1}{64}-\tfrac{\theta}{16})\tau^4}, 
		\end{align*}
		we conclude the relations
		\begin{align}
			\lvert\varphi_{\tau}(\ii x_{2/3})\rvert&\leq\frac{\Big(1+\tfrac{\tfrac{1}{4}+\theta}{(\tfrac{1}{64}-\tfrac{\theta}{16})\tau^4}\Big)^{1/8+\theta/2}}{\Big(\tfrac{1}{2(\tfrac{1}{4}-\theta)}-\tfrac{1}{512(\tfrac{1}{4}-\theta)}\Big)^{1/4}}
			\leq\frac{16^{1/8+\theta/2}}{(\tfrac{255}{512})^{1/4}}(\tfrac{1}{4}-\theta)^{\tfrac{1}{8}-\tfrac{\theta}{2}}\tau^{-\tfrac{1}{2}-2\theta}.	\label{eq 7.15}
		\end{align}
		
		The domain $\cD(A^2)$ being dense in $X$, (\ref{eq 7.11})--(\ref{eq 7.12}) and (\ref{eq 7.14})--(\ref{eq 7.15}) imply the assertion.
	\end{proof}
\end{lemma}

For the next statement, recall the analytic framework from Section~\ref{section 2.1}.

\begin{corollary}\label{corollary 7.4}
	Let $\eps$ and $\mu$ satisfy (\ref{assumptions parameters}), $\theta\in(1/2,1)$, $w_1\in\cD(M_2)=\cD(M^3)\cap X_2$, and $\tau\in(0,1/4]$.
	The estimate
	\begin{align*}
		\lVert(I-\tfrac{\tau}{2}A_{-1})^{-1}A_{-1}BM_2\Lambda_1(\tau)w_1\rVert_{L^2(Q)}\leq C\tau^{-\theta}\lVert w_1\rVert_{\cD(M_2)}
	\end{align*}
	is valid with a uniform constant $C=C(\theta,\eps,\mu,Q)>0$.
	\begin{proof}
		Let $\hat{w}_1:=\Lambda_1(\tau)Mw_1$.
		Proposition~\ref{proposition 6.15} and (\ref{eq 7.3}) then imply that $\hat{w}_1\in PH^{5/3}(Q)^6$ with
		\begin{align}
			\lVert\hat{w}_1\rVert_{PH^{5/3}(Q)}\leq C\lVert w_1\rVert_{\cD(M_2)},\label{eq 7.16}
		\end{align} 
		involving a constant $C=C(\eps,\mu,Q)>0$.		
		In the following, we extend the operator $B$ to an operator $\tilde{B}$ by defining
		\begin{align*}
			(\tilde{B}(\mE,\mH))^{(i)}&:=(-\tfrac{1}{\eps^{(i)}}\cC_2\mH^{(i)},-\tfrac{1}{\mu^{(i)}}\cC_1\mE^{(i)}),\qquad i\in\{1,\dots,N\},\\
			(\mE,\mH)\in\cD(\tilde{B})&:=\{(\mE,\mH)\in L^2(Q)^6\ |\ (\cC_2\mH^{(i)},\cC_1\mE^{(i)})\in L^2(Q_i)^6,\; i\in\{1,\dots,N\}\}.
		\end{align*}
		Via interpolation theory, we infer that $\tilde{B}$ is bounded from $PH^{5/3}(Q)^6$ into $PH^{2/3}(Q)^6$.
		Consequently, the vector $\check{w}_1:=B\hat{w}_1=\tilde{B}\hat{w}_1$ belongs to $PH^{2/3}(Q)^6$.
		Employing Lemma~\ref{lemma 7.3} and (\ref{eq 7.16}), we arrive at the desired relations
		\begin{align*}
			\lVert(I-\tfrac{\tau}{2}A_{-1})^{-1}A_{-1}\check{w}_1\rVert_{L^2(Q)}\leq C\tau^{-\theta}\lVert\check{w}_1\rVert_{PH^{1-\theta}(Q)}\leq C\tau^{-\theta}\lVert w_1\rVert_{\cD(M_2)},
		\end{align*} 
		with a constant $C=C(\eps,\mu,\theta,Q)>0$.
	\end{proof}
\end{corollary}

\subsection{Conclusion of the error bound}\label{section 6.3}\

We finally establish the time-discrete error bound for the Peaceman-Rachford ADI-scheme (\ref{eq 7.1}). 
Note that arguments from \cite{HansenOstermann}, Theorem~4.2 in \cite{HoJaSc}, Theorem~5.1 in \cite{EiSc18}, and Theorem~10.7 in \cite{ZerullaDiss} are used.
The major difference to \cite{HansenOstermann,HoJaSc,EiSc18} is the following:
The embedding of $\cD(M_2)$ into $\cD(AB)$ is not valid, see also Remark~10.5 in \cite{ZerullaDiss}. 
Thus we extrapolate $A$ to the operator $A_{-1}$, see Section~\ref{section2}. 
In this respect, we additionally note that we cannot estimate the expression from Corollary~\ref{corollary 7.4} without a loss of convergence order.

Throughout, the solution of (\ref{Maxwell system}) is denoted by $w=(\mE,\mH)$, and the approximation from (\ref{eq 7.1}) to $w$ at time $t_n=n\tau$ is $w_n$.
For the external current $\mJ$ in (\ref{Maxwell system}), we use the space
\begin{align*}
	W_T&:=C([0,T],\cD(M_2))\cap W^{2,1}([0,T],X_2),\\
	\norm{\cdot}_{W_T}&:=\norm{\cdot}_{C([0,T],\cD(M_2))}+\norm{\cdot}_{W^{2,1}([0,T],X_2)},
\end{align*}
for a fixed final time $T>0$. Recall also the framework from Section~\ref{section 2.2}.

\begin{theorem}\label{theorem 7.5}
	Let $\eps$ and $\mu$ satisfy (\ref{assumptions parameters}), $\theta\in(1/2,1)$, $T\geq1$ and $\tau\in(0,1/4)$.
	Let furthermore $w^0=(\mE^0,\mH^0)\in\cD(M_2)=\cD(M^3)\cap X_2$ be the initial data for (\ref{Maxwell system}), and let $(\tfrac{1}{\eps}\mJ,0)\in W_T$.
	There is a uniform constant $C=C(\eps,\mu,\theta,Q)>0$ with
	\begin{align*}
		\norm{w_n-w(t_n)}_{L^2(Q)}\leq CT\tau^{2-\theta}(\norm{w_0}_{\cD(M_2)}+\norm{(\tfrac{1}{\eps}\mJ,0)}_{W_T})
	\end{align*} 
	for all $n\in\NN_0$ with $n\tau\leq T$.
	\begin{proof}
		1) Throughout the proof, $C=C(\eps,\mu,\theta,Q)>0$ is a uniform constant that is allowed to change from line to line. 
		We first analyze the local error of (\ref{eq 7.1}). 
		Combining Corollary~\ref{corollary 6.17} with Lemma~\ref{lemma 7.2}, the function $M^k\Lambda_j(\tau)w(t)$ belongs to $\cD(A_{-1}B)\cap\cD(A)\cap\cD(B)$ for $k\in\{0,1\}$, $j\in\NN_0$, and $t\geq0$.
		Moreover, $\Lambda_3(\tau)w(t)$ is contained in $\cD(M^3)$ by definition of $X_2$.
		
		Let $k\in\NN_0$ with $(k+1)\tau\leq T$.
		We first derive a convenient representation formula for the local error.
		Equations~(5.2) and (5.3) from \cite{EiSc18} are still valid in our setting, and we recall them for reference as
		\begin{align}
			\tfrac{1}{\eps}\mJ(t_k+s)&=\tfrac{1}{\eps}\mJ(t_k)+\tfrac{s}{\eps}\mJ'(t_k)+\int_{t_k}^{t_k+s}\!\!(t_k+s-r)\tfrac{1}{\eps}\mJ''(r)\dd r,\quad s\in[0,\tau],\label{eq 7.17}\\
			w(t_{k+1})&=\Lambda_0(\tau)w(t_k)+\tau\Lambda_1(\tau)(-\tfrac{1}{\eps}\mJ(t_k),0)+\tau^2\Lambda_2(\tau)(-\tfrac{1}{\eps}\mJ'(t_k),0)\nonumber\\
			&\quad+R_k(\tau),\label{eq 7.18}
		\end{align}
		with the remainder
		\begin{align*}
			R_k(\tau)=\int_0^{\tau}\e^{(\tau-s)M}\big(\int_{t_k}^{t_k+s}(t_k+s-r)(-\tfrac{1}{\eps}\mJ''(r),0)\dd r\big)\dd s.
		\end{align*}
		
		Similar to (5.4) in \cite{EiSc18}, we obtain the formula
		\begin{align}\label{eq 7.19}
			\cT_{\tau}(w(t_k))&=(I-\tfrac{\tau}{2}B)^{-1}\Big[(I-\tfrac{\tau}{2}A_{-1})^{-1}(I+\tfrac{\tau}{2}A_{-1})(I+\tfrac{\tau}{2}B)w(t_k)\nonumber\\
			&\hphantom{=(I-\tfrac{\tau}{2}B)^{-1}\big[}+(I+\tfrac{\tau}{2}A)\big(
			\tau(-\tfrac{1}{\eps}\mJ(t_k),0)+\tfrac{\tau^2}{2}(-\tfrac{1}{\eps}\mJ'(t_k),0)\nonumber\\
			&\hphantom{=(I-\tfrac{\tau}{2}B)^{-1}\big[}+\tfrac{\tau}{2}\int_{t_k}^{t_{k+1}}(t_{k+1}-r)(-\tfrac{1}{\eps}\mJ''(r),0)\dd r\big)\Big],
		\end{align}
		where we extrapolate $A$, since $Bw(t_k)$ is in general not contained in $\cD(A)$.
		Subtracting (\ref{eq 7.18}) and (\ref{eq 7.19}), we arrive at
		\begin{align}
			&\cT_{\tau}(w(t_k))-w(t_{k+1})\nonumber\\
			&=(I-\tfrac{\tau}{2}B)^{-1}(I-\tfrac{\tau}{2}A_{-1})^{-1}\big[(I+\tfrac{\tau}{2}A_{-1})(I+\tfrac{\tau}{2}B)-(I-\tfrac{\tau}{2}A_{-1})(I-\tfrac{\tau}{2}B)\e^{\tau M}\big]w(t_k)\nonumber\\
			&\quad+(I-\tfrac{\tau}{2}B)^{-1}\big[\tau(I+\tfrac{\tau}{2}A)-\tau(I-\tfrac{\tau}{2}B)\Lambda_1(\tau)\big](-\tfrac{1}{\eps}\mJ(t_k),0)\nonumber\\
			&\quad+(I-\tfrac{\tau}{2}B)^{-1}\big[\tfrac{\tau^2}{2}(I+\tfrac{\tau}{2}A)-\tau^2(I-\tfrac{\tau}{2}B)\Lambda_2(\tau)\big](-\tfrac{1}{\eps}\mJ'(t_k),0)\nonumber\\
			&\quad+\tfrac{\tau}{2}(I-\tfrac{\tau}{2}B)^{-1}(I+\tfrac{\tau}{2}A)\int_{t_k}^{t_{k+1}}(t_{k+1}-r)(-\tfrac{1}{\eps}\mJ''(r),0)\dd r-R_k(\tau)\nonumber\\
			&=:e_{1,k}(\tau)+e_{2,k}(\tau)+e_{3,k}(\tau)+e_{4,k}(\tau)-R_k(\tau).\label{eq 7.20}
		\end{align} 
		In the sequel, we analyze the summands on the right hand side of (\ref{eq 7.20}).
		
		1.a) We rewrite the defining relation of $e_{1,k}(\tau)$ as
		\begin{align}\label{eq 7.21}
			e_{1,k}(\tau)=&(I-\tfrac{\tau}{2}B)^{-1}(I-\tfrac{\tau}{2}A_{-1})^{-1}\nonumber\\
			&\cdot[I+\tfrac{\tau}{2}M+\tfrac{\tau^2}{4}A_{-1}B-(I-\tfrac{\tau}{2}M+\tfrac{\tau^2}{4}A_{-1}B)\Lambda_0(\tau)]w(t_k).
		\end{align}
		Applying on the other hand (\ref{eq 7.3}) twice, the formulas
		\begin{align*}
			(I+\Lambda_0(\tau))w(t_k)&=(2I+\tau M\Lambda_1(\tau))w(t_k)=(2I+\tau M+\tau^2M^2\Lambda_2(\tau))w(t_k),\\
			(I-\Lambda_0(\tau))w(t_k)&=-\tau M\Lambda_1(\tau)w(t_k),
		\end{align*}
		follow, while (\ref{eq 7.3.1}) yields
		\begin{align*}
			(I-\Lambda_0(\tau))w(t_k)=(-\tau M-\tfrac{\tau^2}{2}M^2-\tau^3 M^3\Lambda_3(\tau))w(t_k).
		\end{align*} 
		Inserting these supplementary formulas into (\ref{eq 7.21}), we deduce the identities
		\begin{align*}
			e_{1,k}(\tau)&=(I-\tfrac{\tau}{2}B)^{-1}(I-\tfrac{\tau}{2}A_{-1})^{-1}\big[I-\Lambda_0(\tau)+\tfrac{\tau}{2}M(I+\Lambda_0(\tau))\\
			&\hphantom{=(I-\tfrac{\tau}{2}B)^{-1}(I-\tfrac{\tau}{2}A_{-1})^{-1}[}+\tfrac{\tau^2}{4}A_{-1}B(I-\Lambda_0(\tau))\big]w(k\tau)\\
			&=(I-\tfrac{\tau}{2}B)^{-1}(I-\tfrac{\tau}{2}A_{-1})^{-1}\big[-\tau^3 M^3\Lambda_3(\tau)+\tfrac{\tau^3}{2}M^3\Lambda_2(\tau)\\
			&\hphantom{=(I-\tfrac{\tau}{2}B)^{-1}(I-\tfrac{\tau}{2}A_{-1})^{-1}\big[}-\tfrac{\tau^3}{4}A_{-1}BM\Lambda_1(\tau)\big]w(k\tau).
		\end{align*}
		To estimate $e_{1,k}(\tau)$, we note that $\lVert M^3w(t_k)\rVert\leq\norm{Mw(k\tau)}_{X_2}$, see Section~\ref{section 2.1}.
		Combining (\ref{eq 7.3}), Lemma~\ref{lemma 7.1} and Corollary~\ref{corollary 7.4}, we derive the bound 
		\begin{align}\label{eq 7.22}
			\norm{e_{1,k}(\tau)}+\lVert(I+\tfrac{\tau}{2}B)e_{1,k}(\tau)\rVert\leq C\tau^{3-\theta}\lVert w(t_k)\rVert_{\cD(M_2)}.
		\end{align}
		
		1.b) We next deal with $e_{2,k}(\tau)$.
		Simple rearranging of operators first yields
		\begin{align*}
			(I-\tfrac{\tau}{2}A)^{-1}&(I-\tfrac{\tau}{2}B)\Lambda_1(\tau)(-\tfrac{1}{\eps}\mJ(t_k),0)\\
			&=\big[(I-\tfrac{\tau}{2}A)^{-1}(I-\tfrac{\tau}{2}M)-I+(I-\tfrac{\tau}{2}A)^{-1}\big]\Lambda_1(\tau)(-\tfrac{1}{\eps}\mJ(k\tau),0).
		\end{align*}
		Using additionally the identity $(I-\tfrac{\tau}{2}A)^{-1}(I+\tfrac{\tau}{2}A)=-I+2(I-\tfrac{\tau}{2}A)^{-1}$, we infer
		\begin{align*}
			e_{2,k}(\tau)&=\tau(I-\tfrac{\tau}{2}B)^{-1}(I-\tfrac{\tau}{2}A)\big[(I-\tfrac{\tau}{2}A)^{-1}(I+\tfrac{\tau}{2}A)\nonumber\\
			&\hphantom{=\tau(I-\tfrac{\tau}{2}B)^{-1}(I-\tfrac{\tau}{2}A)\big[}-(I-\tfrac{\tau}{2}A)^{-1}(I-\tfrac{\tau}{2}B)\Lambda_1(\tau)\big](-\tfrac{1}{\eps}\mJ(t_k),0)\nonumber\\
			&=\tau(I-\tfrac{\tau}{2}B)^{-1}(I-\tfrac{\tau}{2}A)\big[-I+\Lambda_1(\tau)+2(I-\tfrac{\tau}{2}A)^{-1}(I-\Lambda_1(\tau))\nonumber\\
			&\hphantom{=\tau(I-\tfrac{\tau}{2}B)^{-1}(I-\tfrac{\tau}{2}A)\big[}+\tfrac{\tau}{2}(I-\tfrac{\tau}{2}A)^{-1}M\Lambda_1(\tau)\big](-\tfrac{1}{\eps}\mJ(t_k),0).
		\end{align*}
		Applying now three times (\ref{eq 7.3}), the formula
		\begin{align*}
			e_{2,k}&(\tau)\\
			&=\tau(I-\tfrac{\tau}{2}B)^{-1}(I-\tfrac{\tau}{2}A)\big[\tau M\Lambda_2(\tau)-2\tau(I-\tfrac{\tau}{2}A)^{-1}M\Lambda_2(\tau)+\tfrac{\tau}{2}(I-\tfrac{\tau}{2}A)^{-1}M\\
			&\hphantom{=\tau(I-\tfrac{\tau}{2}B)^{-1}(I-\tfrac{\tau}{2}A)\big[ }+\tfrac{\tau^2}{2}(I-\tfrac{\tau}{2}A)^{-1}M^2\Lambda_2(\tau)\big](-\tfrac{1}{\eps}\mJ(t_k),0)
		\end{align*} 
		is obtained.
		Multiplying $\tau M\Lambda_2(\tau)$ by $I=(I-\tfrac{\tau}{2}A)^{-1}(I-\tfrac{\tau}{2}A)$ and using (\ref{eq 7.3}), we then conclude
		\begin{align*}
			e_{2,k}(\tau)&=\tau(I-\tfrac{\tau}{2}B)^{-1}(I-\tfrac{\tau}{2}A)\big[-\tau(I-\tfrac{\tau}{2}A)^{-1}M(\Lambda_2(\tau)-\tfrac{1}{2}I)\\
			&\hphantom{=\tau(I-\tfrac{\tau}{2}B)^{-1}(I-\tfrac{\tau}{2}A)\big[}+\tfrac{\tau^2}{2}(I-\tfrac{\tau}{2}A)^{-1}BM\Lambda_2(\tau)\big](-\tfrac{1}{\eps}\mJ(t_k),0)\\
			&=\tau(I-\tfrac{\tau}{2}B)^{-1}\big[-\tau^2M^2\Lambda_3(\tau)+\tfrac{\tau^2}{2}BM\Lambda_2(\tau)\big](-\tfrac{1}{\eps}\mJ(t_k),0).
		\end{align*}
		Lemmas~\ref{lemma 7.1}--\ref{lemma 7.2}, as well as (\ref{eq 7.3}) now give rise to
		\begin{align}
			\norm{e_{2,k}(\tau)}+\lVert(I+\tfrac{\tau}{2}B)e_{2,k}(\tau)\rVert\leq C\tau^3\lVert(\tfrac{1}{\eps}\mJ(t_k),0)\rVert_{\cD(M_2)}.\label{eq 7.23}
		\end{align}
		
		1.c) We next analyze $e_{3,k}(\tau)$, $e_{4,k}(\tau)$ and $R_k(\tau)$.
		Relation (\ref{eq 7.3}) leads to
		\begin{align*}
			e_{3,k}(\tau)=\tau^2(I-\tfrac{\tau}{2}B)^{-1}\big[\tfrac{\tau}{4}A-\tau M\Lambda_3(\tau)+\tfrac{\tau}{2}B\Lambda_2(\tau)\big](-\tfrac{1}{\eps}\mJ'(t_k),0).
		\end{align*}
		Combining Lemmas~\ref{lemma 7.1}--\ref{lemma 7.2} with (\ref{eq 7.3}), we then infer the bounds
		\begin{align}
			\lVert e_{3,k}(\tau)\rVert+\lVert(I+\tfrac{\tau}{2}B)e_{3,k}(\tau)\rVert&\leq C\tau^3\lVert(\tfrac{1}{\eps}\mJ'(t_k),0)\rVert_{X_2},\label{eq 7.24}\\
			\lVert e_{4,k}(\tau)\rVert+\lVert(I+\tfrac{\tau}{2}B)e_{4,k}(\tau)\rVert&\leq C\tau^2\lVert(\tfrac{1}{\eps}\mJ,0)\rVert_{W^{2,1}([t_k,t_{k+1}],X_2)},\nonumber\\
			\lVert R_k(\tau)\rVert+\lVert(I+\tfrac{\tau}{2}B)R_k(\tau)\rVert&\leq C\tau^2\lVert(\tfrac{1}{\eps}\mJ,0)\rVert_{W^{2,1}([t_k,t_{k+1}],X_2)}.\nonumber
		\end{align} 
		
		Altogether, we have estimated the local error $\cT_{\tau}(w(t_k))-w(t_{k+1})$ as well as the vector $(I+\tfrac{\tau}{2}B)(\cT_{\tau}(w(t_k))-w(t_{k+1}))$. 
		Estimates on the latter expression are crucial to control the error propagation.
		
		2) The global error is now estimated by means of the unconditional stability of the ADI scheme, as well as the bounds on the local error from part~1). 
		Similar to (\ref{eq 7.20}), we obtain the useful representation
		\begin{align*}
			w_n-w(t_n)&=(I-\tfrac{\tau}{2}B)^{-1}(I-\tfrac{\tau}{2}A_{-1})^{-1}\big[(I+\tfrac{\tau}{2}A_{-1})(I+\tfrac{\tau}{2}B)(w_{n-1}-w(t_{n-1}))\\
			&\qquad+((I+\tfrac{\tau}{2}A_{-1})(I+\tfrac{\tau}{2}B)-(I-\tfrac{\tau}{2}A_{-1})(I-\tfrac{\tau}{2}B)\e^{\tau M})w(t_{n-1})\big]\\
			&\quad+\sum_{l=2}^4e_{l,n-1}(\tau)-R_{n-1}(\tau).
		\end{align*}
		This recursive formula can also be transformed into
		\begin{align*}
			w_n-w(t_n)&=\sum_{k=0}^{n-1}\big[(I-\tfrac{\tau}{2}B)^{-1}(I+\tfrac{\tau}{2}A)(I-\tfrac{\tau}{2}A)^{-1}(I+\tfrac{\tau}{2}B)\big]^{n-1-k}\\
			&\hphantom{=\sum_{k=0}^{n-1}\big[}\cdot\Big(\sum_{l=1}^{4}e_{l,k}(\tau)-R_k(\tau)\Big).
		\end{align*}
		
		Similar arguments are also employed in the proof of Theorem~9.3 in \cite{Eilinghoff Diss}. 
		Combining (\ref{eq 7.22})--(\ref{eq 7.24}), Lemma~\ref{lemma 7.1} and the assumption $\tau<1/4$, we arrive at the relations
		\begin{align*}
			&\lVert w_n-w(n\tau)\rVert\leq\sum_{k=0}^{n-2}\lVert(I-\tfrac{\tau}{2}B)^{-1}\rVert\lVert(\gamma_{\tau}(A)\gamma_{\tau}(B))^{n-2-k}\gamma_{\tau}(A)\rVert\\
			&\quad\hphantom{w(n\tau)\rVert\leq\sum_{k=0}^{n-2}}\cdot\lVert(I+\tfrac{\tau}{2}B)\big(\sum_{l=1}^4e_{l,k}(\tau)-R_k(\tau)\big)\rVert+\lVert\sum_{l=1}^4e_{l,n-1}(\tau)-R_{n-1}(\tau)\rVert\\
			&\leq C\tau^3\sum_{k=0}^{n-1}\Big(\tau^{-\theta}\lVert w(t_k)\rVert_{\cD(M_2)}+\lVert(\tfrac{1}{\eps}\mJ(t_k),0)\rVert_{\cD(M_2)}\\
			&\hphantom{\leq\tau^3\sum_{k=0}^{n-1}\Big(}+\lVert(\tfrac{1}{\eps}\mJ,0)\rVert_{W^{1,\infty}([0,T],X_2)}+\tau^{-1}\lVert(\tfrac{1}{\eps}\mJ,0)\rVert_{W^{2,1}([t_k,t_{k+1}],X_2)}\Big).
		\end{align*}
		In view of Corollary~\ref{corollary 6.17}, we finally arrive at the desired error estimate.
	\end{proof}
\end{theorem}

\section{Numerical example}\label{section 7} 

In this section we illustrate Theorem~\ref{theorem 7.5} by a numerical example
implemented in MATLAB.
We consider the Maxwell equations \eqref{Maxwell system} 
with perfectly conducting boundary conditions
on the unit cube, i.e.\ with $a_i^-=0$ and $a_i^+=1$ for $i=1,2,3$.
We choose a constant magnetic permeability $ \mu\equiv 1$ and a discontinuous, piecewise constant electric permittivity
\begin{align*}
    \eps(x)=\begin{cases} 1.1 & \text{if } x_1 > 0.5 \text{ and } x_2 > 0.5
    \\
    0.1 & \text{else.}
    \end{cases}
\end{align*}
The Maxwell equations were solved on the time interval $[0,1]$ with initial data 
\begin{align*}
    \mE_{0}(x) &= \frac{1}{\eps(x)}(x_1-\tfrac{1}{2})x_1^2 (x_1-1)^2 (x_2-\tfrac{1}{2})^2 \sin(\pi x_2) \sin(\pi x_3) \begin{pmatrix}1 \\ 0 \\ 0 \end{pmatrix},
    \\[-2mm]
    \mH_0(x) &= 0,
\end{align*}
which were chosen in such a way that the boundary, transmission and divergence conditions
are fulfilled.

Our error analysis only refers to the semi-discretization in time, but for numerical computations the problem has to be discretized in space, too. 
For this purpose, we used the standard finite difference method on the Yee grid 
with 150 $\times$ 150 $\times$ 75 grid points\footnote{Repeating the experiment with a coarser space discretization (100 $\times$ 100 $\times$ 50 grid points) gave nearly the same result (data not shown).}; see \cite{Yee}, Section~3.6 in \cite{Taflove-Hagness}, or Section 4.4 in \cite{HoJaSc} for details.

Since the exact solution of this problem setting is not known, a reference solution was computed with the ADI scheme with a very fine step-size $\tau=10^{-4}$.
It would have been preferable, of course, to compute the reference solution with some other method, but the size of the problem is way too large for MATLAB's ODE solver, and computing a reference solution with the standard Yee scheme is questionable because no rigorous convergence analysis for this method in presence of discontinuous 
electric permittivity is known (to us).

\begin{figure}
     \includegraphics[angle = 0, origin = c, width = 0.7\textwidth]{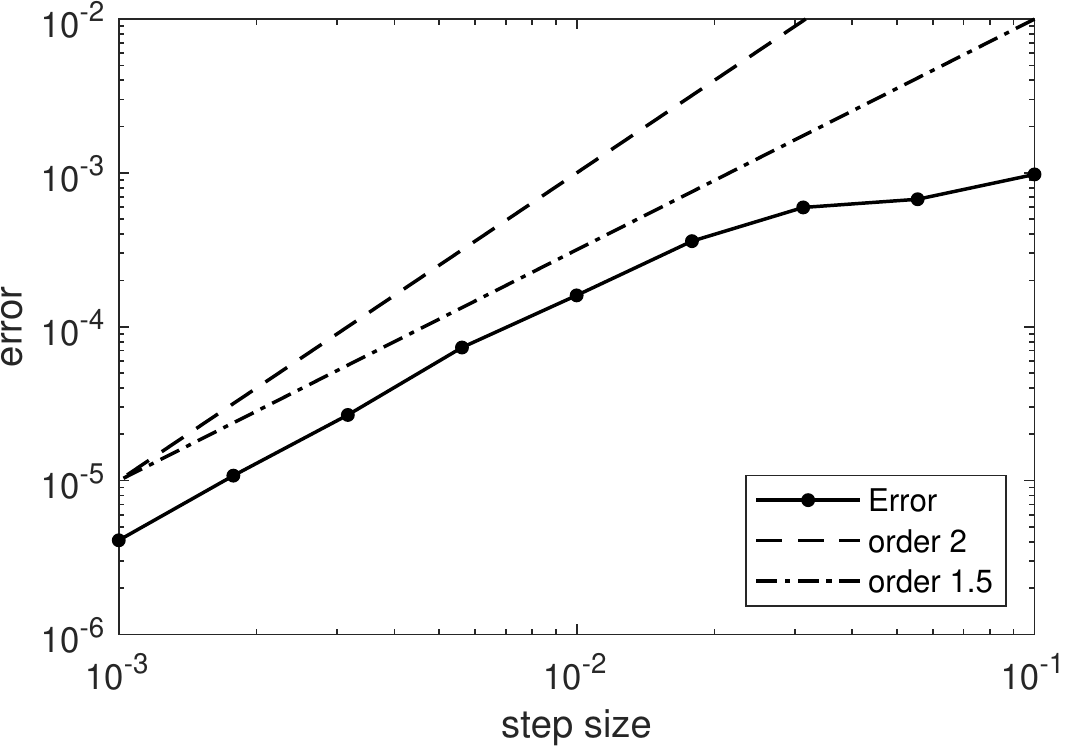}
     \caption{Error of the ADI scheme measured at $t=1$ in the discrete counterpart of $ \|\cdot\|_{L^2}.$}
\label{Fig.plot-convergence}
\end{figure}

Figure~\ref{Fig.plot-convergence} shows the error of the ADI scheme at the final time
$t=1$ for different step-sizes $\tau \in [0.001,0.1]$.
The error was measured in the discrete counterpart of $ \|\cdot\|_{L^2}, $ i.e. with integrals approximated by quadrature formulas.
The plot clearly shows that for sufficiently small $\tau$ the method converges, 
but with a reduced order of approximately 1.5 instead of the classical order 2.
This is precisely the convergence behaviour predicted by Theorem~\ref{theorem 7.5}.

\end{document}